\documentclass[12pt]{article}
\usepackage{amsmath,amsthm,amssymb,color,hyperref}

\newcommand{\ds}{\displaystyle}

\newcommand{\dfr}[2]{\dfrac{#1}{#2}}
\newcommand{\cd}{\cdot}
\newcommand{\cds}{\cdots}
\newcommand{\dsum}{\displaystyle \sum}

\newcommand{\simto}{\stackrel{\sim}{\to}}

\renewcommand{\l}{\left}
\renewcommand{\r}{\right}

\newcommand{\vsv}{\vspace{5mm}}
\newcommand{\vsb}{\vspace{2mm}}

\newcommand{\q}{\quad}

\newcommand{\maru}[1]{{\ooalign{\hfil#1\/\hfil\crcr
\raise.167ex\hbox{\mathhexbox20D}}}}

\newcommand{\ruby}[2]{%
 \leavevmode
 \setbox0=\hbox{#1}%
 \setbox1=\hbox{\tiny #2}%
 \ifdim\wd0>\wd1 \dimen0=\wd0 \else \dimen0=\wd1 \fi
 \hbox{%
   \kanjiskip=0pt plus 2fil
   \xkanjiskip=0pt plus 2fil
   \vbox{%
     \hbox to \dimen0{%
       \tiny \hfil#2\hfil}%
     \nointerlineskip
     \hbox to \dimen0{\mathstrut\hfil#1\hfil}}}}

\newcommand{\la}{\langle}
\newcommand{\ra}{\rangle}

\newcommand{\abs}[1]{\lvert{#1}\rvert}

\DeclareMathOperator*{\tensor}{\otimes}
\DeclareMathOperator*{\fusion}{\boxtimes}

\newcommand{\Z}{\mathbb{Z}}
\newcommand{\C}{\mathbb{C}}
\newcommand{\R}{\mathbb{R}}
\newcommand{\N}{\mathbb{N}}

\newcommand{\M}{\mathbb{M}}
\newcommand{\RM}{\mathrm{RM}}

\newcommand{\End}{\mathrm{End}}
\newcommand{\vir}{\mathrm{Vir}}
\newcommand{\aut}{\mathrm{Aut}}
\newcommand{\wt}{\mathrm{wt}}
\newcommand{\tr}{\mathrm{tr}}
\renewcommand{\hom}{\mathrm{Hom}}
\newcommand{\ch}{\mathrm{ch}}

\newcommand{\id}{\mathrm{id}}
\newcommand{\ind}{\mathrm{Ind}}

\newcommand{\be}{\beta}
\newcommand{\al}{\alpha}

\newcommand{\pii}{\pi \sqrt{-1}\, }

\newcommand{\Span}{\mathrm{Span}}

\newcommand{\w}{\omega}
\newcommand{\vacuum}{\mathrm{1\hspace{-3.2pt}l}}
\newcommand{\irr}{\mathrm{Irr}}

\makeatletter \@addtoreset{equation}{section}

\theoremstyle{plain}
\newtheorem{theorem}{Theorem}[section]
\newtheorem{proposition}[theorem]{Proposition}
\newtheorem{lemma}[theorem]{Lemma}
\newtheorem{corollary}[theorem]{Corollary}
\newtheorem{remark}[theorem]{Remark}

\theoremstyle{definition}
\newtheorem{definition}[theorem]{Definition}

\newtheorem{conditions}{Condition}

\title{On the structure of framed vertex operator algebras
and their pointwise frame stabilizers}

\author{
  Ching Hung Lam\footnote{Partially supported by NSC grant 94-2115-M-006-001
  of Taiwan, R.O.C.}%
  \\
  {\it \small Department of Mathematics, National Cheng Kung University}
  \\ 
  {\it \small Tainan, Taiwan 701}
  \vsb\\
  Hiroshi Yamauchi\footnote{Supported by JSPS Research
  Fellowships for Young Scientists.}
  \\
  {\it \small Graduate School of Mathematical ciences, the University of Tokyo}
  \\
  {\it \small 3-8-1 Komaba, Meguro-ku, Tokyo 153--8914, Japan}
}

\date{}

\newcommand{\sfr}[2]{\leavevmode\kern-.1em
  \raise.5ex\hbox{\the\scriptfont0 #1}\kern-.1em
  /\kern-.15em\lower.25ex\hbox{\the\scriptfont0 #2}}

\newcommand{\shf}{\sfr{1}{2}}

\DeclareMathOperator*{\PI}{\textstyle \prod}

\newcommand{\pf}{\noindent {\bf Proof:}\q }
\newcommand{\supp}{\mathrm{supp}}
\newcommand{\del}{\partial}
\newcommand{\stab}{\mathrm{Stab}}
\newcommand{\pstab}{\mathrm{Stab}^{\mathrm{pt}}}
\newcommand{\h}{\mathfrak{h}}

\begin{document}

\baselineskip 6mm

\maketitle

\vspace{-0.5cm}

\begin{center}

\vspace{0.3cm}

\noindent
{\sl Dedicated to Professor Koichiro Harada on
his 65th birthday}

\end{center}

\vspace{0.5cm}

\begin{abstract}
  In this paper, we study the structure of a general
  framed vertex operator algebra (VOA).
  We show that the structure codes $(C,D)$ of a framed VOA
  $V$ satisfy certain duality conditions.
  As a consequence, we  prove that every framed VOA
  is a simple current extension of the associated binary
  code VOA $V_C$.
  This result suggests the feasibility of classifying
  framed vertex operator algebras, at least if the
  central charge is small.
  In addition, the pointwise frame stabilizer of $V$
  is studied.
  We completely determine all automorphisms in the
  pointwise stabilizer, which are of order 1, 2 or 4.
  The 4A-twisted sector and the 4A-twisted orbifold theory
  of the famous moonshine VOA $V^\natural$ are also
  constructed explicitly.
  We verify that the top module of this twisted sector
  is of dimension 1 and of weight 3/4 and the VOA obtained
  by 4A-twisted orbifold construction of $V^\natural$ is
  isomorphic to $V^\natural$ itself.
\end{abstract}

\newpage

\tableofcontents

\newpage

\section{Introduction}

A framed vertex operator algebra $V$ is a simple vertex operator
algebra (VOA) which contains a  sub VOA $F$ called a {\it Virasoro
frame} isomorphic to a tensor product of $n$-copies of the simple
Virasoro VOA $L(\shf,0)$ such that the conformal element of $F$ is
the same as the conformal element of $V$.
There are many important examples such as the moonshine VOA
$V^\natural$ and the Leech lattice VOA. In \cite{DGH}, a basic
theory of framed VOAs was established. A general structure theory
about the automorphism group and the {\it frame stabilizer}, the
subgroup which stabilizes $F$ setwise, was also included. Moreover,
Miyamoto\,\cite{M3} showed that if $V=\oplus_{n\in \Z}V_n$ is a
framed VOA over $\R$ such that $V$ has a positive definite invariant
bilinear form and $V_1=0$, then the full automorphism group
$\aut(V)$ is finite (see also \cite{M1,M2}). Hence, the theory of
framed VOA is very useful in studying certain finite groups such as
the Monster.

It is well-known (cf.\ \cite{DMZ,DGH,M3}) that for any framed VOA $V$
with a frame $F$, one can associate two binary codes $C$ and $D$ to $V$
as follows.

Since $F\simeq L(\shf,0)^{\otimes n}$ is a rational vertex operator
algebra, $V$ is completely reducible as an $F$-module.
That is,
\begin{equation*}
  V = \bigoplus_{h_i\in\{0,1/2,1/16\}}
  m_{h_1,\ldots, h_n} L(\shf,h_1)\tensor \cds \tensor L(\shf,h_n),
\end{equation*}
where $m_{h_1,\ldots,h_r}$ is the multiplicity of the $F$-module 
$L(\shf,h_1)\tensor \cds \tensor L(\shf,h_n)$ in $V$. 
In particular, all the multiplicities are
finite and $m_{h_1,\ldots,h_r}$ is at most $1$ if all $h_i$ are
different from 1/16.

Let $M=L(\shf,h_1)\tensor \cdots \tensor L(\shf, h_n)$ be an
irreducible module over $F$.
The \emph{1/16-word (or $\tau$-word)} $\tau(M)$ of $M$ is a
binary codeword $\beta=(\be_1, \dots,\be_n)\in \Z_2^n$ such that
\begin{equation}\label{eq:1.1}
  \be_i=
  \begin{cases}
    0 & \text{ if }\ h_i=0 \text{ or } 1/2,
    \vsb\\
    1 & \text{ if }\ h_i=1/16.
   \end{cases}
\end{equation}
For any $\alpha\in \Z_2^n$, define $V^\alpha$ as the sum of all
irreducible submodules $M$ of $V$ such that $\tau(M)=\alpha$. Denote
$ D:=\{\alpha \in \Z_2^n \mid V^\alpha \ne 0\}$. Then $D$ is an even
linear subcode of $\Z_2^n$ and we obtain a $D$-graded structure on
$V=\oplus_{\alpha \in D} V^\alpha$ such that $V^\alpha\cd
V^\beta=V^{\alpha+\beta}$. In particular, $V^0$ itself is a
subalgebra and $V$ can be viewed as a $D$-graded extension of $V^0$.

For any $\gamma=(\gamma_1,\dots,\gamma_n)\in \Z_2^n$, denote
$V(\gamma):=L(\shf,h_1)\tensor \cds \tensor L(\shf,h_n)$ where
$h_i=1/2$ if $\gamma_i=1$ and $h_i=0$ elsewhere.
Set
\begin{equation*}
  C:=\{\gamma \in \Z_2^n \mid m_{\gamma_1/2,\dots,\gamma_n/2} \neq 0\}.
\end{equation*}
Then $V(0)=F$ and $V^0=\oplus_{\gamma \in C}V(\gamma )$.
The sub VOA $V^0$ forms a $C$-graded simple current extension of $F$
which has a unique simple VOA structure \cite{M2}.
A VOA of the form $V^0=\oplus_{\gamma\in C}V(\gamma)$ is often referred
to as a {\it code VOA} associated to $C$.

The codes $C$ and $D$ are very important parameters for $V$ and we
shall call them the \emph{structures codes} of $V$ with respect to
the frame $F$. One of the main purposes of this paper is to study
the precise relations between the structure codes $C$ and $D$. As
our main result, we shall show in Theorem \ref{thm:5.5} that for any
$\alpha\in D$, the subcode $C_\alpha:=\{\beta \in C \mid
\supp(\beta)\subset \supp(\alpha)\}$ contains a doubly even subcode
which is self-dual with respect to $\alpha$. From this we can prove
that every framed VOA forms a $D$-graded simple current extension of
a code VOA associated to $C$ in Theorem \ref{thm:5.6}. This shows
that one can obtain any framed VOA by performing simple current
extensions in two steps: first extend $F$ to a code VOA $V_C$
associated to $C$, then form a $D$-graded simple current extension
of $V_C$ by adjoining suitable irreducible $V_C$-modules. The
structure and representation theory of simple current extensions is
well-developed by many authors \cite{DM1,M2,L3,Y1,Y2}. It is known
that a simple current extension has a unique structure of a simple
vertex operator algebra. Since $F$ is rational, this implies there
exist only finitely many inequivalent framed VOAs with  a given
central charge. Therefore, together with the conditions on $(C,D)$
in Theorem \ref{thm:5.5}, our results provide a method for
determining all framed VOAs with a fixed central charge, at least if
the central charge is small. It is well-known that the structure
codes $(C,D)$ of a holomorphic framed VOA must satisfy $C=D^\perp$
(cf.\ \cite{DGH,M3}). In this case, we shall describe some necessary
and sufficient conditions which $C$ has to satisfy. Namely, we shall
show in Theorem \ref{thm:5.17} that there exists a holomorphic
framed VOA with structure codes $(C,C^\perp)$ if and only if $C$
satisfies the following.
\vsb\\
\begin{tabular}{lp{300pt}}
  (1) & The length of $C$ is divisible by 16.
  \vsb\\
  (2) & $C$ is even, every codeword of $C^\perp$ has a weight divisible by
  8, and \\ & $C^\perp\subset C$.
  \vsb\\
  (3) & For any $\alpha \in C^\perp$, the subcode $C_\alpha$ of $C$ contains
  a doubly even subcode\\ & which is self-dual with respect to $\alpha$.
\end{tabular}
\vsb\\ We shall call such a code an \emph{$F$-admissible code}.

Since the conditions above provide quite strong restrictions on a
code $C$, it is possible to classify all the codes satisfying these
conditions if the length is small. Once the classification of the
$F$-admissible codes of a fixed length is done, one can consider the
classification of holomorphic framed VOAs with the corresponding
central charge since a holomorphic framed VOA is always a simple
current extension of a code VOA. Based on the results of the present
paper, one can also characterize the moonshine vertex operator
algebra as the unique holomorphic framed vertex operator algebra of
central charge $24$ with trivial weight one subspace (cf.\
\cite{LY}, see also Remark \ref{rem:5.18}). It is a special case of
the famous uniqueness conjecture of Frenkel-Lepowsky-Meurman
\cite{FLM}.

In our argument,  doubly even self-dual codes play an important role
in prescribing structures of framed VOAs, and it is also revealed
that if we omit the doubly even property, then we lose the
self-duality of certain summands $V^\alpha$ of $V$ which will give
rise to an involutive symmetry analogous to the lift of the
$(-1)$-isometry on  a lattice VOA $V_L$.
By the standard notation as in \cite{FLM}, a lattice VOA has a form
\begin{equation}\label{eq:1.2}
  V_L=\bigoplus_{\alpha \in L} M_{\h}(\alpha) ,
\end{equation}
where $M_{\h}(\alpha)$ denotes the irreducible highest weight
representation over the free bosonic vertex operator algebra
$M_{\h}(0)$ associated to the vector space $\h = \C\tensor_\Z L$
with highest weight $\alpha \in \h^*=\h$. Since the fusion algebra
associated to $M_{\h}(0)$ is canonically isomorphic to the group
algebra $\C [\h]$, one has a duality relation
$M_{\h}(\alpha)^*\simeq M_{\h}(-\alpha)$. This shows that there
exists an order two symmetry inside the decomposition
\eqref{eq:1.2}, namely, we can define an involution $\theta\in
\aut(V_L)$ such that $\theta M_{\h}(\alpha)=M_{\h}(-\alpha)$ which
is an extension of an involution on $M_{\h}(0)$. However, since a
framed VOA $V$ has a decomposition $V=\oplus_{\alpha\in D}V^\alpha$
graded by an elementary abelian 2-group $D$, one cannot see the
analogous symmetry directly from the decomposition. We shall show
that by breaking the doubly even property in $(C,D)$, we can find a
pair of structure subcodes $(C^0,D^0)$ with $[C:C^0]=[D:D^0]=2$ such
that one can obtain a decomposition
\begin{equation}\label{eq:1.3}
  V=\left (\bigoplus_{\alpha \in D^0} V^{\alpha+}\oplus
  V^{\alpha-}\right)
  \bigoplus \left(\bigoplus_{\alpha \in D^1} V^{\alpha+}\oplus
  V^{\alpha-}\right)
\end{equation}
which forms a $(D^0\oplus \Z_4)$-graded simple current extension of
a code VOA associated to $C^0$, where $D^1$ is the complement of
$D^0$ in $D$.
Actually, the main motivation of the present work is to obtain the
decomposition above.
In the study of McKay's $E_8$-observation on the Monster simple group
\cite{LYY1,LYY2}, the authors found that McKay's $E_8$-observation is
related the conjectural $\Z_p$-orbifold construction of the moonshine VOA
from the Leech lattice VOA for $p>2$, where the case $p=2$ is solved in
\cite{FLM,Y3}.
Based on the decomposition \eqref{eq:1.3}, we can perform a $\Z_4$-twisted
orbifold construction on $V^\natural$.

The order four symmetry defined by the decomposition in
\eqref{eq:1.3} can be found as an automorphism fixing $F$ pointwise.
The group of automorphisms which fixes $F$ pointwise is referred to
as the \emph{pointwise frame stabilizer} of $V$. We shall show that
the pointwise frame stabilizer only has elements of order 1, 2 or 4
and it is completely determined by the structure codes $(C,D)$. As
an example, we compute the pointwise stabilizer of the Moonshine VOA
$V^\natural$ associated with a frame given in \cite{DGH,M3}. A
4A-element of the Monster is described as an element of the
pointwise frame stabilizer and the associated McKay-Thompson series
is computed in the proof of Theorem \ref{thm:7.5}. In addition, the
4A-twisted sector and the 4A-twisted orbifold theory of $V^\natural$
are constructed. We shall verify that the lowest degree subspace of
this twisted sector is of dimension 1 and of weight 3/4, and the VOA
obtained by the 4A-twisted orbifold construction of $V^\natural$ is
isomorphic to $V^\natural$ itself.

\paragraph{Acknowledgment}
The authors thank Masahiko Miyamoto for discussions and valuable
comments on the proof of Theorem \ref{thm:5.5}. They also thank
Masaaki Kitazume and Hiroki Shimakura for discussions on binary
codes. The second-named author wishes to thank Markus Rosellen for
information on the associativity and the locality of vertex
operators. Part of the work was done when the second author was
visiting the National Center for Theoretical Sciences, Taiwan on
February 2006. He thanks the staff of the center for their
hospitality.

\paragraph{Notation and Terminology}
In this article, $\N$, $\Z$ and $\C$ denote the set of non-negative
integers, integers, and the complex numbers, respectively. For
disjoint subsets $A$ and $B$ of a set $X$, we use $A\sqcup B$ to
denote the disjoint union. Every vertex operator algebra is defined
over the complex number field $\C$ unless otherwise stated. A VOA
$V$ is called {\it of CFT-type} if it has the grading
$V=\oplus_{n\geq 0}V_n$ with $V_0=\C \vacuum$. For a VOA structure
$(V,Y(\cd,z),\vacuum,\w)$ on $V$, the vector $\w$ is called the {\it
conformal vector}\footnote{We have changed the definition of the
conformal vector and the Virasoro vector. In our past works, their
definitions are opposite.} of $V$. For simplicity, we often use
$(V,\w)$ to denote the structure $(V,Y(\cd,z),\vacuum,\w)$. The
vertex operator $Y(a,z)$ of $a\in V$ is expanded as
$Y(a,z)=\sum_{n\in \Z}a_{(n)} z^{-n-1}$. For subsets $A\subset V$
and $B\subset M$ of a $V$-module $M$, we set
$$
  A\cd B :=\Span_\C \{ a_{(n)} v\mid a\in A, v\in B, n\in \Z\} .
$$
If $M$ has an $L(0)$-weight space decomposition
$M=\oplus_{n=0}^\infty M_{n+h}$ with $M_h\neq 0$, we call $M_h$ the
{\it top level} or {\it top module} of $M$ and $h$ the {\it top
weight} of $M$. The top level and top weight of a twisted module
 can be defined similarly.

For $c,h\in \C$, let $L(c,h)$  be the irreducible highest weight
module over the Virasoro algebra with central charge $c$ and highest
weight $h$. It is well-known that $L(c,0)$ has a simple VOA
structure. An element $e\in V$ is referred to as a {\it Virasoro
vector with central charge} $c_e\in \C$ if $e\in V_2$ and it
satisfies $e_{(1)}e=2e$ and $e_{(3)}e=(1/2)c_e\vacuum$. It is
well-known that by setting $L^e(n):=e_{(n+1)}$, $n\in \Z$, we obtain
a representation of the Virasoro algebra on $V$ (cf.\ \cite{M1}),
i.e.,
$$
  [L^e(m),L^e(n)]=(m-n)L^e(m+n)+\delta_{m+n,0}\dfr{m^3-m}{12}c_e.
$$
Therefore, a Virasoro vector together with the vacuum element
generates a Virasoro VOA inside $V$.
We shall denote this subalgebra by $\vir(e)$.

In this paper,  we define a sub VOA of $V$ to be a pair $(U,e)$ such
that $U$ is a subalgebra of $V$ containing the vacuum element
$\vacuum$ and $e$ is the conformal vector of $U$. Note that $(U,e)$
inherits the grading of $V$, that is, $U=\oplus_{n\geq 0} U_n$ with
$U_n=V_n\cap U$, but $e$ may not be the conformal vector of $V$. In
the case that $e$ is also the conformal vector of $V$, we shall call
the sub VOA $(U,e)$ a \emph{full} sub VOA\footnote{It is also called
a {\it conformal} sub VOA in the literature.}.

For a positive definite even lattice $L$, we shall denote the
lattice VOA associated to $L$ by $V_L$ (cf.\ \cite{FLM}). We adopt
the standard notation for $V_L$ as in \cite{FLM}. In particular,
$V_L^+$ denotes the fixed point subalgebra of $V_L$ by a lift of
$(-1)$-isometry on $L$. The letter $\Lambda$ always denotes the
Leech lattice, the unique even unimodular lattice of rank 24 without
roots.

Given an automorphism group $G$ of $V$, we denote by $V^G$ the fixed
point subalgebra of $G$ in $V$.
The subalgebra $V^G$ is called the {\it $G$-orbifold} of $V$ in
the literature.
For a $V$-module $(M,Y_M(\cd,z))$ and $\sigma\in \aut (V)$,
we set $Y_M^\sigma(a,z):= Y_M(\sigma a,z)$ for $a\in V$.
Then the {\it $\sigma$-conjugate module $M^\sigma$} of
$M$ is defined to be the module $(M,Y_M^\sigma(\cd,z))$.

We denote the ring $\Z/p\Z$ by $\Z_p$ with $p\in \Z$ and often
identify the integers $0,1,\dots, p-1$ with their images in $\Z_p$.
An additive subgroup $C$ of $\Z_2^n$ together with the standard
$\Z_2$-bilinear form is called a {\it linear code}. For a codeword
$\alpha=(\alpha_1,\dots,\alpha_n)\in C$, we define the {\it support} of
$\alpha$ by $\supp(\alpha):=\{ i \mid \alpha_i=1\}$ and the {\it weight}
by $\wt(\alpha):= \abs{\supp(\alpha)}$. For a subset $A$ of $C$, we
define $\supp(A):= \cup_{\alpha\in A} \supp(\alpha)$. For a binary
codeword $\gamma\in \Z_2^n$ and for any linear code $C\subset
\Z_2^n$, we denote $C_\gamma:=\{ \alpha \in C \mid \supp(\alpha)
\subset \supp(\gamma)\}$ and $C^{\perp_\gamma}:=\{ \beta\in C^\perp
\mid \supp(\beta)\subset \supp(\gamma)\}$, where $C^\perp=\{ \delta
\in \Z_2^n \mid \la \alpha,\delta\ra=0 ~ \text{for all}~ \alpha\in
\C\}$. A subcode $H$ of $C$ is called {\it self-dual with respect to
$\beta\in C$} if $\supp(H)=\supp(\beta)$ and $H^{\perp_\beta}=H$
(see also Definition \ref{df:4.14}). The \emph{all-one vector} is a
codeword $(11\dots1)\in \Z_2^n$.
For $\alpha=(\alpha_1,\dots,\alpha_n)$ and
$\beta\in (\beta_1,\dots,\beta_n)\in \Z_2^n$, we define
$$
  \alpha\cd \beta := (\alpha_1\beta_1,\dots,\alpha_n\beta_n) \in \Z_2^n .
$$
That is, the product $\alpha\cd \beta$ is taken in the ring $\Z_2^n$.
Note that $\alpha\cd \beta \in (\Z_2^n)_\alpha\cap (\Z_2^n)_\beta$.

\section{Preliminaries on simple current extensions}

We shall present some basic facts on simple current extensions of
a rational $C_2$-cofinite vertex operator algebra of CFT-type.

\subsection{Fusion algebras}

We recall the notion of the fusion algebra associated to a rational
VOA $V$.
It is known that a rational VOA $V$ has finitely many inequivalent
irreducible modules (cf.\ \cite{DLM2}).
Let $\irr(V)=\{ X^i \mid 1\leq i\leq r\}$ be the set of inequivalent
irreducible $V$-modules.
It is shown in \cite{HL} that the fusion product $X^i\fusion_V X^j$
exists for a rational VOA $V$.
The irreducible decomposition of $X^i\fusion_V X^j$ is referred to as
the {\it fusion rule}
\begin{equation}\label{eq:2.1}
  X^i \fusion_V X^j = \bigoplus_{k=1}^r N_{ij}^k \, X^k,
\end{equation}
where the integer $N_{ij}^k\in \Z$ denotes the multiplicity of
$X^k$ in the fusion product, and is called the {\it fusion
coefficient} which is also the dimension of the space of all
$V$-intertwining operators of type $X^i\times X^j\to X^k$.
We shall denote by $\binom{X^k}{X^i\ X^j}_V$ the space of
$V$-intertwining operators of type $X^i\times X^j\to X^k$.
We define the {\it fusion algebra} (or the {\it Verlinde algebra})
associated to $V$ by the linear space
$\mathcal{V}(V)=\oplus_{i=1}^r \C X^i$ spanned by a formal basis
$\{ X^i \mid 1\leq i\leq r\}$ equipped with a product defined by
the fusion rule \eqref{eq:2.1}.
By the symmetry of fusion coefficients, the fusion algebra $\mathcal{V}(V)$
is commutative (cf.\ \cite{FHL}).
Moreover, it is shown in \cite{H2} that if $V$ is rational, $C_2$-cofinite
and of CFT-type, then $\mathcal{V}(V)$ is associative.
In this subsection, we assume that $V$ is rational, $C_2$-cofinite and
of CFT-type.

A $V$-module $M$ is called a {\it simple current} if for any irreducible
$V$-module $X$, the fusion product $M\fusion_V X$ is also irreducible.
In other words, a simple current $V$-module $M$ induces a permutation on
$\irr(V)$ via $X \mapsto M\fusion_V X$ for $X\in \irr(V)$.
Note that $V$ itself is a simple current $V$-module.

Next we shall recall the notion of the dual module.
For a graded $V$-module $M=\oplus_{n\in \N} M_{n+h}$ such
that $\dim M_{n+h}<\infty$, define its restricted dual by
$M^*=\oplus_{n\in \N} M^*_{n+h}$, where
$M_{n+h}^*:=\hom_\C(M_{n+h},\C)$ is the dual space of $M_{h+n}$.
Let $Y_M(\cd,z)$ be the vertex operator on $M$.
We can introduce the {\it contragredient} vertex operator $Y_M^*(\cd,z)$
on $M^*$ defined by
\begin{equation}\label{eq:2.2}
  \la Y_M^*(a,z) x,v\ra
  := \la x,Y_M(e^{zL(1)}(-z^{-2})^{L(0)}a,z^{-1}) v\ra
\end{equation}
for $a\in V$, $x\in M^*$ and $v\in M$ (cf.\ \cite{FHL}).
The module $(M^*,Y_M^*(\cd,z))$ is called the {\it dual}
(or {\it contragredient}) {\it module} of $M$.

Note that if the dual module $M^*$ of $M$ is isomorphic to $N$,
there exists a $V$-isomorphism $f\in \hom_V(N,M^*)$.
Then $f$ induces a $V$-intertwining operator of type $V\times N\to M^*$.
This implies that $\binom{V^*}{M\ N}_V\ne 0$ or equivalently
$M\fusion_V N\supset V^*$.
A $V$-module $M$ is called {\it self-dual} if $M^*\simeq M$.
It is obvious that the space of $V$-invariant bilinear forms on an
irreducible self-dual $V$-module is one-dimensional.

\begin{lemma}\label{lem:2.1}
  (\cite{Y2})
  Let $U,W$ be $V$-modules such that $U\fusion_V W=V$ in the fusion algebra.
  Then both $U$ and $W$ are simple current $V$-modules.
\end{lemma}

\pf First, we show that $U\fusion_V M\ne 0$ for any $V$-module $M$.
We may assume that $M$ is irreducible as $V$ is rational. Since the
fusion product is commutative and associative, we have $(U\fusion_V
M)\fusion_V W=(U\fusion_V W)\fusion_V M=V\fusion_V M=M$. This shows
that $U\fusion_V M\ne 0$. Similarly, $W\fusion_V M\ne 0$. Now assume
that $U\fusion_V M=M^1\oplus M^2$ for $V$-submodules $M^1$ and
$M^2$. Then $W\fusion_V(U\fusion_V M)= (W\fusion_V M^1)\oplus
(W\fusion_V M^2)$. On the other hand, $$W\fusion_V(U\fusion_V
M)=(W\fusion_V U)\fusion_V M=V\fusion_V M=M. $$ Therefore,
$(W\fusion_V M^1)\oplus (W\fusion_V M^2)= M$. Since $W\fusion_V
M^i\ne 0$ if $M^i\ne 0$, we see that $U\fusion_V M$ is irreducible
if $M$ is. This shows that $U$, and also $W$, are simple current
$V$-modules. \qed

\begin{corollary}\label{cor:2.2}
  Assume that $V$ is simple, rational, $C_2$-cofinite, of CFT-type and self-dual.
  Then the following hold.
  \\
  (1) Every simple current $V$-module is irreducible.
  \\
  (2) A $V$-module $U$ is simple current if and only if $U\fusion_V U^*=V$.
  \\
  (3) The set of simple current $V$-modules forms a multiplicative abelian
  group in $\mathcal{V}(V)$ under the fusion product.
\end{corollary}

\pf
Let $U$ be a simple current $V$-module.
Then $U=V\fusion_V U$ is irreducible as $V$ is simple.
By the symmetry of fusion rules
$\binom{U}{V\ U}_V \simeq \binom{U}{U\ V}_V \simeq \binom{V^*}{U\ U^*}_V$
(cf.\ \cite{FHL}) and the assumption $V^*\simeq V$,
we have $U\fusion_V U^*\supset V$.
Since $U$ and $U^*$ are irreducible, we have $U\fusion_V U^*=V$.
This shows (1) and (2).
Now let $\mathcal{A}$ be the subset of $\mathcal{V}(V)$ consisting of
all the (inequivalent) simple current $V$-modules.
Since a fusion product of simple current modules is again a simple current,
$\mathcal{A}$ is closed under the fusion product.
Clearly $V\in \mathcal{A}$ so that $\mathcal{A}$ contains a unit element.
Finally, if $U\in \mathcal{A}$, then $U\fusion_V U^*=V$ so that the inverse
$U^*\in \mathcal{A}$ by (2).
This completes the proof.
\qed

\subsection{Simple current extensions}

We review some basic results about simple current extensions from
\cite{L3,Y1}.
Let $V^0$ be a simple rational $C_2$-cofinite VOA of CFT type and
let $\{ V^\alpha \mid \alpha \in D\}$ be a set of inequivalent
irreducible $V^0$-modules indexed by an abelian group $D$.
A simple VOA $V_D=\oplus_{\alpha\in D} V^\alpha$ is called a
{\it $D$-graded extension} of $V^0$ if $V^0$ is a full sub VOA of
$V_D$ and $V_D$ carries a $D$-grading, i.e.,
$V^\alpha \cd V^\beta \subset V^{\alpha+\beta}$ for $\alpha,\beta\in D$.
In this case, the dual group $D^*$
of $D$ acts naturally and faithfully on $V_D$. If all $V^\alpha$,
$\alpha\in D$, are simple current $V^0$-modules, then $V_D$ is
referred to as a {\it $D$-graded simple current extension} of $V^0$.
The abelian group $D$ is automatically finite since $V^0$ is
rational (cf.\ \cite{DLM2}).

\begin{proposition}\label{prop:2.3}
  (\cite{ABD,DM2,L3,Y1})
  Let $V^0$ be a simple rational $C_2$-cofinite VOA of CFT type.
  Let $V_D=\oplus_{\alpha\in D} V^\alpha$ be a $D$-graded simple current
  extension of $V^0$.
  Then
  \\
  (1) $V_D$ is rational and $C_2$-cofinite.
  \\
  (2) If $\tilde{V}_D=\oplus_{\alpha\in D}\tilde{V}^\alpha$ is another
  $D$-graded simple current extension of $V^0$ such that
  $\tilde{V}^\alpha\simeq V^\alpha$ as $V^0$-modules for all $\al\in D$,
  then $V_D$ and $\tilde{V}_D$ are isomorphic VOAs over $\C$.
   \\
  (3) For any subgroup $E$ of $D$, a subalgebra
  $V_E:=\oplus_{\alpha \in E} V^\alpha$ is an $E$-graded simple current
  extension of $V^0$.
  Moreover, $V_D$ is a $D/E$-graded simple current extension of $V_E$.
\end{proposition}

A representation theory of simple current extensions is developed in
\cite{L3,Y1}.
It is shown that each irreducible module over a simple current extension
corresponds to an irreducible module over a finite dimensional semisimple
associative algebra.
Moreover, it is also proved that any $V^0$-module can be extended
to certain twisted modules over $V_D$.

Let $M$ be an irreducible $V_D$-module.
Since $V^0$ is rational, we can take an irreducible
$V^0$-submodule $W$ of $M$.
Define
$D_W:=\{\alpha \in D \mid V^\alpha \fusion_{V^0} W\simeq_{V^0} W\}$.
Then $D_W$ is a subgroup of $D$.
Note that the subgroup $D_W$ is independent of the choice
of the irreducible $V^0$-module $W$ in $M$.
In other words, $D_W=D_{W'}$ for any irreducible
$V^0$-submodules $W$ and $W'$ of $M$.
We call $M$ \emph{ $D$-stable} if $D_W=0$.
In this case,
$V^\alpha \fusion_{V^0} W\simeq_{V^0} V^\beta \fusion_{V^0} W$
if and only if $\alpha =\beta$ and by setting
$M^\alpha:=V^\alpha \fusion_{V^0} W$, we have a
$D$-graded isotypical decomposition
$M=\oplus_{\alpha\in D} M^\alpha (\simeq V \fusion_{V^0} W)$
as a $V^0$-module.

\begin{theorem}\label{thm:2.4}
  (\cite{L3,Y1})
  Let $W$ be an irreducible $V^0$-module.
  Then there exists a unique $\chi_W \in D^*\subset \aut(V_D)$
  such that $W$ can be extended to an irreducible
  $\chi_W$-twisted $V_D$-module.
  If $D_W=0$, then the extension of $W$ to an irreducible
  $\chi_W$-twisted $V_D$-module is unique and $D$-stable.
  Moreover, the extension of $W$ is given by
  $V_D\fusion_{V^0}W$ as a $V^0$-module.
\end{theorem}

One can easily compute fusion rules among
irreducible $D$-stable modules.

\begin{proposition}\label{prop:2.5}
  (\cite{SY,Y1})
  Let $V_D$ be a $D$-graded simple current extension of a simple
  rational $C_2$-cofinite VOA $V^0$ of CFT-type.
  Let $M^i$, $i=1,2,3$ be irreducible $D$-stable $V_D$-modules.
  Denote by $M^i=\oplus_{\alpha \in D} (M^i)^{\alpha}$ a $D$-graded
  isotypical decomposition of $M^i$.
  Then the following linear isomorphism holds:
  $$
    \binom{M^3}{M^1\q M^2}_{V_D}\simeq
    \binom{(M^3)^\gamma}{(M^1)^\alpha \q (M^2)^\beta}_{V^0},
  $$
  where $\alpha,\beta,\gamma\in D$ are arbitrary.
\end{proposition}

We shall need the following result on $\Z_2$-graded simple current
extensions.

\begin{proposition}\label{prop:2.6}
  Let $V^0$ be a simple rational $C_2$-cofinite self-dual VOA of CFT-type.
  Let $V^1$ be a simple current $V^0$-module not isomorphic to $V^0$ such that
  $V^1\fusion_{V^0} V^1 =V^0$.
  Assume that $V^1$ has  an integral top weight and the invariant bilinear form
  on $V^1$ is symmetric.
  Then there exists a unique simple VOA structure on $V=V^0\oplus V^1$ as
  a $\Z_2$-graded simple current extension of $V^0$.
\end{proposition}

\pf For $a,b\in V^0$ and $u,v\in V^1$, define a vertex operator
$Y(\cd,z)$ as follows:
$$
\begin{array}{l}
  Y(a,z)b:=Y_{V^0}(a,z)b,\ Y(a,z)u:=Y_{V^1}(a,z)u,
  \ Y(u,z)a:=e^{zL(-1)}Y(a,-z)u,
\end{array}
$$
and $Y(u,z)v$ is defined by means of the matrix coefficients
$$
   \la Y(u,z)v,a\ra_{V^0}=\la v,Y(e^{zL(1)}(-z^{-2})^{L(0)}u,z^{-1})a\ra_{V^1},
$$
where $\la\cd,\cd\ra_{V^i}$ denotes the invariant bilinear form on
$V^i$, $i=0,1$. Since the invariant bilinear form on $V^1$ is
symmetric, we have the skew-symmetry $Y(u,z)v=e^{zL(-1)}Y(v,-z)u$
for any $u,v\in V^1$ by Proposition 5.6.1 of \cite{FHL}. It is also
shown in \cite{FHL,Li2} that $(V^0\oplus V^1,Y(\cd,z))$ forms a
$\Z_2$-graded simple vertex operator algebra if and only if we have
a locality for any three elements in $V^1$, that is, for any $u,v\in
V^1$, there exists $N\in \N$ such that for any $w\in V^1$ we have
$$
  (z_1-z_2)^N Y(u,z_1)Y(v,z_2)w = (z_1-z_2)^N Y(v,z_2)Y(u,z_1)w .
$$
By Huang \cite[Theorem 3.5]{H2} (see also Theorem 3.2 and 3.5 of
\cite{H}), it is shown that there exists $\lambda \in \C^*$ such
that for any $u,v,w\in V^1$ and sufficiently large $k\in \Z$,  we
have
\begin{equation}\label{eq:2.3}
  (z_0+z_2)^k Y(u,z_0+z_2)Y(v,z_2)w = \lambda (z_2+z_0)^k Y(Y(u,z_0)v,z_2)w .
\end{equation}
We shall show that the associativity above leads to the locality.
The idea of the following argument comes from \cite{R}. Let $N\in
\Z$ such that $z^N Y(u,z)v \in V^0[\![z]\!]$. Take sufficiently
large $s,t\in \Z$ such that $z^s Y(v,z)w \in V^0[\![ z]\!]$ and
\eqref{eq:2.3} holds for $(u,v,w)$ and $(v,u,w)$ with $k=t,s$. Then
$$
\begin{array}{l}
  z_1^t z_2^s (z_1-z_2)^N Y(u,z_1)Y(v,z_2)w
  \vsb\\
  = e^{-z_2\del_{z_1}}\l( (z_1+z_2)^t z_2^s  z_1^N Y(u,z_1+z_2)Y(v,z_2)w\r)
  \vsb\\
  = \lambda e^{-z_2\del_{z_1}}\l( (z_2+z_1)^t z_2^s  z_1^N Y(Y(u,z_1)v,z_2)w\r)
  \vsb\\
  = \lambda e^{-z_2\del_{z_1}}\l( (z_2+z_1)^t z_2^s  z_1^N
     Y(e^{z_1L(-1)}Y(v,-z_1)u,z_2)w\r)
  \vsb\\
  = \lambda e^{-z_2\del_{z_1}} e^{z_1\del_{z_2}} \l( z_2^t (z_2-z_1)^s z_1^N
     Y(Y(v,-z_1)u,z_2)w\r) .
\end{array}
$$
Define $p(z_1,z_2):=z_2^t (z_2-z_1)^s  z_1^N Y(Y(v,-z_1)u,z_2)w$.
The equations above show that $p(z,w)\in V^1[\![ z_1,z_2]\!]$.
On the other hand,
$$
\begin{array}{l}
  z_1^t z_2^s (-z_2+z_1)^N Y(v,z_2) Y(u,z_1) w
  \vsb\\
  = e^{-z_1\del_{z_2}}\l( z_1^t (z_2+z_1)^s (-z_2)^N Y(v,z_2+z_1)Y(u,z_1)w \r)
  \vsb\\
  = \lambda e^{-z_1\del_{z_2}}\l( z_1^t (z_1+z_2)^s (-z_2)^N Y(Y(v,z_2)u,z_1)w \r)
  \vsb\\
  = \lambda e^{-z_1\del_{z_2}} p(-z_2,z_1).
\end{array}
$$
Thus the locality follows from
$$e^{-w\del_z} e^{z\del_w}p(z,w)=e^{-w\del_z} p(z,w+z) = e^{-w\del_z} p(z,z+w)
= p(z-w,z)$$ and $e^{-z\del_w} p(-w,z)=p(-w+z,z)=p(z-w,z)$.

The uniqueness has already been shown in \cite{DM2} in a general
fashion.
\qed
\vsb

Later, we shall consider a construction of framed VOAs. The
following extension property will be used frequently.

\begin{theorem}\label{thm:2.7}
  (Extension property \cite[Theorem 4.6.1]{Y2})
  Let $V^{(0,0)}$ be a simple rational $C_2$-cofinite VOA of CFT-type,
  and let $D_1$, $D_2$ be finite abelian groups.
  Assume that we have a set of inequivalent irreducible simple current
  $V^{(0,0)}$-modules $\{ V^{(\alpha,\beta)} \mid (\alpha,\beta) \in
  D_1\oplus D_2\}$ with $D_1\oplus D_2$-graded  fusion rules
  $V^{(\alpha_1,\beta_1)} \fusion_{V^{(0,0)}} V^{(\alpha_2,\beta_2)}=
  V^{(\alpha_1 +\alpha_2,\beta_1+\beta_2)}$ for any
  $(\alpha_1,\beta_1), (\alpha_2,\beta_2) \in D_1\oplus D_2$.
  Further assume that all $V^{(\alpha,\beta)}$, $(\alpha,\beta) \in D_1\oplus
  D_2$, have integral top weights and we have $D_1$- and $D_2$-graded simple
  current extensions $V_{D_1}=\oplus_{\alpha \in D_1} V^{(\alpha,0)}$ and
  $V_{D_2}= \oplus_{\beta \in D_2} V^{(0,\beta)}$.
  Then $V_{D_1\oplus D_2}:=\oplus_{(\alpha,\beta) \in D_1\oplus D_2}
  V^{(\alpha,\beta)}$ possesses a unique structure of a simple vertex operator
  algebra as a $(D_1\oplus D_2)$-graded simple current extension of $V^{(0,0)}$.
\end{theorem}

\section{Ising frame and framed VOA}

We shall review the notion of an Ising frame and a framed vertex
operator algebra.

\subsection{Miyamoto involutions}

We begin by the definition of an Ising vector.

\begin{definition}\label{df:3.1}
  A Virasoro vector $e$ is called an {\it Ising vector} if $\vir(e)\simeq L(\shf,0)$.
  Two Virasoro vectors $u,v\in V$ are called {\it orthogonal} if
  $[Y(u,z_1), Y(v,z_2)]=0$.
  A decomposition $\w=e^1+\cds +e^n$ of the conformal vector $\w$ of $V$ is
  called {\it orthogonal} if  $e^i$ are mutually orthogonal Virasoro vectors.
\end{definition}

Let $e\in V$ be an Ising vector. By definition, $\vir(e)\simeq
L(\shf,0)$. It is well-known that $L(\shf,0)$ is rational,
$C_2$-cofinite and has three irreducible modules $L(\shf,0)$,
$L(\shf,\shf)$ and $L(\shf,\sfr{1}{16})$. The fusion rules of
$L(\shf,0)$-modules are computed in \cite{DMZ}:
\begin{equation}\label{eq:3.1}
\begin{array}{l}
  L(\shf,\shf)\fusion L(\shf,\shf)=L(\shf,0),
  \q
  L(\shf,\shf)\fusion L(\shf,\sfr{1}{16})=L(\shf,\sfr{1}{16}),
  \vsb\\
  L(\shf,\sfr{1}{16})\fusion L(\shf,\sfr{1}{16})=L(\shf,0)\oplus L(\shf,\shf).
\end{array}
\end{equation}
By \eqref{eq:3.1}, one can define some involutions in the following way.
Let $V_e(h)$ be the sum of all irreducible $\vir(e)$-submodules of $V$ isomorphic
to $L(\shf,h)$ for $h=0,1/2,1/16$.
Then one has the isotypical decomposition
$$
  V=V_e(0)\oplus V_e(\shf)\oplus V_e(\sfr{1}{16}).
$$
Define a linear automorphism $\tau_e$ on $V$ by
$$
  \tau_e=
  \begin{cases}
    \ \  1 & \text{ on }\  V_e(0)\oplus V_e(\shf),
    \vsb\\
    -1& \text{ on }\  V_e(\sfr{1}{16}).
\end{cases}
$$
Then by the fusion rules in \eqref{eq:3.1}, $\tau_e$ defines an automorphism on
the VOA $V$ (cf.\ \cite{M1}).
On the fixed point subalgebra $V^{\la \tau_e\ra}=V_e(0)\oplus V_e(\shf)$,
one can define another linear automorphism $\sigma_e$ by
$$
  \sigma_e=
  \begin{cases}
    \ \  1 & \text{ on }\  V_e(0),
    \vsb\\
    -1& \text{ on }\  V_e(\sfr{1}{2}).
  \end{cases}
$$
Then $\sigma_e$ also defines an automorphism on $V^{\la \tau_e\ra}$ (cf.\ \cite{M1}).
The automorphisms $\tau_e\in \aut(V)$ and $\sigma_e\in \aut(V^{\la \tau_e\ra})$ are
often called {\it Miyamoto involutions}.

\subsection{Framed VOAs and their structure codes}\label{sec:3.2}

Let us define the notion of a framed VOA.

\begin{definition}\label{df:3.2}
  (\cite{DGH,M3})
  A simple vertex operator algebra $(V,\w)$ is called \emph{framed}
  if there exists a set $\{e^1, \dots,e^n\}$ of Ising vectors of $V$ such that
  $\w=e^1+\cds +e^n$ is an orthogonal decomposition.
  The full sub VOA  $F$ generated by $e^1,\dots,e^n$ is called
  an \emph{Ising frame} or simply  a \emph{frame} of $V$.
  By abuse of notation, we sometimes call the set of Ising vectors
  $\{ e^1,\dots, e^n\}$ a {\it frame}, also.
\end{definition}

Let $(V,\w)$ be a framed VOA with an Ising frame $F$.
Then
$$
  F \simeq \vir(e^1)\tensor \cds \tensor \vir(e^i)
  \simeq  L(\shf,0)^{\tensor n}
$$
and $V$ is a direct sum of irreducible $F$-submodules
$\tensor_{i=1}^n L(\shf,h_i)$ with $h_i\in \{ 0,1/2,1/16\}$.
For each irreducible $F$-module $W=\tensor_{i=1}^n L(\shf,h_i)$,
we define its binary {\it $1/16$-word} (or {\it $\tau$-word})
$\tau(W)=(\alpha_1,\cds,\alpha_n)\in \Z_2^n$ by $\alpha_i=1$
if and only if $h_i=1/16$.
For $\alpha \in \Z_2^n$, denote by $V^\alpha$ the sum of
all irreducible $F$-submodules of $V$ whose $1/16$-words
are equal to $\alpha$.
Define a linear code
$D\subset \Z_2^n$ by $D=\{ \alpha \in \Z_2^n \mid V^\alpha\ne 0\}$.
Then we have the {\it 1/16-word decomposition}
$V=\oplus_{\alpha \in D} V^\alpha$.
By the fusion rules of $L(\shf,0)$-modules,
it is easy to see that
$V^\alpha\cd V^\beta \subset V^{\alpha +\beta}$.
Hence, the dual group $D^*$ of $D$ acts on $V$.
In fact, the action of $D^*$ coincides with
the action of the elementary abelian 2-group generated
by Miyamoto involutions $\{ \tau_{e^i} \mid 1\leq i \leq n\}$.
Therefore, all $V^\alpha$, $\alpha \in D$, are
irreducible $V^0$-modules by \cite{DM1}.
Since there is no $L(\shf,\sfr{1}{16})$-component in $V^0$,
the fixed point subalgebra $V^{D^*}=V^0$ has
the following shape:
$$
  V^0= \bigoplus_{h_i\in \{ 0,1/2\}} m_{h_1,\dots,h_n}
       L(\shf,h_1)\tensor \cds \tensor L(\shf,h_n),
$$
where $m_{h_1,\dots,h_n}\in \N$ denotes the multiplicity.
On $V^0$ we can define Miyamoto involutions $\sigma_{e^i}$
for $i=1,\dots,n$.
Denote by $Q$ the elementary abelian 2-subgroup of
$\aut (V^0)$ generated by $\{ \sigma_{e^i} \mid 1\leq i\leq n\}$.
Then the fixed point subalgebra
$(V^0)^Q=F$ and each
$m_{h_1,\cds,h_n} L(\shf,h_1)\tensor \cds \tensor L(\shf,h_n)$
is an irreducible $F$-submodule again by \cite{DM1}.
Thus $m_{h_1,\cds,h_n}\in \{ 0,1\}$ and we obtain an even linear code
$C:= \{ (2h_1,\cds ,2h_n)\in \Z_2^n \mid h_i\in
\{ 0,\shf\},\ m_{h_1,\cds,h_n}\ne 0\}$,
namely,
\begin{equation}\label{eq:3.2}
  V^0=\bigoplus_{\alpha =(\alpha_1,\cds,\alpha_n)\in C}
  L(\shf,\alpha_1/2) \tensor \cds \tensor L(\shf,\alpha_n/2).
\end{equation}
Since $L(\shf,0)$ and $L(\shf,\shf)$ are simple current
$L(\shf,0)$-modules, $V^0$ is a $C$-graded simple current extension
of $F$. By Proposition \ref{prop:2.3}, the simple VOA structure on
$V^0$ is unique. The simple VOA $V^0$ of the form \eqref{eq:3.2} is
called the {\it code VOA associated to $C$} and denoted by $V_C$. It
is clear that $V_C$ is  simple, rational, $C_2$-cofinite and of
CFT-type. Since $L(1)(V_C)_1=0$ by its $F$-module structure, $V_C$
has a non-zero invariant form and thus is self-dual as a
$V_C$-module by \cite{Li1}. Similarly, we also have $L(1)V_1=0$ and
$V$ is self-dual as a $V$-module.

Summarizing, there exists a pair $(C,D)$ of even linear codes
such that $V$ is an $D$-graded extension of a code VOA $V_C$
associated to $C$.
We call the  pair $(C,D)$ the {\it structure codes} of
a framed VOA $V$ associated with the frame $F$.
Since the powers of $z$ in an $L(\shf,0)$-intertwining operator
of type $L(\shf,\shf)\times L(\shf,\shf)\to L(\shf,\sfr{1}{16})$
are half-integral, the structure codes $(C, D)$ satisfy
$C\subset D^\perp$.

\paragraph{Notation}
Let $V$ be a framed VOA with the structure codes $(C,D)$,
where $C,D\subset \Z_2^n$.
For a binary codeword $\beta \in \Z_2^n$, we define:
\begin{equation}\label{eq:3.3}
  \ds \sigma_\beta:=\PI_{i\in \supp(\beta)}\sigma_{e^i}\in \aut(V^0)
  \q \mbox{ and }\q \tau_\beta:=\PI_{i\in \supp(\beta)}\tau_{e^i} \in \aut(V).
\end{equation}
Namely, by associating Miyamoto involutions to a codeword of $\Z_2^n$,
$\sigma : \Z_2^n\to \aut(V^0)$ and $\tau : \Z_2^n\to \aut(V)$
define group homomorphisms.
It is also clear that $\ker \sigma=C^\perp$ and $\ker \tau=D^\perp$.

\section{Representation of code VOAs}

Since every framed VOA is an extension of its code sub VOA, it is
quite natural to study a framed VOA as a module over its code sub VOA.
Let us first review a structure theory for the irreducible modules
over a code VOA.

\subsection{Central extension of codes}

Let $\nu^1=(10\dots0)$, $\nu^2=(010\dots0),\dots,\nu^n=(0\dots01)\in
\Z_2^n$. Define $\varepsilon :\Z_2^n\times \Z_2^n \to \C^*$ by
\begin{equation}\label{eq:4.1}
  \varepsilon (\nu^i,\nu^j):=-1\q \mbox{ if }\q i> j
  \q \mbox{ and }\q 1\q \mbox{ otherwise},
\end{equation}
and extend to $\Z_2^n$ linearly. Then $\varepsilon$ defines a
2-cocycle in $Z^2(\Z_2^n,\C^*)$. By definition,
\begin{equation}\label{eq:4.2}
  \varepsilon(\alpha,\beta)\varepsilon(\beta,\alpha)
  = (-1)^{\la \alpha,\beta\ra +\wt(\alpha)\wt(\beta)}
    \ \mbox{and} \
    \varepsilon(\alpha,\alpha)
    = (-1)^{\wt(\alpha) (\wt(\alpha)-1)/2}
\end{equation}
for all $\alpha,\beta\in \Z_2^n$. In particular,
$\varepsilon(\alpha,\alpha)=(-1)^{\wt(\alpha)/2}$  and
$\varepsilon(\alpha,\beta)\varepsilon(\beta,\alpha)
  = (-1)^{\la \alpha,\beta\ra }$ if
$\alpha, \beta \in \Z_2^n$ are  even.

Let $G$ be the central extension of $\Z_2^n$ by $\C^*$ with
associated 2-cocycle $\varepsilon$.
Recall that $G=\Z_2^n \times \C^*$ as a set, but the group
operation is given by
$$
  (\alpha,u)(\beta,v)=(\alpha\beta, \varepsilon(\alpha,\beta)uv)
$$
for all $\alpha,\beta \in \Z_2^n$ and $u,v\in \C^*$.
Let $C$ be a binary even linear code of $\Z_2^n$.
Since $\varepsilon$ takes values in $\{ \pm 1\}$, we can
take a subgroup
$\tilde{C}=\{ (\alpha, u)\in G \mid \alpha \in C,\, u\in \{ \pm 1\}\}$
of $G$ so that $\tilde{C}$ forms  a central extension of
$C$ by $\{ \pm 1\}$:
\begin{equation}\label{eq:4.3}
  1\longrightarrow \{ \pm 1\} \longrightarrow \tilde{C}
  \stackrel{\pi}{\longrightarrow} C \longrightarrow 1.
\end{equation}
We shall note that the radical of the standard bilinear form $\la\cd
,\cd \ra$ on $C$ is given by $R=C\cap C^\perp$ and thus by
\eqref{eq:4.2}, the preimage $\tilde{R}=\pi^{-1}(C\cap C^\perp)$ is
the center of $\tilde{C}$. Take a subgroup $D$ of $C$ such that $
C=R\oplus D$. Then the form $\la\cd,\cd\ra$ is non-degenerate on
$D$.
It follows from \eqref{eq:4.2} that the preimage
$\tilde{D}:=\pi^{-1}(D)$ is an extra-special 2-subgroup of
$\tilde{C}$.
The central extension $\tilde{C}$ is then isomorphic\footnote{Note
that the isomorphism type of $\tilde{D}$ is determined by the
dimension of maximal isotropic subspaces of $D$ with respect to the
quadratic form $q(\al)=\varepsilon(\al,\al)$ (cf. \cite{Go} and
\cite[Section 5.3]{FLM}), which depends on the choice of the
complement $D$ if $R$ is not doubly even. For example, we can take
$C=\Span_{\Z_2}\{ (11000),(00110),(00101)\}$. Then the radical
$R=\{(00000), (11000)\}$. Set $D=\{(00000),
(00110),(00101),(00011)\}$ and
$D'=\{(00000),(11110),(00101),(11011)\}$. Then both $D$ and $D'$ are
complement of $R$ in $C$ but $\tilde{D}\not\simeq \tilde{D}'$, for
$\tilde{D}$ is a quaternion group whereas $\tilde{D}'$ is a dihedral
group of order 8. Nevertheless, the central product $\tilde{D}
*_{\{\pm 1\}} \tilde{R}$ is still uniquely determined by $C$ up to
isomorphisms.} to the central product of $\tilde{D}$ and $\tilde{R}$
over $\{\pm 1\}\subset \C^*$ which we shall denote by $\tilde{D}
*_{\{\pm 1\}} \tilde{R}$.

We identify the multiplicative group $\C^*$ with the central
subgroup $(0,\C^*)=\{ (0,u ) \in G \mid u\in \C^*\}$ of $G$, and let
$\C^*\tilde{C}=\{ (\alpha,u) \in G \mid \alpha \in C, u\in \C^*\}$
be the subgroup of $G$ generated by $\C^*=(0,\C^*)$ and $\tilde{C}$.
Then we have the exact sequence:
\begin{equation}\label{eq:4.4}
  1\longrightarrow \C^* \longrightarrow \C^*\tilde{C}
  \stackrel{\pi_{\C^*}}{\longrightarrow} C \longrightarrow 1.
\end{equation}
Since $\C^*$ is injective in the category of abelian groups,
the preimage of $C\cap C^\perp$ in $\C^*\tilde{C}$ splits and one has
an isomorphism
$$
  \C^*\tilde{C}
  \simeq (C\cap C^\perp) \times (\C^* *_{\{\pm 1\}}\tilde{D}) .
$$

Now let $\psi : C\to \End (V)$ be a $\varepsilon$-projective
representation of $C$ on $V$, that is,
$\psi(\alpha)\psi(\beta)=\varepsilon(\alpha,\beta)\psi(\alpha+\beta)$
for $\alpha,\beta\in C$. Then one defines a linear representation
$\tilde{\psi}$ of $\C^*\tilde{C}$ via
$\tilde{\psi}(\alpha,u):=u\psi(\alpha)\in \End (V)$ for $\alpha \in
C$ and $u\in \C^*$. Since $\C^*\tilde{C}$ is isomorphic to a direct
product of $R=C\cap C^\perp$ and $\C^* *_{\{\pm 1\}} \tilde{D}$,
$\tilde{\psi}$ is a tensor product of a linear character of $R$ and
an irreducible non-linear character of $\tilde{D}$ if $\tilde{\psi}$
is irreducible. Since $\tilde{D}$ is an extra-special $2$ group, $D$
has only one non-linear irreducible character up to isomorphisms
(cf.\ \cite{Go} and \cite[Theorem 5.5.1]{FLM}). Therefore, the
number of inequivalent irreducible $\varepsilon$-projective
representation of $C$ is equal to the order of $R=C\cap C^\perp$.

Let us review the structure of the irreducible non-linear
character of $\tilde{D}$ in more detail.
Let $H$ be a maximal self-orthogonal subcode of $D$.
Then by \eqref{eq:4.2} the preimage
$\pi^{-1}_{\C^*}(H)$ of $H$ in $\C^*\tilde{C}$ splits.
Hence, there exists a map $\iota : H \to \C^*$ such that
$\varepsilon(\alpha,\beta)=(\partial \iota)(\alpha,\beta)
=\iota(\alpha)\iota(\beta)/\iota(\alpha+\beta)$ for all
$\alpha,\beta\in H$.
Since $\varepsilon (\alpha,\beta)\in\{ \pm 1\}$, one also
has $\varepsilon(\alpha,\beta)=\varepsilon(\alpha,\beta)^{-1}
=\iota(\alpha+\beta)/\iota(\alpha)\iota(\beta)$.
Then the section map
$H\ni \alpha \mapsto (\alpha,\iota(\alpha))\in \pi_{\C^*}^{-1}(H)$
is a group homomorphism.
Let $\chi$ be a linear character of $H$ and define a linear
character $\tilde{\chi}$ of $\pi_{\C^*}^{-1}(H)$ by
$\tilde{\chi}(\alpha, \iota(\alpha)u)=u \chi(\alpha)$ for
$\alpha \in H$ and $u\in \C^*$.
Since the preimage $\tilde{H}:=\pi^{-1}(H)$ is a subgroup of
$\pi_{\C^*}^{-1}(H)$, we may view $\tilde{\chi}$ as a linear
character of $\tilde{H}$.
Then the irreducible non-linear character of $\tilde{D}$ is
realized by the induced module
$\ind_{\tilde{H}}^{\tilde{D}}\tilde{\chi}$
(cf.\ Theorem 5.5.1 of \cite{FLM}).
Summarizing, we have:

\begin{proposition}\label{prop:4.1}
  (Theorem 5.5.1 of \cite{FLM})
  Let $\psi$ be an irreducible $\varepsilon$-projective representation of $C$.
  Then the associated linear representation $\tilde{\psi}$ of $\C^*\tilde{C}$
  is of the form
  $\lambda \tensor_\C \ind_{\tilde{H}}^{\tilde{D}} \tilde{\chi}$,
  where $\lambda$ is a linear character of $C\cap C^\perp$, $\tilde{H}$ is
  the preimage of a maximal self-orthogonal subcode $H$ of $D$ in $\tilde{C}$,
  and $\tilde{\chi}$ is a linear character of $\tilde{H}$ such that
  $\tilde{\chi}(0,-1)=-1$.
  In particular, $\tilde{\psi}$ is induced from a linear character of
  a maximal abelian subgroup of $\tilde{C}$.
\end{proposition}

\subsection{Structure of modules}

Let $C$ be an even linear code of $\Z_2^n$.
For a codeword $\alpha =(\alpha_1,\dots,\alpha_n)\in C$, we set
$$
   V(\alpha):=L(\shf,\alpha_1/2)\tensor \cds \tensor L(\shf,\alpha_n/2).
$$
Let $V_C=\oplus_{\alpha\in C} V(\alpha)$ be the code VOA
associated to $C$.
Since  $V(0)= L(\shf,0)^{\tensor n}$ is a rational full sub VOA of $V_C$,
every $V_C$-module is completely reducible as a $V(0)$-module.
We shall review the structure theory of irreducible $V_C$-modules
from \cite{M2,L3,Y1,Y2}.

Let $M$ be an irreducible $V_C$-module.
Take an irreducible $V(0)$-submodule $W$ of $M$, which is possible as
$V(0)$ is rational.
Let $\tau(W) \in \Z_2^n$  be the binary $1/16$-word of $W$ as defined
in \eqref{eq:1.1} (see also Section \ref{sec:3.2}).
Then it follows from the fusion rules of $L(\shf,0)$-modules that
$\tau(W)\in C^\perp$ and $\tau(W)=\tau(W')$ for any irreducible
$V(0)$-submodule $W'$ of $M$.
Set $C_W:=\{\alpha \in C \mid V(\alpha) \fusion_{V(0)} W\simeq W\}$.
Then $C_W= \{ \alpha \in C \mid \supp(\alpha)\subset \supp(\tau(W))\}$
and $C_{W'}=C_W$ for any irreducible $V(0)$-submodule $W'$ of $M$.
Let $\{ \alpha_i \mid 1\leq i\leq r\}$ be the coset representatives
for $C_W$ in $C$.
By the definition of $C_W$, it follows
$V(\alpha_i)\fusion_{V(0)}W\not\simeq V(\alpha_j)\fusion_{V(0)} W$
if $i\neq j$, because if
$V(\beta)\fusion_{V(0)} W=V(\gamma)\fusion_{V(0)} W$
in the fusion algebra then
$W=V(\beta)\fusion_{V(0)}V(\gamma)\fusion_{V(0)} W
= V(\beta+\gamma)\fusion_{V(0)}W$
for $\beta,\gamma\in C$.
Note that the fusion algebra associated to $V(0)$ is associative
and $V(\beta)\fusion_{V(0)}V(\gamma)=V(\beta+\gamma)$.
For simplicity, we set $W^i:=V(\alpha_i) \fusion_{V(0)} W$.
Then we have the following isotypical decomposition:
$$
  M = \bigoplus_{i=1}^r W^i \tensor \hom_{V(0)}(W^i,M).
$$
In the decomposition above, each homogeneous component
$$
  W^i \tensor \hom_{V(0)}(W^i,M)
$$
of $M$ forms an irreducible $V_{C_W}$-submodule,
where $V_{C_W}$ is the code VOA associated to $C_W$.
Let $U:=\hom_{V(0)}(W,M)$. It is shown in \cite{M2,L3,Y2} that
$U$ is an irreducible $\varepsilon$-projective representation
of $C_W$ so that $U$ is also an irreducible $\C^*\tilde{C}_W$-module.
Moreover, the $V_C$-module structure on $M$ is uniquely determined
by the $\C^*\tilde{C}_W$-module structure on $U$.

\begin{theorem}\label{thm:4.2}
  (\cite{M2,L3,Y2})
  Let $C$ be an even linear code and
  $V_C=\oplus_{\alpha \in C} V(\alpha)$ the associated code VOA.
  Let $W$ be an irreducible $V(0)$-module such that
  $\tau(W)\in C^\perp$.
  Then there is a one to one correspondence between
  the isomorphism classes of irreducible $\varepsilon$-projective
  representations   of $C_W$ and the isomorphism classes of
  irreducible $V_C$-modules containing $W$ as a $V^0$-submodule.
\end{theorem}

In the following, we shall give an explicit construction of
irreducible $V_C$-modules from irreducible $\varepsilon$-projective
$C_W$-modules.

\paragraph{An explicit construction}

Let $W$ be an irreducible $V(0)$-module such that the $1/16$-word
$\tau(W)\in C^\perp$.
Let $H$ be  a maximal self-orthogonal subcode of
$C_W= \{ \alpha \in C \mid \supp(\alpha)\subset \supp(\tau(W))\}$.
Since the preimage $\pi_{\C^*}^{-1}(H)$ of $H$ in \eqref{eq:4.4}
splits, there is a map $\iota: H \to \C^*$ such that
$(\alpha,\iota(\alpha))(\beta,\iota(\beta))
=(\alpha+\beta,\iota(\alpha+\beta))$ for all $\alpha,\beta\in H$.
Let $\chi$ be a linear character of $H$.
Then we can define a linear character
$\tilde{\chi}_\iota$ of $\pi_{\C^*}^{-1}(H)$ by
\begin{equation}\label{eq:4.5}
   \tilde{\chi}_\iota (\alpha, \iota(\alpha)u)
  = u \chi(\alpha)\quad  \text{ for }\ \alpha \in H,\ u\in \C^*.
\end{equation}
In this case, $\tilde{\chi}_\iota$ is also a linear character on
the preimage $\tilde{H}=\pi^{-1}(H)$ of $H$ in \eqref{eq:4.3}.
Let $\C^\varepsilon[C]$ be the twisted group algebra associated
to the 2-cocycle $\varepsilon\in Z^2(C,\C^*)$ defined
in \eqref{eq:4.1}.
That means
$\C^\varepsilon[C]=\Span_{\C}\{ e^\alpha \mid \alpha \in C\}$
as a linear space and
$e^\alpha e^\beta=\varepsilon(\alpha,\beta)e^{\alpha+\beta}$.
By \eqref{eq:4.2}, we have
\begin{equation}\label{eq:4.6}
  e^\alpha e^\beta=(-1)^{\la \alpha,\beta\ra} e^\beta e^\alpha .
\end{equation}
It is clear that $\C^\varepsilon[C_W]=\oplus_{\alpha\in C_W}\C
e^\alpha$ and $\C^\varepsilon[H]=\oplus_{\beta\in H}\C e^\beta$ are
subalgebras of $\C^\varepsilon[C]$.
Moreover,
$\C^\varepsilon[H]\simeq \C [H]$ as $\C$-algebras.
Let $\{ \alpha_1,\dots,\alpha_r\}$ be a set of coset representatives for
$C_W$ in $C$ and let $\{ \beta_1,\dots,\beta_s\}$ be a set of coset
representatives for $H$ in $C_W$.
Consider an induced module
$\ind_{\tilde{H}}^{\tilde{C}}\tilde{\chi}_\iota$.
As a linear space, it is defined by
$$
  \ind_{\tilde{H}}^{\tilde{C}}\tilde{\chi}_\iota
  = \bigoplus_{i=1}^r \bigoplus_{j=1}^s \C \, e^{\alpha_i+\beta_j}\!
    \tensor_{\C^\varepsilon [\tilde{H}]}\! v_{\tilde{\chi}_\iota},
$$
where $\C v_{\tilde{\chi}_\iota}$ is a $\C^\varepsilon[H]$-module
affording the character $\tilde{\chi}_\iota$, that is,
$\iota(\alpha)e^\alpha\cd v_{\tilde{\chi}_\iota}
= \chi(\alpha) v_{\tilde{\chi}_\iota}$ for all $\alpha \in H$.
Note also that the components
$$
  U^i:=\bigoplus_{j=1}^s \C
  e^{\alpha_i+\beta_j}\tensor_{\C^\varepsilon[\tilde{H}]}
  v_{\tilde{\chi}_\iota},  \quad 1\leq i\leq r,
$$
are irreducible $\C^{\varepsilon}[C_W]$-modules. Set $W^i:=
V(\alpha_i)\fusion_{V(0)} W$ for $1\leq i\leq r$. Let
$I^{\alpha,i}(\cd,z)$ be a $V(0)$-intertwining operator of type
$V(\alpha) \times W^i \to V(\alpha) \fusion_{V(0)} W^i$. Since all
$V(\alpha)$, $\alpha\in C$, are simple current $V(0)$-modules,
$I^{\alpha,i}(\cd,z)$ are unique up to scalars. It is possible to
choose these intertwining operators such that
$$
\begin{array}{l}
  (z_0+z_2)^m I^{\alpha,j'}(x^\alpha,z_0+z_2)
              I^{\beta,j}(x^\beta,z_2)w^j
  \vsb\\
  =  \varepsilon(\alpha,\beta)(z_2+z_0)^m
    I^{\alpha+\beta,j}(Y_{V_C}(x^\alpha,z_0)x^\beta,z_2)w^j
\end{array}
$$
for $x^\alpha\in V^\alpha$, $x^\beta\in V^\beta$, $w^j\in W^j$,
$\alpha_{j'}+C_W=\beta+\alpha_j+C_W $ and $m\gg 0$ (cf.\
\cite{M2,Y2}). We can also choose $I^{0,i}(\cd,z)$ so that
$I^{0,i}(\vacuum,z)=\id_{W^i}$. Now put
$$
  M=\ind_{V_H}^{V_C}(W,\tilde{\chi}_\iota)
  := \bigoplus_{i=1}^r W^i\tensor_\C  U^i
$$
and define a vertex operator $Y(\cd,z):V_C\times M\to M(\!(z)\!)$ by
$$
  Y(x^\alpha,z)w^i\tensor_\C u^i
  := I^{\alpha,i}(x^\alpha,z)w^i\tensor_\C (e^\alpha \cd u^i)
$$
for $x^\alpha\in V^\alpha$, $w^i\in W^i$ and $u^i\in U^i$.

\begin{theorem}\label{thm:4.3}
  (\cite{M2,L3,Y1})
  The induced module $\ind_{V^0}^{V_C}(W,\tilde{\chi}_\iota)$ equipped
  with the vertex operator defined above is  an irreducible $V_C$-module.
  Moreover, every irreducible $V_C$-module is isomorphic to
  an induced module.
\end{theorem}


\begin{remark}\label{rem:4.4}
  Even if $\tau(W)\not\in C^\perp$, one can still define an
  irreducible   $\Z_2$-twisted $V_C$-module structure on
  $\ind_{V_H}^{V_C}(W,\tilde{\chi})$ (cf.\ \cite{L1,Y1}).
\end{remark}

\paragraph{Parameterization by a pair of binary codewords}
The irreducible $V_C$-modules can also be parameterized
by a pair of binary codewords.
For given $\beta\in C^\perp$ and $\gamma\in \Z_2^n$,
we define a weight vector
$h_{\beta,\gamma}=(h^1_{\beta,\gamma},\dots,h^n_{\beta,\gamma})$,
$h_{\beta,\gamma}^i \in \{ 0,1/2,1/16\}$ by
$$
  h^i_{\beta,\gamma}:=
  \begin{cases}
  \dfrac{1}{16} & \text{  if } \beta_i=1,
  \vsb\\
  \dfrac{\gamma_i}2 & \text{ if } \beta_i=0.
  \end{cases}
$$
Let
$$
  L(h_{\beta,\gamma}):= L(\shf,h^1_{\beta,\gamma})\tensor \cds
  \tensor L(\shf,h^n_{\beta,\gamma})
$$
be the irreducible $L(\shf,0)^{\tensor n}$-module with the weight
$h_{\beta,\gamma}$. Set $C_\beta:=\{ \alpha \in C \mid
\supp(\alpha)\subset \supp(\beta)\}$ and let $R^\beta = C_\beta\cap
(C_\beta)^\perp$
be the radical of $C_\beta$. Fix a map $\iota: R^\beta \to \C^*$
such that the section map $R^\beta\ni \alpha \mapsto (\alpha,
\iota(\alpha)) \in \pi_{\C^*}^{-1}(R^\beta)$ is a group
homomorphism. Take a maximal self-orthogonal subcode $H$ of
$C_\beta$. Then $R^\beta\subset H$ and we can extend $\iota$ to $H$
such that the section map $H\ni \alpha \mapsto
(\alpha,\iota(\alpha))\in \pi^{-1}_{\C^*}(H)$ is a group
homomorphism. For, there exists a map $\jmath: H\to \C^*$ such that
$H\ni \alpha \mapsto (\alpha,\jmath(\alpha))\in \pi_{\C^*}^{-1}(H)$
is a group homomorphism. Then $\mu=\iota/\jmath$ restricted on
$R^\beta$ is a character since $\partial \mu=\partial \iota/\partial
\jmath =\varepsilon/\varepsilon=1$ on $R^\beta$. Take a complement
$K$ such that $H=R^\beta\oplus K$ and extend $\mu$ to $H$ by letting
$\mu(K)=1$. Then $\mu \jmath$ coincides with $\iota$ on $R^\beta$ as
desired. Define a character $\chi_\gamma$ of $H$ by
$\chi_\gamma(\alpha):=(-1)^{\la \gamma,\alpha\ra}$ for $\alpha\in H$
and extend to a character $\tilde{\chi}_{\gamma;\iota}$ of
$\pi_{\C^*}^{-1}(H)$ by $
\tilde{\chi}_{\gamma;\iota}(\alpha,\iota(\alpha)u)
:=\chi_\gamma(\alpha)u$ for $\alpha \in H$ and $u\in \C^*$. Note
that $\tilde{\chi}_{\gamma;\iota}$ also defines a linear character
on $\tilde{H}$. Moreover, every character $\varphi$ of $\tilde{H}$
such that $\varphi(0,-1)=-1$ is of the form $\tilde{\chi}_{\gamma;
\iota}$ for some $\gamma\in \Z_2^n$. Then by Theorem \ref{thm:4.3},
the pair $(\beta,\gamma)$ determines an irreducible $V_C$-module
\[
  M_C(\beta,\gamma; \iota)
  :=\ind_{V_H}^{V_C}(L(h_{\beta,\gamma}),\tilde{\chi}_{\gamma;\iota}).
\]
Note that if $C$ is self-orthogonal and $\supp(C)\subset
\supp(\beta)$, then $ M_C(\beta,\gamma;\iota)\simeq L(
h_{\beta,\gamma})$ as a $V(0)$-module. If $\be=0$, then $H=0$ and
$\iota$ is trivial. We shall simply denote $M_C(0,\gamma; \iota)$ by
$M_C(0,\gamma)$. It is also obvious that
$$
  M_C(0,\gamma)
  = \bigoplus_{\alpha=(\alpha_1,\dots,\alpha_n)\in C+\gamma}
    L(\shf,\alpha_1/2)\tensor \cds \tensor L(\shf,\alpha_n/2).
$$
This module is called a {\it coset module} in
\cite{M2}\footnote{This name has nothing to do with so-called the
{\it coset construction} (cf.\ \cite{FZ,GKO}) of a commutant
subalgebra.}. We sometimes denote $M_C(0,\gamma)$  by
$V_{C+\gamma}$, also. \vsb

We shall review some basic properties of $M_C(\beta,\gamma;\iota)$.

\begin{lemma}\label{lem:4.5}  (\cite{DGL})
  The module structure of $M_C(\beta,\gamma; \iota)$ is
  independent of the choice of the   maximal self-orthogonal
  subcode $H$ of $C_\beta$ and the choice of the extension
  of $\iota$ from $R^\beta$ to $H$.
\end{lemma}

\pf
Let $H'$ be a maximal self-orthogonal subcode of $C_\beta$ and
$\iota': H' \to \C^*$ an extension of $\iota$ to $H'$ such that
the section map
$H'\ni \alpha\mapsto (\alpha,\iota'(\alpha))\in \pi_{\C^*}^{-1}(H')$
is a group homomorphism. Define a linear character
$\tilde{\chi}'_{\gamma;\iota'}$ of $\pi_{\C^*}^{-1}(H')$ by
$\tilde{\chi}'_{\gamma;\iota'}(\alpha,\iota'(\alpha)u)
:= (-1)^{\la \gamma,\alpha\ra}u$ for $\alpha \in H'$ and $u\in \C^*$.
Then $\tilde{\chi}'_{\gamma;\iota'}$ is also a linear character
of the preimage $\tilde{H}'$ of $H'$ in \eqref{eq:4.3}.
We shall show that
$$
  \ind_{V_{H'}}^{V_C}(L(h_{\beta,\gamma}),\tilde{\chi}_{\gamma;\iota'})
  \simeq
  \ind_{V_H}^{V_C}(L(h_{\beta,\gamma}),\tilde{\chi}_{\gamma,\iota}).
$$
For this, it suffices to show that
$\ind_{\tilde{H}}^{\tilde{C}_\beta}\tilde{\chi}_{\gamma;\iota}
\simeq \ind_{\tilde{H}'}^{\tilde{C}_\beta}\tilde{\chi}'_{\gamma;\iota'}$
by Theorem \ref{thm:4.2} and the construction of induced modules.
By definition, it is clear that
$\tilde{\chi}_{\gamma;\iota}|_{R^\beta}
= \tilde{\chi}'_{\gamma;\iota'}|_{R^\beta}$.
For simplicity, we denote it by $\lambda$.

Let $K$ be a compliment of $R^\beta$ in $H+H'$ and
take $D$ be a complement of $R^\beta$ in $C_\beta$ such that
$K\subset D$. Then $H=R^\beta\oplus (H\cap K)$, $H'=R^\beta\oplus
(H'\cap K)$ and $\tilde{C}_\beta\simeq \tilde{R}^\beta *_{\{\pm 1\}}
\tilde{D}$. It is obvious that both $H_1=H\cap K$ and $H_2=H'\cap K$
are maximal self-orthogonal subcodes of $D$. Therefore, by
Proposition \ref{prop:4.1}, we have
$\ind_{\tilde{H}}^{\tilde{C}_\beta}\tilde{\chi}_{\gamma; \iota}
\simeq \lambda \tensor_\C \ind_{\tilde{H}_1}^{\tilde{D}}
\tilde{\chi}_{\gamma;\iota}|_{\tilde{H}_1}$ and
$\ind_{\tilde{H}'}^{\tilde{C}_\beta}\tilde{\chi'}_{\gamma;\iota'}
\simeq \lambda \tensor_\C \ind_{\tilde{H}_2}^{\tilde{D}}
\tilde{\chi}'_{\gamma;\iota'}|_{\tilde{H}_2}$.
Since there is only one linear representation of $\tilde{D}$
such that $(0, -1)$ acts non-trivially, we have
$\ind_{\tilde{H}_1}^{\tilde{D}} \tilde{\chi}_{\gamma;\iota}
\simeq \ind_{\tilde{H}_2}^{\tilde{D}} \tilde{\chi}'_{\gamma;\iota'}$
and
$\ind_{\tilde{H}}^{\tilde{C}_\beta}\tilde{\chi}_{\gamma;\iota}
\simeq \ind_{\tilde{H}'}^{\tilde{C}_\beta}\tilde{\chi}'_{\gamma;\iota}$
as desired
\qed

\vsb

\begin{remark}\label{rmk2}
If we choose another map $\iota': R^\beta \to \C^*$ such that the
section map $R^\beta\ni \alpha \mapsto (\alpha,\iota'(\alpha))\in
\pi_{\C^*}^{-1}(R^\beta)$ is a group homomorphism, then
$\iota/\iota'$ is a linear character of $R^\beta$. Thus, there
exists $\xi\in (\Z_2^n)_\beta$ such that $\iota(\al)/\iota'(\al)=
(-1)^{\la \alpha,\xi\ra}$ for $\alpha\in R^\beta$. Hence, for any
$\alpha \in R^\beta$ and $u\in \C^*$, we have
$$
\begin{array}{ll}
  \tilde{\chi}_{\gamma;\iota'}(\alpha, \iota(\alpha)u)
  &= \tilde{\chi}_{\gamma;\iota'}
     (\al, (-1)^{\la\alpha,\xi\ra}\iota'(\alpha)u)
  = (-1)^{\la\alpha, \xi\ra}\cd (-1)^{\la \alpha,\gamma\ra} u
  \vsb\\
  &= (-1)^{\la \alpha, \gamma+\xi\ra } u
  =\tilde{\chi}_{\gamma+\xi;\iota}(\al, \iota(\al)u) .
\end{array}
$$
Hence, $\tilde{\chi}_{\gamma;\iota'}= \tilde{\chi}_{\gamma+\xi;\iota}$
on $\tilde{R}^\beta$ and we have
$M_C(\be,\gamma; \iota')\simeq M_C(\be,\gamma+\xi; \iota)$.
\end{remark}

Similarly, one can show the following by considering linear
characters of $\tilde{H}$.

\begin{lemma}\label{lem:4.6}  (\cite{DGL})
    Let $\beta,\beta'\in C^\perp$ and $\gamma,\gamma'\in \Z_2^n$.
  Then the irreducible $V_C$-modules $M_C(\beta,\gamma;\iota)$ and
  $M_C(\beta',\gamma';\iota)$ are isomorphic if and only if
  $$
    \beta=\beta' \q \mbox{and} \q \gamma+\gamma' \in C+ H^{\perp_\beta},
  $$
  where $H$ is a maximal self-orthogonal subcode of $C_{\beta}$ and
  $H^{\perp_\beta}=\{ \alpha \in \Z_2^n \mid \supp(\alpha)\subset
  \supp(\beta) \mbox{ and } \la \alpha,\delta\ra =0 \mbox{ for all }
  \delta \in H\}$.
\end{lemma}

\pf
Assume that $M_C(\beta,\gamma;\iota)\simeq M_C(\beta',\gamma';\iota)$.
Then clearly $\beta=\beta'$ by 1/16-word decompositions.
It is also obvious from the definition of $M_C(\beta,\gamma;\iota)$
that $M_C(\beta,\gamma;\iota)\simeq M_C(\beta,\gamma+\delta;\iota)$
for any $\delta\in H^{\perp_\beta}$.
Let $\{\alpha_1,\dots,\alpha_r\}$ and $\{\delta_1,\dots,\delta_s\}$
be transversals for $C_\beta$ in $C$ and $H$ in $C_\beta$,
respectively.
Then by definition we have a decomposition
$$
  M_C(\beta,\gamma;\iota)
  =\bigoplus_{i=1}^r \bigoplus_{j=1}^s
  \left(V(\alpha_i)\fusion_{V(0)} L(h_{\beta,\gamma})\right)
  \tensor_\C \left(e^{\alpha_i+\delta_j}
  \tensor_{\C^\varepsilon[H]}\tilde{\chi}_{\gamma;\iota}\right).
$$
It follows from \eqref{eq:4.6} that
\begin{equation}\label{proof4.6}
  \left( V(\alpha_i)\fusion_{V(0)} L(h_{\beta,\gamma})\right)
  \tensor_\C \left(e^{\alpha_i+\delta_j}
  \tensor_{\C^\varepsilon[H]}\tilde{\chi}_{\gamma; \iota}\right)
  \simeq M_H(\beta,\gamma+\alpha_i+\delta_j;\iota)
\end{equation}
as a $V_H$-submodule. Therefore, we have the following
decompositions:
$$
\begin{array}{lll}
  M_C(\beta,\gamma;\iota)
  &=& \ds \bigoplus_{\delta+H\in C/H}M_H(\beta,\gamma+\delta;\iota),
  \vsb\\
  M_C(\beta,\gamma';\iota)
  &=& \ds \bigoplus_{\delta+H\in C/H}M_H(\beta,\gamma'+\delta;\iota).
\end{array}
$$
Since $H=C_\beta\cap H^{\perp_\beta}$ by the maximality of $H$, all
$M_H(\beta,\gamma+\delta;\iota)$, $\delta\in C/H$, are inequivalent
irreducible $V_H$-submodules.
Thus, if $M_C(\beta,\gamma;\iota)\simeq M_C(\beta,\gamma';\iota)$,
then
$\tilde{\chi}_{\gamma';\iota}=\tilde{\chi}_{\gamma+\delta;\iota}$
for some $\delta\in C$.
This is possible if and only if $\gamma+\gamma'\in C+H^{\perp_\beta}$.
Conversely, if $\gamma+\gamma'\in C+H^{\perp_\beta}$, then
$M_C(\beta,\gamma;\iota)$ and $M_C(\beta,\gamma';\iota)$
contain isomorphic irreducible $V_H$-submodules.
Since $V_C$-module structures on $M_C(\beta,\gamma;\iota)$ and
$M_C(\beta,\gamma';\iota)$ are uniquely determined by their
$V_H$-module structures, they are isomorphic.
\qed
\vsb

In the proof above, we have shown the following useful fact.

\begin{corollary}\label{cor:4.7}
  Let $M_C(\beta,\gamma; \iota)$ be an irreducible $V_C$-module.
  Let $H$ be a maximal self-orthogonal subcode of $C_\beta$.
  Then as a $V_H$-module,
  $$
    M_C(\beta,\gamma;\iota)
    = \bigoplus_{\delta+H\in C/H} M_H(\beta,\gamma+\delta; \iota).
  $$
  In particular, every irreducible $V_H$-submodule of
  $M_C(\beta,\gamma; \iota)$ is multiplicity-free.
\end{corollary}

\begin{lemma}
Let $R$ be the radical of $C_\be$ with respect to the standard
bilinear form and $H$ a maximal self-orthogonal subcode of $C_\be$.
Then $C_\be+ H^{\perp_{\be}}= R^{\perp_{\be}}$ and hence the code
$C+ H^{\perp_\beta}= C+R^{\perp_\beta}$ is again independent of the
choice of the maximal self-orthogonal subcode $H$.
\end{lemma}

\pf Since $R\subset H$, we have $H^{\perp_{\be}}\subset
R^{\perp_{\be}}$. By definition, it is also clear that $C_\beta
\subset R^{\perp_{\be}}$ and we have $C_\be+ H^{\perp_{\be}} \subset
R^{\perp_{\be}}$.  Since the bilinear form $\la \cdot, \cdot\ra$
restricted on the quotient space $C_\beta /R$ is non-degenerate and
$H/R$ is a maximal self-orthogonal subspace, we have $\dim
C_\beta/R= 2\dim H/R$ and hence $\dim C_\be = 2\dim H-\dim R$. On
the other hand,
\[
\begin{split}
\dim (C_\be+H^{\perp_{\be}})&= \dim C_\be+\dim
H^{\perp_{\be}}-\dim (C_\be \cap H^{\perp_{\be}})\\
&=\dim C_\be+ \wt\, \be -2\dim H= \wt\, \be - \dim R=\dim
R^{\perp_{\be}}.
\end{split}
\]
Thus, we have $C_\be+ H^{\perp_{\be}}= R^{\perp_{\be}}$ \qed

\subsection{Dual module}

We shall determine the structure of the dual module of a
$V_C$-module $M_C(\beta,\gamma;\iota)$. Recall that  the dual module
of a $V$-module $M=\oplus_{n\in \N} M_{n+h}$ is defined to be its
restricted dual $M^*=\oplus_{n\in \N} M^*_{n+h}$ equipped with a
vertex operator $Y_M^*(\cd,z)$ defined by \eqref{eq:2.2}.

First, we consider the case when the code is self-orthogonal.
Let $H$ be a self-orthogonal code of $\Z_2^n$.
In this case, one can define a character $\varphi$ of $H$ by
$\varphi(\alpha)=(-1)^{\wt(\alpha)/2}$ for $\alpha \in H$.
So there exists a codeword $\kappa \in \Z_2^n$ such that
$\varphi(\alpha)=(-1)^{\la \kappa,\alpha\ra}$ for all $\alpha \in H$.

\begin{lemma}\label{lem:4.8}
  Let $H\subset \Z_2^n$ be a self-orthogonal code.
  For any $\gamma \in \Z_2^n$, the dual module of
  $M_H((1^n),\gamma; \iota)$ is isomorphic to
  $M_H((1^n),\gamma+\kappa;\iota)$, where $\kappa\in \Z_2^n$ is
  given by $(-1)^{\la \kappa,\alpha\ra}=(-1)^{\wt (\alpha)/2}$
  for all $\alpha \in H$.
\end{lemma}

\pf By assumption, $M_H((1^n),\gamma; \iota)
= L(\shf,\sfr{1}{16})^{\tensor n} \tensor \tilde{\chi}_{\gamma;\iota}$
is an irreducible $V(0)=L(\shf,0)^{\tensor n}$-module.
Therefore, $$M_H((1^n),\gamma;\iota)^*\simeq M_H((1^n),\gamma';
\iota)\quad  \text{ for some }\gamma'\in \Z_2^n.$$ Since
$L(\shf,\sfr{1}{16})^{\tensor n}$ is a self-dual $V(0)$-module, we
have a $V(0)$- isomorphism $f: M\to M^*$. Let $Y(\cd,z)$ and
$Y^*(\cd,z)$ be the vertex operators on $M$ and $M^*$, respectively.
For $\alpha \in H$, let $x^\alpha\in V(\alpha)$ be a non-zero
highest weight vector of weight $\wt(\alpha)/2$. Then
$Y^*(x^\alpha,z)f$ and $fY(x^\alpha,z)$ are $V(0)$-intertwining
operators of type
$$
  V(\alpha) \times M_H((1^n),\gamma; \iota) \to M_H((1^n),\gamma;\iota)^*.
$$
Since the space of $V(0)$-intertwining operators of this type is
one-dimensional, there exists a scalar $\lambda_\alpha\in \C^*$ such
that $Y^*(x^\alpha,z)f=\lambda_\alpha fY(x^\alpha,z)$. Let $v$ be a
non-zero highest weight vector of $M_H((1^n),\gamma;\iota )$. Then
by \eqref{eq:2.2}, one has
$$
\begin{array}{ll}
  \la Y^*(x^\alpha,z)fv, v\ra
  &= (-1)^{\wt (\alpha)/2} z^{-\wt(\alpha)}
     \la fv,Y(x^\alpha,z^{-1}) v\ra
  \vsb\\
  &= (-1)^{\wt(\alpha)/2} z^{-\wt(\alpha)/2}
     \la fv, x^\alpha_{(\wt(\alpha)/2-1)}v\ra.
\end{array}
$$
On the other hand,
$$
\begin{array}{ll}
  \la Y^*(x^\alpha,z)fv,v\ra
  &= \lambda_\alpha \la fY(x^\alpha,z)v,v\ra
  \vsb\\
  &= \lambda_\alpha z^{-\wt(\alpha)/2}
     \la fx^\alpha_{(\wt(\alpha)/2-1)}v,v\ra.
\end{array}
$$
Since $x^\alpha_{(\wt(\alpha)/2-1)}v=tv$ for some $t\in \C^*$ and
$\la fv,v\ra\ne 0$, we have
$\lambda_\alpha=(-1)^{\wt(\alpha)/2}=(-1)^{\la \kappa,\alpha\ra}$.
Therefore, by considering the linear character associated to
$M_H((1^n),\gamma; \iota)^*$, we see that $M_H((1^n),\gamma;
\iota)^* \simeq M_H((1^n),\gamma+\kappa;\iota)$. \qed

\begin{proposition}\label{prop:4.9}
  Let $C$ be an even linear code, $\beta\in C^\perp$ and
  $\gamma\in \Z_2^n$.
  Let $H$ be a maximal self-orthogonal subcode of $C_\beta$.
  Then the dual module $M_C(\beta,\gamma;\iota)^*$ is isomorphic to
  $M_C(\beta,\gamma+\kappa_H; \iota)$ where $\kappa_H \in \Z_2^n$
  is such that $\supp(\kappa_H)\subset \supp(\beta)$  and
  $(-1)^{\la \kappa_H,\alpha\ra}=(-1)^{\wt(\alpha)/2}$
  for all $\alpha \in H$.
\end{proposition}

\pf
By Corollary \ref{cor:4.7}, $M_C(\beta,\gamma;\iota)$ contains a
$V_H$-submodule
$$
  M_H(\beta,\gamma;\iota)
  \simeq L(h_{\beta,\gamma})\tensor \tilde{\chi}_{\gamma;\iota}.
$$
By the previous lemma, the dual module $M_C(\beta,\gamma;\iota)^*$
contains a $V_H$-submodule isomorphic to
$$
  M_H(\beta,\gamma;\iota)^*
  \simeq
  M_H(\beta,\gamma+\kappa_H;\iota)
  = L(h_{\beta,\gamma})\tensor\tilde{\chi}_{\gamma+\kappa_H;\iota}.
$$
Therefore, by the structure of irreducible $V_C$-modules,
$$
  M_C(\beta,\gamma;\iota)^*
  \simeq \ind_{V_H}^{V_C}(L(h_{\beta,\gamma}),
    \tilde{\chi}_{\gamma+\kappa_H;\iota})
$$
and hence $M_C(\beta,\gamma;\iota)^*\simeq
M_C(\beta,\gamma+\kappa_H;\iota)$. \qed \vsb

As an immediate corollary, we have:

\begin{corollary}\label{cor:4.10}
  With reference to the proposition above,
  $M_C(\beta,\gamma;\iota)$ is self-dual if and only if $\kappa_H\in C$.
  In particular, $M_C(0,\gamma)$ is self-dual for all $\gamma\in \Z_2^n$.
\end{corollary}

\subsection{Fusion rules}

We shall compute the fusion rules among some irreducible
$V_C$-modules. First, we recall a result from \cite{M2}
which is a direct consequence of Proposition \ref{prop:2.5}.

\begin{lemma}\label{lem:4.11}
  (\cite{M2})
  For $\alpha,\beta\in \Z_2^n$,
  $M_C(0,\alpha)\fusion_{V_C} M_C(0,\beta)=M_C(0,\alpha+\beta)$.
\end{lemma}

By the lemma above, we see that $M_C(0,\alpha)\fusion_{V_C}
M_C(0,\alpha)=M_C(0,0)\simeq V_C$. Therefore, all $M_C(0,\alpha)$,
$\alpha \in \Z_2^n$, are simple current $V_C$-modules by Corollary
\ref{cor:2.2}. It also follows that $M_C(0,\alpha)\fusion_{V_C}
M_C(\beta,\gamma;\iota)$ is an irreducible $V_C$-module with the
1/16-word $\beta$.
The corresponding fusion rules are also computed
by Miyamoto \cite{M3} in the case $\supp(\alpha)\subset \supp(\beta)$.

\begin{lemma}\label{lem:4.12}
  (\cite{M3})
  Let $\alpha,\beta,\gamma \in \Z_2^n$ with $\beta\in C^\perp$.
  Then
  $$
    M_C(0,\alpha)\fusion_{V_C}M_C(\beta,\gamma;\iota)
    =M_C(\beta,\alpha+\gamma;\iota).
  $$
  Moreover, the difference of the top weight of $M_C(\beta,\gamma;\iota)$
  and the top weight of $M_C(\beta,\alpha+\gamma;\iota)$ is congruent to
  $\la \al, \al+\be\ra/2 $ modulo $\Z$.
\end{lemma}

\pf
The assertion is proved in Lemma 3.27 of \cite{M3} in the case
$\supp(\alpha)\subset \supp(\beta)$.
We generalize his argument to obtain the desired fusion rule.
Since $M_C(0,\alpha)$ is a simple current $V_C$-module, we know that
there exists $\gamma'\in \Z_2^n$ such that
$$
  M_C(0,\alpha)\fusion_{V_C} M_C(\beta,\gamma;\iota)
  \simeq M_C(\beta,\gamma';\iota) .
$$
Therefore, if we can construct a non-zero
$V_C$-intertwining operator of type $M_C(0,\alpha)\times
M_C(\beta,\gamma;\iota)\to M_C(\beta,\alpha+\gamma;\iota)$,
then we are done.
To do this, we have to extend $V_C$ to a larger algebra.
The case $\alpha\in C$ is trivial so that we assume that $\alpha
\not\in C$. Set $C'=C\sqcup (C+\alpha)$. We can define a simple
vertex operator (super)algebra structure on the space $V_{C'}=V_C
\oplus V_{C+\alpha}=M_C(0,0)\oplus M_C(0,\alpha)$. This is a VOA if
$\alpha$ is even, and an SVOA if $\alpha$ is odd.

Set $H':=C'_\beta\cap H^\perp$, which is the unique maximal subcode
of $C'_\beta$ containing $H$ such that its preimage $\tilde{H}'$ is
a maximal abelian subgroup of $\tilde{C}'_\beta$ in the extension
\eqref{eq:4.3}.
We can take
$\jmath:H'\to \C^*$ such that $\jmath|_H=\iota$ and the section map
$H'\ni \delta\mapsto (\delta,\jmath(\delta))\in \pi_{\C^*}^{-1}(H')$
defines a group homomorphism.  In the definition of the induced
module
$M_C(\beta,\gamma;\iota)=
\ind_{V_H}^{V_C}L(h_{\beta,\gamma},\tilde{\chi}_{\gamma;\iota})$,
if we use $\C^\varepsilon[C']$ instead of $\C^\varepsilon[C]$ and
replace $\iota$ by $\jmath$,  then we obtain an irreducible
$V_{C'}$-module
$$
M_{C'}(\beta,\gamma; \jmath) :=
\ind_{V_{H'}}^{V_{C'}}(L(h_{\beta,\gamma},\tilde{\chi}_{\gamma,
\jmath})
$$
which contains $M_C(\beta,\gamma; \iota)=\ind_{V_{H}}^{V_{C}}
(L(h_{\beta,\gamma},\tilde{\chi}_{\gamma, \jmath|_H}) $ as a
$V_C$-submodule. The induced module $M_{C'}(\beta,\gamma; \jmath)$
is an untwisted $V_{C'}$-module if $\la \alpha,\beta\ra=0$ and
otherwise it is a $\Z_2$-twisted $V_{C'}$-module (cf.\ \cite{Y2}).
Nevertheless, the subspace $M=V_{C+\alpha}\cd
M_C(\beta,\gamma;\iota)$ of $M_{C'}(\beta,\gamma;\jmath)$ is an
irreducible $V_C$-submodule. It follows from \eqref{eq:4.6} and
\eqref{proof4.6} that there exists an irreducible $V_H$-submodule of
$M$ isomorphic to $M_H(\beta,\alpha+\gamma;\jmath|_H)$. Then
$M\simeq M_C(\beta,\alpha+\gamma; \jmath|_H)$ by the structure of an
irreducible $V_C$-module. Since $M\simeq
M_C(0,\alpha)\fusion_{V_C}M_C(\beta,\gamma;\iota)$, we obtain the
desired fusion rule.

Since the $L(\shf,0)$-intertwining operators of types
$L(\shf,h)\times L(\shf,\sfr{1}{16})\to L(\shf,\sfr{1}{16})$, $h\in
\{0,1/2\}$, keep the top weights but the $L(\shf,0)$-intertwining
operators of type $L(\shf,1/2)\times L(\shf,0)\to
L(\shf,\sfr{1}{2})$ and $L(\shf,1/2)\times L(\shf,\sfr{1}{2})\to
L(\shf,0)$ change the top weights by $1/2$ or $-1/2$, the difference
of top weights is as in the assertion.  \qed \vsb

By this lemma, we can compute the following fusion rule.

\begin{proposition}\label{prop:4.13}
  Let $\beta\in C^\perp$ and $\gamma\in \Z_2^n$.
  Let $H$ be a maximal self-orthogonal subcode of $C_\beta$.
  Then
  $$
    M_C(\beta,\gamma;\iota)\fusion_{V_C} M_C(\beta,\gamma;\iota)^*
    = \sum_{\delta+C\in C+H^{\perp_\beta}} M_C(0,\delta),
  $$
  where $\delta\in \Z_2^n$ runs over a transversal
  for $C$ in $C+H^{\perp_\beta}$.
\end{proposition}

\pf
It follows from the 1/16-word decomposition that
the fusion product
$$
  M_C(\beta,\gamma; \iota)\fusion_{V_C}M_C(\beta,\gamma;\iota)^*
$$
contains only modules of type $M_C(0,\delta)$.
Now assume that
$\binom{M_C(0,\delta)}{M_C(\beta,\gamma;\iota)\q
M_C(\beta,\gamma;\iota)^*}_{V_C}\ne 0$. Then by the symmetry of
fusion rules, we have
$$
  \binom{M_C(0,\delta)}{M_C(\beta,\gamma;\iota)\q M_C(\beta,\gamma;\iota)^*}_{V_C}
  \simeq \binom{M_C(\beta,\gamma;\iota)}{M_C(\beta,\gamma;\iota)\q M_C(0,\delta)}_{V_C}
  \ne 0.
$$
Since $M_C(\beta,\gamma;\iota)\fusion_{V_C}
M_C(0,\delta)=M_C(\beta,\gamma+\delta;\iota)$ by the previous lemma,
this is possible if and only if $\delta\in C+H^{\perp_\beta}$ by
Lemma \ref{lem:4.6}.
Therefore, we have the fusion rule as stated.
\qed
\vsb

By the lemma above, we introduce the following definition.

\begin{definition}\label{df:4.14}
  Let $\beta\in \Z_2^n$ and $H$ a subcode with
  $\supp (H) \subset \supp (\beta)$.
  $H$ is said to be \emph{self-dual with respect to $\beta$}
  if $H=H^{\perp_\beta}$.
\end{definition}

\begin{remark}\label{rem:4.15}
  Note that if $H$ is a self-dual subcode of $C_\be$ w.r.t.\ $\be$
  then $C+H^{\perp_\beta}=C$.
\end{remark}

By  Corollary \ref{cor:2.2} and Proposition \ref{prop:4.13},
we have

\begin{corollary}\label{cor:4.16}
  $M_C(\beta,\gamma;\iota)$ is a simple current module if and only if
  $C_\beta$ contains a self-dual subcode w.r.t.\ $\beta$.
\end{corollary}

\begin{remark}\label{rem:4.17}
Now suppose $M_C(\beta,0;\iota)\fusion_{V_C}
M_C(\beta,0;\iota)^*=\sum_{i=1}^p M_C(0,\delta_i)$. Let $H$ be a
maximal self-orthogonal subcode of $C_\beta$, and let $\kappa_H \in
(\Z_2^n)_\beta$ such that $\la \kappa_H,\alpha\ra=\la
\alpha,\alpha\ra/2 \mod 2$ for all $\alpha \in H$ as in Proposition
\ref{prop:4.9}. Then
$$M_C(\beta,0;\iota)^*=M_C(\beta,\kappa_H;\iota)
=M_C(0,\kappa_H)\fusion_{V_C}M_C(\beta,0;\iota)$$ and thus
$M_C(\beta,0;\iota)=M_C(0,\kappa_H)\fusion_{V_C}M_C(\beta,0;\iota)^*$.
Using this, we can compute the following fusion rule:
$$
\begin{array}{l}
  \ds M_C(\beta,\gamma_1;\iota)\fusion_{V_C}
   M_C(\beta,\gamma_2;\iota)
  \\
  \ds
  = \l\{ M_C(0,\gamma_1)\fusion_{V_C} M_C(\beta,0;\iota) \r\}
     \fusion_{V_C} \l\{ M_C(0,\gamma_2) \fusion_{V_C}
     M_C(\beta,0;\iota)\r\}
  \vspace{3mm}\\
  \ds
  = M_C(0,\gamma_1+\gamma_2) \fusion_{V_C} M_C(\beta,0;\iota)
    \fusion_{V_C} M_C(\beta,0;\iota)
  \vsb\\
  \ds
  = M_C(0,\gamma_1+\gamma_2) \fusion_{V_C} M_C(\beta,0;\iota)
    \fusion_{V_C}
    \l\{ M_C(0,\kappa_H)\fusion_{V_C} M_C(\beta,0;\iota)^*\r\}
  \vspace{3mm}\\
  \ds
  = M_C(0,\gamma_1+\gamma_2+\kappa_H) \fusion_{V_C} M_C(\beta,0;\iota)
    \fusion_{V_C} M_C(\beta,0;\iota)^*
  \vsb\\
  \ds
  = M_C(0,\gamma_1+\gamma_2+\kappa_H) \fusion_{V_C}
    \l\{\dsum_{i=1}^p M_C(0,\delta_i)\r\}
  \\
  \ds
  = \dsum_{i=1}^p M_C(0,\gamma_1+\gamma_2+\kappa_H +\delta_i).
\end{array}
$$
\end{remark}

\section{Structure of framed VOAs}

We shall prove that every framed VOA is a simple current extension
of a code VOA.
This result has many fruitful consequences.
For  example, the irreducible representations of a framed VOA can be
determined by a notion of induced modules.
Another interesting result is the conditions on possible structure
codes of holomorphic framed VOAs, namely we obtain a necessary
and sufficient condition for a pair of codes $(C,D)$ to be
the structure codes of some holomorphic framed VOAs
in Theorem \ref{thm:5.17}.

\subsection{Simple current structure}

In this subsection we discuss how a code VOA can be extended
to a framed VOA.
First, we give a construction of a non-trivial simple current extension.

\begin{lemma}\label{lem:5.1}
  Let $C$ be an even linear subcode of $\Z_2^n$ and
  $\beta\in C^\perp$ a non-zero codeword. Let $\gamma\in \Z_2^n$
  be a binary codeword such that the irreducible $V_C$-module
  $M_C(\beta,\gamma;\iota)$ has an integral top weight.
  If $C_\beta$ contains a doubly even self-dual subcode
  w.r.t.\ $\beta$,
  then there exists a unique structure of a framed VOA on
  $V_C\oplus M_C(\beta,\gamma;\iota)$
  which forms a $\Z_2$-graded simple current extension of $V_C$.
\end{lemma}

\pf Let $H$ be a doubly even self-dual subcode of $C_\beta$ w.r.t.\
$\beta$. By Proposition \ref{prop:4.9}, $M_C(\beta,\gamma;\iota)$ is
self-dual, and by Corollary \ref{cor:4.16},
$M_C(\beta,\gamma;\iota)$ is a simple current $V_C$-module. By
Corollary \ref{cor:4.7}, $M_C(\beta,\gamma;\iota)$ has a
$V_H$-module structure
$$
  M_C(\beta,\gamma;\iota)
  = \bigoplus_{\delta+H\in C/H} M_H(\beta,\gamma+\delta;\iota)
$$
where all irreducible $V_H$-submodules
$M_H(\beta,\gamma+\delta;\iota)$ are self-dual by Proposition
\ref{prop:4.9}. It is clear that a $V_C$-invariant bilinear form
on $M_C(\beta,\gamma;\iota)$ induces a non-degenerate $V_H$-invariant
bilinear form on $M_H(\beta,\gamma+\delta;\iota)$.
It is shown in \cite{Li1} that a $V_C$-invariant bilinear form on
$M_C(\beta,\gamma;\iota)$ is either symmetric or skew-symmetric.
Since the top level of $M_H(\beta,\gamma+\delta;\iota)$ is
one-dimensional, the $V_C$-invariant bilinear form on
$M_C(\beta,\gamma;\iota)$ must be symmetric. Therefore, $V_C\oplus
M_C(\beta,\gamma;\iota)$ forms a $\Z_2$-graded simple current
extension of $V_C$ by Proposition \ref{prop:2.6}.
\qed
\vsb

\begin{lemma}\label{lem:5.2}
  Let $\beta\in C^\perp$ and $\gamma \in \Z_2^n$ with $\beta\ne 0$.
  Assume that $V=V_C\oplus M_C(\beta,\gamma;\iota)$ forms a framed VOA.
  Then there exists a maximal self-orthogonal subcode $K$ of $C_\beta$
  which is doubly even .
\end{lemma}

\pf Let $H$ be a maximal self-orthogonal subcode of $C_\beta$. If
$H$ is doubly even, then we are done. So we assume that $H$ contains
a codeword whose weight is congruent to 2 modulo 4. Since
$V=V_C\oplus M_C(\beta,\gamma;\iota)$ forms a simple VOA,
$M_C(\beta,\gamma;\iota)$ is a self-dual $V_C$-module. Therefore, by
Corollary \ref{cor:4.10}, there exists a codeword $\kappa_H\in
C_\beta$ such that $(-1)^{\la
\kappa_H,\alpha\ra}=(-1)^{\wt(\alpha)/2}$ for all $\alpha \in H$. By
Corollary \ref{cor:4.7}, $M_C(\beta,\gamma;\iota)$ has the following
decomposition as a $V_H$-module:
$$
  M_C(\beta,\gamma;\iota)
  = \bigoplus_{\delta+H\in C/H} M_H(\beta,\gamma+\delta;\iota).
$$
By our choice of $H$, $\kappa_H\not\in H^{\perp_\beta}$ so that
$M_H(\beta,\gamma;\iota)$ and its dual
$M_H(\beta,\gamma+\kappa_H;\iota)$ are inequivalent irreducible
$V_H$-submodules of $V$. We consider a sub VOA $U$ generated by
$M_H(\beta,\gamma;\iota)\oplus M_H(\beta,\gamma+\kappa_H;\iota)$. By
the fusion rule given in Proposition \ref{prop:4.13}, $U$ has the
following shape as a $V_H$-module:
\begin{equation}\label{eq:5.1}
  U= M_H(0,0)\oplus M_H(0,\kappa_H)\oplus M_H(\beta,\gamma;\iota)
     \oplus M_H(\beta,\gamma+\kappa_H;\iota).
\end{equation}
Note that $H=C \cap H^{\perp_\beta}$ by the maximality of $H$.
Set $H':=H\sqcup (H+\kappa_H)$, $H_0:=H\cap \la \kappa_H\ra^\perp$
and take any $\alpha' \in H\setminus H_0$. Then $$H'=H_0\sqcup
(H_0+\alpha')\sqcup(H_0+\kappa_H)\sqcup(H_0+\alpha'+\kappa_H).$$ We
set $K:=H_0\sqcup(H_0+\kappa_H)$. It is clear that $U$ also
possesses a symmetric invariant bilinear form which we shall denote
by $\la \cd,\cd\ra_U$. Since $M_H(\beta,\gamma;\iota)$ and
$M_H(\beta,\gamma+\kappa_H;\iota)$ are dual to each other, we have
$$
  \la M_H(\beta,\gamma;\iota),M_H(\beta,\gamma;\iota)\ra_U
  = \la M_H(\beta,\gamma+\kappa_H;\iota),
    M_H(\beta,\gamma+\kappa_H;\iota)\ra_U
  =0.
$$
By Lemma \ref{lem:4.6}, $M_H(\beta,\gamma;\iota)$ and
$M_H(\beta,\gamma+\kappa_H;\iota)$ are isomorphic irreducible
$V_{H_0}$-modules and there exists a $V_{H_0}$-isomorphism
$\varphi:M_H(\beta,\gamma;\iota) \to M_H(\beta,\gamma+\kappa_H;\iota)$.
Let $u$ be a non-zero highest weight vector of
$M_H(\beta,\gamma;\iota)$.
Since the top level of $M_H(\beta,\gamma;\iota)$ is
one-dimensional, we may assume that
$\la u,\varphi(u)\ra_U=1$. Now consider the decomposition
\eqref{eq:5.1} of $U$ with respect to a series of sub VOAs
$V_{H_0}\subset V_K \subset V_{H'}$ of $U$. It is clear that
$M_H(0,0)\oplus M_H(0,\kappa_H)=M_K(0,0)\oplus M_K(0,\alpha')$.
Therefore, there exists a decomposition of $U$ as a $V_K$-module
$$
  U=M_K(0,0)\oplus M_K(0,\alpha')\oplus W
$$
with $W=M_H(\beta,\gamma;\iota)\oplus
M_H(\beta,\gamma+\kappa_H;\iota)$. Let $W^0$ be an irreducible
$V_K$-submodule of $W$. Since $K$ is a self-orthogonal subcode of
$C_\beta$, the top level of $W^0$ is one-dimensional. Let $v\in W^0$
be a non-zero highest weight vector. As we mentioned,
$M_H(\beta,\gamma;\iota)$ and $M_H(\beta,\gamma+\kappa_H;\iota)$ are
isomorphic $V_{H_0}$-submodules. But $M_H(\beta,\gamma;\iota)$ and
$M_H(\beta,\gamma+\kappa_H;\iota)$ {\it cannot} form
$V_K$-submodules by the fusion rule of $V_H$-modules. Therefore, we
can write $v=c_1 u +c_2 \varphi(u)$ with $c_1,c_2\ne 0$. This shows
that $\la v,v\ra_U = 2c_1c_2\ne 0$ so that $W^0$ is a self-dual
$V_K$-submodule. Then $K$ is a doubly even code by Corollary
\ref{cor:4.10}. Since $\abs{K}=\abs{H}=2\abs{H_0}$, $K$ is a maximal
self-orthogonal subcode of $C_\beta$. Therefore, $C_\beta$ contains
the desired subcode $K$.
\qed
\vsb

We recall the following fact from the coding theory.

\begin{theorem}\label{thm:5.3}
  (\cite{McST})
  Let $n$ be divisible by $8$ and $H$ a doubly even code of $\Z_2^n$
  containing the all-one vector $(11\dots1)\in \Z_2^n$.
  Then there exists a doubly even self-dual code $H'$ such that $H\subset H'$.
\end{theorem}

Now we begin to prove that every framed VOA is a simple current extension
of a code VOA.
For this, it suffices to show the following proposition.

\begin{proposition}\label{prop:5.4}
  Let $\beta\in C^\perp$ and $\gamma\in \Z_2^n$.
  Assume that $V=V_C\oplus M_C(\beta,\gamma;\iota)$ forms a framed VOA.
  Then $C_\beta$ contains a doubly even self-dual subcode w.r.t.\ $\beta$.
\end{proposition}

\pf By Lemma \ref{lem:5.2}, $C_\beta$ contains a maximal
self-orthogonal subcode $H$ which is doubly even. By Corollary
\ref{cor:4.7}, $V$ has a decomposition
$$
  V=\bigoplus_{\delta+H\in C/H} M_H(0,\delta)\oplus M_H(\beta,\gamma+\delta;\iota)
$$
as a $V_H$-module. By the fusion rule in Proposition \ref{prop:4.13},
the subspace
$$
  U:=V_H \oplus M_H(\beta,\gamma;\iota)
$$
forms a sub VOA of $V$, since $H=C_\beta \cap H^{\perp_\beta}$. If
$H$ is not self-dual, then there exists a doubly even self-dual
subcode $H'$ of $\Z_2^n$ w.r.t.\ $\beta$ such that $H\cup
(H+\beta)\subset H'$ by Theorem \ref{thm:5.3}. Note that the weight
of $\beta$ is divisible by 8 since $M_C(\beta,\gamma;\iota)$ has an
integral top weight. Let us consider the code VOA $V_{H'}$
associated to $H'$. Since $H\subset H'$, it is clear that $V_{H'}$
contains $V_H$ as a sub VOA.  We can also take a map $\jmath: H'\to
\C^*$ such that $\jmath|_H=\iota$ and the section map $H'\ni \alpha
\mapsto (\alpha,\jmath(\alpha))\in \pi_{\C^*}^{-1}(H')$ is a group
homomorphism. By the structure theory in Theorem \ref{thm:4.3}, we
can define an irreducible $V_{H'}$-module
$M_{H'}(\beta,\gamma;\jmath)$ such that
$M_{H'}(\beta,\gamma;\jmath)|_{V_H}\simeq M_H(\beta,\gamma;\iota)$
as a $V_H$-module. For simplicity, we shall denote
$M_{H'}(\beta,\gamma;\jmath)$ by $W$.  Since the top level of $W$ is
one-dimensional, the $V_{H'}$-invariant bilinear form on $W$ is
symmetric. Therefore, by Proposition \ref{prop:2.5}, we can define a
framed VOA structure on
$$
  U':= V_{H'}\oplus W .
$$
We denote the vertex operator on $U'$ by $Y'(\cd,z)$.
Now suppose $H$ is a proper subcode of $H'$.
Then
$$
  V_{H'}=\bigoplus_{\delta+H\in H'/H} V_{\delta+H},\q
  V_{\delta+H}= M_H(0,\delta),
$$
as a $V_H$-module.
Let $\pi_{\delta+H}: V_{H'}\to V_{\delta+H}$ be the projection map.
Then for $u,v\in W$, we have
$$
  Y'(u,z)v=\sum_{\delta+H\in H'/H} \pi_{\delta+H} Y'(u,z)v.
$$
Since the simple VOA structure is unique on $V_{H'}\oplus W$,
we may assume that $\pi_H Y'(u,z)v =Y_V(u,z)v$.
Take any $\alpha \in H'\setminus H$ and set $K:=H\sqcup (H+\alpha)$.
We shall show the following claim:
\paragraph{Claim} For $u,v\in W$, there exists $N=N(u,v)\in \N$ such that
\begin{equation}\label{eq:5.2}
  (z_1-z_2)^N Y'(u,z_1) \pi_{H+\alpha} Y'(v,z_2)w
  = (z_1-z_2)^N Y'(v,z_2) \pi_{H+\alpha} Y'(u,z_1)w
\end{equation}
for any $w\in W$.
\vsb\\
Take any $a\in V_{H+\alpha}$. Since $U=V_H\oplus
M_H(\beta,\gamma;\iota)$ forms a framed VOA by assumption, there
exists $N=N(u,v)\in \N$ such that
\begin{equation}\label{eq:5.3}
  (z_1-z_2)^N Y'(u,z_1)\pi_H Y'(v,z_2) w= (z_1-z_2)^N Y'(v,z_2) \pi_H Y'(u,z_1) w.
\end{equation}
Take a sufficiently large $k\in \N$.
Then one has
$$
\begin{array}{l}
  (z_1-z_2)^N (z_0-z_1)^k (z_0-z_2)^k
    Y'(u,z_1)\pi_{H+\alpha} Y'(v,z_2) Y'(a,z_0)w
  \vsb\\
  = (z_1-z_2)^N (z_0-z_1)^k (z_0-z_2)^k Y'(a,z_0) Y'(u,z_1)\pi_H Y'(v,z_2) w
  \vsb\\
  = (z_1-z_2)^N (z_0-z_1)^k (z_0-z_2)^k Y'(a,z_0) Y'(v,z_2)\pi_H Y'(u,z_1) w
  \vsb\\
  = (z_1-z_2)^N (z_0-z_1)^k (z_0-z_2)^k Y'(v,z_2)\pi_{H+\alpha} Y'(u,z_1) Y'(a,z_0)w .
\end{array}
$$
Since the expansions of both sides of the equations have only
finitely many negative powers of $z_0$, we get
$$
  (z_1-z_2)^N Y'(u,z_1)\pi_{H+\alpha} Y'(v,z_2) w
  = (z_1-z_2)^N  Y'(v,z_2)\pi_{H+\alpha} Y'(u,z_1) w .
$$
Note that $Y'(a,z)\pi_H=\pi_{H+\alpha}Y'(a,z)$ on $V_{H'}$ and
$W= V_{H+\alpha} \cd W$.
\vsb

By the Claim above, we can introduce a framed VOA structure on
$$
  X:=V_K\oplus M_K(\beta,\gamma; \jmath|_K)
$$
as follows. Since $W$ as
a $V_K$-module is isomorphic to $M_K(\beta,\gamma;\jmath|_K)$, we
can identify these structures. For $a,b\in V_K$ and $u,v\in
M_K(\beta,\gamma;\jmath|_K)$, we define the vertex operator map
$Y_X(\cd,z)$ by
$$
  Y_X(a,z)b:= Y'(a,z)b,\q Y_X(a,z)u:=Y'(a,z)u,\q Y_X(u,z)a:= Y'(u,z)a,
$$
and
$$
  Y_X(u,z) v:= \pi_H Y'(u,z)v +\pi_{H+\alpha} Y'(u,z)v.
$$
Let $\la \cd,\cd\ra_{U'}$ be an non-zero invariant bilinear form on
$U'$, which is unique up to scalar multiples. Since $V_K$ is a
subalgebra of $U'$, we can define an invariant bilinear form
$\la\cd,\cd\ra_{V_K}$ on $V_K$ by $\la a,b\ra_{V_K}:= \la
a,b\ra_{U'}$. In addition, since $W$ as a $V_K$-module is isomorphic
to $M_K(\beta,\gamma;\jmath|_K)$, we may view $\la \cd ,\cd
\ra_{U'}$ restricted on $W$ as a $V_K$-invariant bilinear form on
$M_K(\beta,\gamma;\jmath|_K)$. Then
$$
\begin{array}{ll}
  \la a, Y_X(u,z)v\ra_{V_K}
  &= \la a, \pi_{H} Y'(u,z)v \ra_{V_K} +\la a,\pi_{H+\alpha} Y'(u,z)v\ra_{V_K}
  \vsb\\
  &= \la a, \pi_{H} Y'(u,z)v \ra_{U'} +\la a,\pi_{H+\alpha} Y'(u,z)v\ra_{U'}
  \vsb\\
  &= \la a, Y'(u,z)v \ra_{U'}
  \vsb\\
  &= \la Y'(e^{zL(1)}(-z^{-2})^{L(0)}u,z^{-1}) a,v\ra_{U'} .
\end{array}
$$
By the equality above, it follows from Section 5.6 of \cite{FHL} and
\cite{Li2} that $Y_X(\cd,z)$ satisfies the Jacobi identity if and
only if we have a locality for any three elements in
$M_K(\beta,\gamma; \jmath|_K)$, which follows from \eqref{eq:5.2}
and \eqref{eq:5.3}. Therefore, $(X,Y_X(\cd,z))$ is also a framed
VOA. In fact, one can define a framed VOA structure on $V_E\oplus
M_E(\beta,\gamma; \jmath|_E)$ for any subcode $E$ of $H'$ containing
$H$ in a similar way. We shall deduce a contradiction from this
observation.

Let $W^1$ and $W^2$ be $V_H$-modules isomorphic to
$M_H(\beta,\gamma; \iota)$.
Since $W^i$ and $M_K(\beta,\gamma;\jmath|_K)$ are isomorphic
$V_H$-modules, we have $V_H$-isomorphisms
$\varphi_i : W^i \to M_K(\beta,\gamma;\jmath|_K) \subset X$
for $i=1,2$.
Set
$$
  X':=V_H\oplus V_{H+\alpha}\oplus W^1\oplus W^2.
$$
We shall define  a vertex operator  $Y_{X'}$ on $X'$ as follows.

For $a^0,b^0\in V_H$, $a^1,b^1\in V_{H+\alpha}$, $u^1,v^1\in W^1$ and
$u^2,v^2\in W^2$, define
$$
\begin{array}{l}
  Y_{X'}(a^0,z)
  :=
  \begin{bmatrix}
    Y_X(a^0,z) & 0 & 0 & 0
    \\
    0 & Y_X(a^0,z) & 0 & 0
    \\
    0 & 0 & \varphi_1^{-1} Y_X(a^0,z)\varphi_1 & 0
    \\
    0 & 0 & 0 & \varphi_2^{-1} Y_X(a^0,z)\varphi_2
  \end{bmatrix},
\end{array}
$$
$$
\begin{array}{l}
  Y_{X'}(a^1,z)
  :=
  \begin{bmatrix}
    0 & Y_X(a^1,z) & 0 & 0
    \\
    Y_X(a^1,z) & 0 & 0 & 0
    \\
    0 & 0 & 0 & \varphi_1^{-1} Y_X(a^1,z)\varphi_2
    \\
    0 & 0 & \varphi_2^{-1} Y_X(a^1,z)\varphi_1 & 0
  \end{bmatrix},
\end{array}
$$
$$
\begin{array}{l}
  Y_{X'}(u^1,z)
  \vsb\\
  :=
  \begin{bmatrix}
    0 & 0 & \pi_H Y_X(\varphi_1 u^1,z)\varphi_1 & 0
    \\
    0 & 0 & 0 & \pi_{H+\alpha} Y_X(\varphi_1 u^1,z) \varphi_2
    \\
    \varphi_1^{-1} Y_X(\varphi_1 u^1,z) & 0 & 0 & 0
    \\
    0 & \varphi_2^{-1} Y_X(\varphi_1 u^1,z) & 0 & 0
  \end{bmatrix},
\end{array}
$$
$$
\begin{array}{l}
  Y_{X'}(u^2,z)
  \vsb\\
  :=
  \begin{bmatrix}
    0 & 0 & 0 & \pi_H Y_X(\varphi_2 u^2,z)\varphi_2
    \\
    0 & 0 & \pi_{H+\alpha} Y_X(\varphi_2 u^2,z) \varphi_1 & 0
    \\
    0 & \varphi_1^{-1} Y_X(\varphi_2 u^2,z) & 0 & 0
    \\
    \varphi_2^{-1} Y_X(\varphi_2 u^2,z) & 0 & 0 & 0
  \end{bmatrix},
\end{array}
$$
on ${}^t[ b^0, b^1,v^1, v^2] \in V_H\oplus V_{H+\alpha}\oplus W^1 \oplus W^2$.
Note that $Y_X(\cd,z)$ is considered as a $V_H$-intertwining operator
and $\pi_H: X\to V_H$ and $\pi_{\alpha+H}: X\to V_{\alpha+H}$ are
$V_H$-projections on $V_H$ and $V_{\alpha+H}$, respectively.
By \eqref{eq:5.2} and \eqref{eq:5.3}, it is straightforward
to check that $Y_{X'}(\cd,z)$ satisfies the locality and hence
$(X',Y_{X'}(\cd,z))$ itself forms a VOA.
In fact, we have defined a VOA structure on $X'$ such that
$$
  V^1\cd W^1=W^2,
  \
  V^1\cd W^2=W^1,
  \
  W^1\cd W^1=W^2\cd W^2=V^0
  \q
  \mbox{and}
  \q
  W^1\cd W^2=V^2
$$
based on the framed VOA structure on $X$.

Now take a subspace
$$
  Z:=\{ a^0+b^0+ u^1+\varphi_2^{-1}\varphi_1 u^1 \in X' \mid
  a^0\in V_H, b^0\in V_{H+\alpha}, u^1\in W^1\} .
$$
Then it follows from the definition of $Y_{X'}(\cd,z)$ that $Z$ is
a subalgebra of $X'$ and the linear isomorphism
$$
  \psi: a^0+b^0+u^1+\varphi_2^{-1}\varphi_1 u^1
  \longmapsto a^0+b^0+\sqrt{2} \varphi_1 u^1,\q
  a^0\in V_H,\ b^0\in V_{H+\alpha},\ u^1\in W^1,
$$
defines a vertex operator algebra isomorphism between $Z$ and $X$.
Since $X$ is a framed VOA and every framed VOA is rational, $X'$ is
a completely reducible $Z$-module. However, since the quotient
$X'/Z$ has no 1/16-word component corresponding to a codeword 0, we
obtain $\psi^{-1}(M_K(\beta,\gamma;\jmath|_K))\cd (X'/Z)=0$ which is
a contradiction by Proposition 11.9 of \cite{DL}. This contradiction
comes from the assumption that $H\ne H'$. Hence, $H=H'$ as desired.
\qed \vsb

Now we present the main theorem of this paper:

\begin{theorem}\label{thm:5.5}
  Let $V=\oplus_{\alpha\in D} V^\alpha$ be a framed VOA with structure codes $(C,D)$.
  Then
  \\
  (1) For every non-zero $\alpha \in D$, the subcode $C_\alpha$ of $C$
  contains a doubly even self-dual subcode w.r.t.\ $\alpha$.
  \\
  (2) $C$ is even, every codeword of $D$ has a weight divisible by 8,
  and $D\subset C\subset D^\perp$.
\end{theorem}

\pf (1) follows from Proposition \ref{prop:5.4} since $V^0\oplus V^\alpha$ is a
framed sub VOA of $V$ for any non-zero $\alpha \in D$.
(2) follows from (1), since a self-dual subcode of $C_\alpha$ w.r.t.\ $\alpha$
always contains the codeword $\alpha$.
\qed
\vsb

As a corollary, we can also prove the following theorem.

\begin{theorem}\label{thm:5.6}
  Let $V=\oplus_{\alpha\in D} V^\alpha$ be a framed VOA with
  structure codes $(C,D)$.
  Then $V=\oplus_{\alpha \in D} V^\alpha$ is a
  $D$-graded simple current extension of the code VOA $V^0=V_C$.
\end{theorem}

\pf The assertion follows immediately from Theorem \ref{thm:5.5} and
Corollary \ref{cor:4.16}. \qed \vsb

There are many applications of Theorems \ref{thm:5.5} and \ref{thm:5.6}.

\begin{corollary}\label{cor:5.7}
  For a positive integer $n$, the number of isomorphism
  classes of framed VOAs with a fixed central charge $n/2$
  is finite.
\end{corollary}

\pf
By Theorem \ref{thm:5.6}, every framed VOA is a simple current
extension of a code VOA.
A code VOA is uniquely determined by its structure code
by Proposition \ref{prop:2.3}, and it has finitely many
irreducible representations as it is rational.
In particular, there are finitely many inequivalent simple current
modules over a code VOA.
Therefore, the number of isomorphism classes of framed VOAs of
given central charge is finite by the uniqueness of a simple
current extension in Proposition \ref{prop:2.3}.
\qed
\vsb

By Theorem \ref{thm:5.6}, we can immediately classify all irreducible
(both untwisted and $\Z_2$-twisted) modules over a framed VOA.

\begin{corollary}\label{cor:5.8}
  Let $V=\oplus_{\alpha\in D} V^\alpha$ be a framed VOA with
  structure codes $(C,D)$.
  Let $W$ be an irreducible $V^0$-module.
  Then there exists $\eta\in \Z_2^n$, which is unique modulo $D^\perp$,
  such that $W$ can be uniquely extended to an irreducible
  $\tau_\eta$-twisted $V$-module which is given by
  $V \fusion_{V^0} W$ as a $V^0$-module.
  In particular, every irreducible untwisted $V$-module is $D$-stable.
\end{corollary}

\pf
Let $\beta\in C^\perp$ be the 1/16-word of $W$.
Since all $V^\alpha$, $\alpha\in D$, are simple current
$V^0$-submodules, the fusion product $W^\alpha:=V^\alpha \fusion_{V^0} W$
is again irreducible.
It is clear that the binary 1/16-word of $W^\alpha$ is $\alpha+\beta$
so that all $W^\alpha$, $\alpha \in D$, are inequivalent $V^0$-modules.
Therefore, there exists a unique untwisted or $\Z_2$-twisted
$V$-module structure on
$\ind_{V^0}^V W =V\fusion_{V^0} W=\oplus_{\alpha\in D} W^\alpha$
by Theorem \ref{thm:2.4}.
Since any element in the dual group $D^*\simeq \Z_2^n/D^\perp$ is
realized as a Miyamoto involution $\tau_\eta$ associated to
a codeword $\eta\in \Z_2^n$, the induced module $\ind_{V^0}^V W$ is
indeed a $\tau_\eta$-twisted $V$-module.
\qed
\vsb

\begin{remark}\label{rem:5.9}
  By the corollary above and Proposition \ref{prop:2.5},
  we can compute the fusion rules of $V$-modules
  from those of $V_C$-modules,
\end{remark}

\begin{corollary}\label{cor:5.10}
  (\cite{DGH,M3})
  A framed VOA $V$ with structure codes $(C,D)$ is holomorphic
  if and only if $C=D^\perp$.
\end{corollary}

\pf
That a framed VOA having a structure code $(D^\perp,D)$ is holomorphic
is proved in \cite{M3} by showing that every module contains
a vacuum-like vector (cf.\ \cite{Li1}).
The converse is also proved in \cite{DGH} by using modular forms.
Here we give another, rather representation-theoretical proof.
Let $V$ be a holomorphic framed VOA with structure codes $(C,D)$
and the 1/16-word decomposition $V=\oplus_{\alpha\in D}V^\alpha$.
Take any codeword $\delta\in D^\perp$.
By the previous corollary, a $V_C$-module $M_C(0,\delta)$ can
be uniquely extended to either an untwisted or $\Z_2$-twisted $V$-module.
As a $V^0$-module, it is given by an induced module
$$
  V\fusion_{V_C} M_C(0,\delta)
  = \bigoplus_{\alpha\in D} V^\alpha \fusion_{V_C} M_C(0,\delta).
$$
By Lemma \ref{lem:4.12}, the top weight of $V^\alpha$ and that of
$V^\alpha \fusion_{V_C} M_C(0,\delta)$ are congruent modulo $\Z$
for all $\alpha \in D$.
Therefore, the induced module $V\fusion_{V_C} M_C(0,\delta)$ is
an irreducible untwisted $V$-module and thus isomorphic to $V$ itself,
as $V$ is holomorphic.
Then by considering the 1/16-word decomposition we see that
$M_C(0,\delta)=V^0=M_C(0,0)$.
Therefore, $\delta\in C$ by Lemma \ref{lem:4.6}
and hence $D^\perp =C$.
\qed


\subsection{Construction of a framed VOA}

In \cite{M3,Y2}, certain constructions of a framed VOA are
discussed. Assume the following:
\vsb\\
(1) $(C,D)$ is a pair of even linear codes of $\Z_2^n$ such that
\vsb\\
\begin{tabular}{lp{290pt}}
  (1-i)
  &  $C\subset D^\perp$,
  \vsb\\
  (1-ii)
  & for each $\alpha \in D$, there is a subcode $E^\alpha\subset C_\alpha$
  such that $E^\alpha$ is a 
  direct sum of the [8,4,4]-Hamming codes.
  \vsb
\end{tabular}
\\
(2) $V^0$ is a code VOA associated to $C$.
\vsb\\
(3) $\{ V^\alpha \mid \alpha\in D\}$ is a set of irreducible $V^0$-modules
such that
\vsb\\
\begin{tabular}{lp{290pt}}
(3-i) &  the 1/16-word of $V^\alpha$ is $\alpha$,
\vsb\\
(3-ii) &  all $V^\alpha$, $\alpha\in D$, have integral top weights,
\vsb\\
(3-iii) & the fusion product $V^\alpha \fusion_{V^0}V^\beta$ contains
$V^{\alpha+\beta}$ for all $\alpha,\beta\in D$.
\end{tabular}
\vsb\\
Then it is shown in \cite{M3,Y2} that the space
$V:=\oplus_{\alpha\in D}V^\alpha$ forms a framed VOA with structure
codes $(C,D)$. Instead of the condition (1-ii), assume that
\vsb\\
  \indent
  (1-iii) for each $\alpha \in D$, $C_\alpha$ contains a doubly even
  self-dual subcode w.r.t.\ $\alpha$.
\vsb\\
Then we have already shown in Lemma \ref{lem:5.1} that
$V^0\oplus V^\alpha$ forms a framed VOA.
So by the extension property of simple current extensions in
Theorem \ref{thm:2.7}, we can again show that
$V=\oplus_{\alpha\in D}V^\alpha$ forms a framed VOA with structure
codes $(C,D)$ under the other conditions.
The key idea in \cite{M3,Y2} is to use a special symmetry of
the code VOA associated to the [8,4,4]-Hamming code to form a minimal
$\Z_2$-graded extension $V^0\oplus V^\alpha$.
Thanks to Lemma \ref{lem:5.1}, we can transcend this step without
the [8,4,4]-Hamming code.

\begin{theorem}\label{thm:5.11}
  With reference to the conditions (1)--(3) above, assume the condition
  (1-iii) instead of (1-ii).
  Then $V=\oplus_{\alpha\in D} V^\alpha$ forms a framed VOA with structure
  codes $(C,D)$.
\end{theorem}

\pf
Let $\{\alpha^i \mid 1\leq i\leq r\}$ be a linear basis of $D$
and set $D^{(j)}:= \Span_{\Z_2}\{ \alpha_j \mid 1\leq j\leq i\}$.
By induction on $i$, we show that the space
$V[i]=\oplus_{\alpha\in D^{(i)}} V^\alpha$ forms a framed VOA with
structure codes $(C,D^{(i)})$. The case $i=0$ is trivial, and the
case $i=1$ is done in Lemma \ref{lem:5.1}. Now assume that $V[i]$
forms a framed VOA for $i\geq 1$. Then we have two simple current
extensions $V[i]=\oplus_{\alpha\in D^{(i)}}V^\alpha$ and
$V^0\oplus V^{\alpha^{i+1}}$. Applying Theorem \ref{thm:2.7} to a
set $\{ V^\alpha\mid \alpha \in D^{(i+1)}\}$ of simple current
$V^0$-modules, we obtain a $D^{(i+1)}$-graded simple current
extension $V[i+1]=\oplus_{\alpha\in D^{(i+1)}} V^\alpha$ of $V^0$.
Repeating this, finally we shall obtain the desired framed VOA
structure on $V[r]=\oplus_{\alpha\in D}V^\alpha$. \qed \vsb

We can also generalize Theorem 7.4.9 of \cite{Y2} as follows:

\begin{theorem}\label{thm:5.12}
  Let $V=\oplus_{\alpha\in D} V^\alpha$ be a framed VOA with
  structure codes $(C,D)$.
  For any even subcode $E$ such that $C\subset E\subset D^\perp$,
  the space
  $$
    \ind_C^E V
    := \bigoplus_{\alpha \in D} \ind_{V_C}^{V_E} V^\alpha
    = \bigoplus_{\alpha\in D} V_E \fusion_{V_C} V^\alpha
  $$
  forms a framed VOA with structure codes $(E,D)$.
\end{theorem}

\pf The idea of the proof is almost the same as that of Theorem
7.4.9 of \cite{Y2}. Let $\{ \gamma^i +C \mid 1\leq i\leq r\}$ be a
transversal for $C$ in $E$. It is clear that $V_E=\oplus_{i=1}^r
V_{C+\gamma^i}$ is an $E/C$-graded simple current extension of $V_C$
by Proposition \ref{prop:2.3}. First, we show that $V^\alpha$ is
uniquely extended to an untwisted $V_E$-module. For this, it
suffices to show that $V_{C+\gamma^i}\fusion_{V_C} V^\alpha$, $1\leq
i\leq r$, are inequivalent $V_C$-module. Assume
$V_{C+\gamma}\fusion_{V_C} V^\alpha \simeq V_{C+\delta}\fusion_{V_C}
V^\alpha$. It follows from a given framed VOA structure and Theorem
\ref{thm:5.5} that $V^\alpha\fusion_{V_C}V^\alpha\simeq V^0\simeq
V_C$. Since the fusion product is commutative and associative, by
multiplying the both sides of
$V_{C+\gamma}\fusion_{V_C}V^\alpha\simeq V_{C+\delta}
\fusion_{V_C}V^\alpha$ by $V^\alpha$ with respect to the fusion
product, we have
$$
  V_{C+\gamma}
  \simeq V_{C+\gamma}\fusion_{V_C} V^\alpha \fusion_{V_C} V^\alpha
  \simeq V_{C+\delta}\fusion_{V_C} V^\alpha \fusion_{V_C} V^\alpha
  \simeq V_{C+\delta} .
$$
Thus, $\gamma\equiv \delta \mod C$ and hence all
$V_{C+\gamma^i}\fusion_{V_C} V^\alpha$, $1\leq i\leq r$,
are inequivalent $V_C$-modules.
Since $E\subset D^\perp$, the top weight of $V^\alpha$ and that of
$V_{C+\gamma^i}\fusion_{V_C} V^\alpha$ are congruent modulo $\Z$
by Lemma \ref{lem:4.12}.
Therefore, $V^\alpha$ is uniquely extended to an irreducible untwisted
$V_E$-module $\ind_{V_C}^{V_E} V^\alpha = V_E\fusion_{V_C} V^\alpha$
by Theorem \ref{thm:2.4}.
Since all $\ind_{V_C}^{V_E} V^\alpha$, $\alpha\in D$, are $E/C$-stable
$V_E$-modules, we have the fusion rule
$$
  \l(\ind_{V_C}^{V_E} V^\alpha \r)\fusion_{V_E} \l(\ind_{V_C}^{V_E} V^\beta\r)
  \simeq \ind_{V_C}^{V_E} \l( V^\alpha\fusion_{V_C}V^\beta\r)
  \simeq \ind_{V_C}^{V_E} V^{\alpha+\beta}
$$
by Proposition \ref{prop:2.5}.
Therefore, the space
$$
  \ind_{V_C}^{V_E} V=\bigoplus_{\alpha\in D} \ind_{V_C}^{V_E} V^\alpha
$$
forms a framed VOA with structure codes $(E,D)$ by Theorem \ref{thm:5.11}.
\qed
\vsb

By Theorem \ref{thm:5.5}, its corollaries and Theorem
\ref{thm:5.11}, a pair of structure codes $(C,C^\perp)$ of a
holomorphic framed VOA satisfies the following conditions.

\begin{conditions}\label{conditions:1}
($F$-admissible condition) ~~\vsb\\
\begin{tabular}{lp{300pt}}
  (1)& The length of $C$ is divisible by 16.
  \vspace{0.1cm}\\
  (2)& $C$ is even, every codeword of $C^\perp$ has a weight divisible by 8,
  and\\ & $C^\perp\subset C$.
  \vspace{0.1cm}\\
  (3)& For any $\alpha \in C^\perp$, the subcode $C_\alpha$ of $C$ contains
  a doubly even self-dual subcode w.r.t.\ $\alpha$.
\end{tabular}
\end{conditions}

For simplicity, we will call a code $C$ \textit{$F$-admissible} if it
satisfies Condition \ref{conditions:1}.
Indeed, we can construct a holomorphic framed VOA starting from an
$F$-admissible code.

\begin{remark}\label{rem:5.13}
  A linear code $C$ is $F$-admissible if and only if
  its dual $C^\perp$ satisfies the following three conditions:
  \vsb\\
  \begin{tabular}{lp{405pt}}
    (i) the length of $C^\perp$ is divisible by 16,
    \vsb\\
    (ii) $C^\perp$ contains the all-one vector,
    \vsb\\
    (iii) $C^\perp$ is {\it triply even},
    that is, $\wt(\alpha)$ is divisible by 8 for any $\alpha \in C^\perp$.
  \end{tabular}
  \vsb\\
  For, let $D$ satisfy the conditions (i), (ii) and (iii) above.
  Then for any $\alpha,\beta \in D$, the weight of their intersection
  $\alpha\cd \beta$ is divisible by 4 and so $\alpha\cd D$ is doubly even.
  Then there exists a doubly even code $E$ containing $\alpha\cd D$
  such that $E$ is self-dual w.r.t.\ $\alpha$ by Theorem \ref{thm:5.3}.
  For any $\delta\in (\alpha\cd D)^{\perp_\alpha}$, we have
  $\la \delta,D\ra=\la \delta\cd \alpha,D\ra=\la \delta,\alpha\cd D\ra = 0$,
  showing $E \subset (\alpha\cd D)^{\perp_\alpha}\subset (D^\perp)_\alpha$.
  Therefore, $D^\perp$ is $F$-admissible.
\end{remark}

Let $C$ be an $F$-admissible code. Then the all-one vector
$\mathbf{1}=(11\dots 1)$ is contained in $C^\perp$. Since
$n=\wt(\mathbf{1})$ is divisible by $16$, all irreducible
$V_C$-modules with the $1/16$-word $\mathbf{1}$ have integral top
weights.  Let $V^\mathbf{1}$ be an irreducible $V_C$-modules with
the $1/16$-word $\mathbf{1}$.  Then $V^{\mathbf{1}}$ is a self-dual
simple current $V_C$-module.

\begin{lemma}\label{lem:5.14}
  Let $C$ be an $F$-admissible code.
  For $\alpha \in C^\perp$, $\alpha\not\in \{ 0, \mathbf{1}\}$,
 there exists an irreducible $V_C$-module $W^\al$
  such that $\tau(W^\al)=\al$ and both $W^\al$ and the fusion product
  $V^{\mathbf{1}}\fusion_{V_C} W^\al$ have
  integral top weights.
\end{lemma}

\pf Clearly, we can find
 an irreducible $V_C$-module $X$ such that $\tau(X)=\alpha$
and $X$ has an integral top weight. The $1/16$-word of the fusion
product $ V^{\mathbf{1}}\fusion_{V_C} X$ is then $\mathbf{1}+\alpha$
and its weight is divisible by 8. Thus, the top weight of
$V^{\mathbf{1}}\fusion_{V_C} X$ is in either $\Z$ or $\Z+1/2$. In
the former case, we just set $W^\alpha=X$. If the top weight is in
$\Z+1/2$, we take a codeword $\delta\in (\Z_2^n)_\alpha$ such that
$\delta$ is odd. By Lemma \ref{lem:4.12}, the $V_C$-module
$(V^{\mathbf{1}}\fusion_{V_C} X)\fusion_{V_C}M_C(0,\delta)$ has an
integral top weight. Set $W^\al= X\fusion_{V_C} M_C(0, \delta)$.
Since $\supp(\delta)\subset \supp(\alpha)$, the top weight of $X$ is
congruent to that of $W^\al$ modulo $\Z$ by Lemma \ref{lem:4.12},
which is integral. Moreover, $V^{\mathbf{1}} \fusion_{V_C} W^\alpha
\simeq V^{\mathbf{1}} \fusion_{V_C} (X \fusion_{V_C}M_C(0,\delta))
\simeq (V^{\mathbf{1}} \fusion_{V_C}X)\fusion_{V_C}M_C(0,\delta)$
also has an integral top weight as desired. \qed \vsb

\begin{proposition}\label{prop:5.15}
  Let $C$ be an $F$-admissible code and $D$ a proper subcode
  of $C^\perp$ containing $\mathbf{1}$.
  Suppose that we have a framed VOA $V = \oplus_{\alpha\in D} V^\alpha$
  with structure codes $(C,D)$.
  Then for $\beta\in C^\perp \setminus D$, there exist a self-dual
  simple current $V$-module $W$ such that $W$ has the 1/16-word
  decomposition $W=\oplus_{\alpha\in D} W^{\alpha+\beta}$ and
  $\tilde{V}=V\oplus W$ forms a framed VOA with structure codes
  $(C,D+\la \beta\ra)$.
\end{proposition}

\pf
By the previous lemma, we can take an irreducible $V_C$-module
$W^\beta$ such that $\tau(W^\be)=\be$ and both $W^\beta$ and
$V^{\mathbf{1}}\fusion_{V_C} W^\beta$ are of integral weights.
Since $\beta\in C^\perp$, it follows from
(3) of Condition \ref{conditions:1} that $W^\beta$ is a self-dual
simple current $V_C$-module. Then the induced module $\ind_{V^0}^V
W^\beta=\oplus_{\alpha \in D} V^\alpha\fusion_{V_C} W^\beta$ is an
irreducible $\tau_\eta$-twisted $V$-module for some $\eta\in \Z_2^n$
by Corollary \ref{cor:5.8}.
If $\eta\in D^\perp$, then the space $V \oplus
\ind_{V^0}^V W^\beta$ forms a framed VOA with structure codes
$(C,D+\la \beta\ra)$ by Theorem \ref{thm:5.11}.

If $\eta\not\in D^\perp$, then $\tau_\eta$ is not trivial.
Set $D^+:=\{ \alpha\in D \mid \la \alpha,\eta\ra=0\}$ and
$D^-:=\{ \alpha \in D \mid \la \alpha,\eta\ra=1\}$.
Then $D=D^+\sqcup D^-$ and $D^\pm \ne \emptyset$.
By our choice of $W^\beta$, the all-one vector $\mathbf{1}$ is
in $D^+$ so that $\eta$ is an even codeword.
We set
$$
  V^\pm :=\bigoplus_{\alpha \in D^\pm} V^\alpha,\q
  W^\pm :=\bigoplus_{\alpha \in D^\pm} V^\alpha \fusion_{V_C} W^\beta.
$$
Then all $V^\pm$, $W^\pm$ are irreducible $V^+$-modules.
The top weight of $W^+$ is integral but the top weight of $W^-$ is
in $\Z+1/2$.
We shall deform $W^-$ so that it has an integral top weight, also.

Since $[(D^+)^\perp \cap \la \beta \ra^\perp,D^\perp\cap \la
\beta\ra^\perp]=2$, there exists a codeword $\gamma \in (D^+)^\perp
\cap \la \beta\ra^\perp$ such that $\la \gamma, D^-\ra=1 \mod 2$.
Then it follows from Corollary \ref{cor:5.8} that
$$
  \tilde{W}^\pm:=V_{C+\gamma}\fusion_{V_C} W^\pm
$$
are irreducible untwisted $V^+$-modules.
Moreover, by our choice of $\gamma$, both of $\tilde{W}^\pm$ have
integral top weights since the top weight of $\tilde{W}^+$ is
congruent to $\la \gamma,\gamma+\beta\ra/2$ modulo $\Z$,
whereas the top weight of $\tilde{W}^-$ is congruent to
$\la \gamma, \beta+\gamma+D^-\ra/2+1/2$ modulo $\Z$.
Therefore, by Theorem \ref{thm:5.11}, we have a framed VOA
structure on
$$
  \tilde{V}:=V^+\oplus V^- \oplus \tilde{W}^+ \oplus \tilde{W}^-
$$
with structure codes $(C,D+\la \beta\ra)$.
Now setting $W:= \tilde{W}^+\oplus \tilde{W}^-$, we have the
desired extension of $V$.
This completes the proof.
\qed

\begin{remark}\label{rem:5.16}
  In the proof above, we can construct another extension of $V^+$ which also
  has the structure codes $(C,D+\la \beta\ra)$ in the following way.
  Take a codeword $\gamma \in D^\perp$ with $\la \gamma, \beta\ra =1$,
  which is  possible as $[D^\perp: D^\perp\cap \la \beta\ra^\perp]=2$,
  and set
  $\tilde{V}^- = V_{C+\gamma}\fusion_{V_C} V^-$ and
  $\tilde{W}^- = V_{C+\gamma}\fusion_{V_C} W^-$.
  Then one can similarly verify that the space
  $V^+\oplus \tilde{V}^-\oplus W^+ \oplus \tilde{W}^-$
  also forms a framed VOA with structure codes $(C,D+\la\beta\ra)$.
\end{remark}

\begin{theorem}\label{thm:5.17}
  There exists a holomorphic framed VOA with structure codes $(C,C^\perp)$
  if and only if $C$ is $F$-admissible, i.e., $C$ satisfies Condition
  \ref{conditions:1}.
\end{theorem}

\pf Let $\{ \alpha_1,\dots,\alpha_r\}$ be a linear basis of
$C^\perp$ with $\alpha_1 = \mathbf{1}$ and set $D[i]:=\Span_{\Z_2}\{
\alpha_j \mid 1\leq j\leq i\}$ for $1\leq i\leq r$. By Lemma
\ref{lem:5.1} we can construct a framed $V[1]:=V_C\oplus
V^\mathbf{1}$ with structure codes $(C,D[1])$. By Proposition
\ref{prop:5.15}, we can construct a framed VOA $V[2]$ with structure
codes $(C,D[2])$ which is a $\Z_2$-graded simple current extension
of $V[1]$. Recursively, we can construct a $\Z_2$-graded simple
current extension $V[i+1]$ of $V[i]$ which has structure codes
$(C,D[i+1])$ and we shall obtain a holomorphic framed VOA $V[r]$
with structure codes $(C,D[r])=(C,C^\perp)$. \qed

\begin{remark}\label{rem:5.18}
  Condition \ref{conditions:1}, especially (2) and (3), give quite strong
  restrictions on a code $C$.
  Roughly speaking, $C$ must be much bigger than its dual $C^\perp$
  by (3) of Condition \ref{conditions:1}.
  In addition, if we assume that the minimum weight of $C$ is greater than 2,
  then the corresponding framed VOA may have a finite full automorphism group
  (cf.\ \cite[Corollary 3.9]{LSY}).
  It seems possible to classify all $F$-admissible codes $C$ if the length
  is small.
  It suggests a possibility for classifying all holomorphic framed
  VOAs of small central charge.
  The most interesting (and the first non-trivial) case would be the
  classification of $c=24$ holomorphic framed VOAs.
  In fact, one can prove that the moonshine vertex operator algebra
  $V^\natural$ is the unique holomorphic framed VOA of central charge 24
  whose weight one subspace is trivial, which is a variant of the famous
  uniqueness conjecture of the moonshine vertex operator algebra proposed
  in \cite{FLM} (see also \cite{DGL}).
  The key point is that the structure codes of $V^\natural$
  (or any holomorphic framed VOA $V$ of central charge 24 and $V_1=0$) are
  closely related to those of the Leech lattice VOA $V_\Lambda$.
  If $V=\oplus_{\al\in C^\perp} V^\al$ with $V^0\simeq V_C$ is
  a holomorphic framed VOA of central charge 24 and $V_1=0$, then the minimal
  weight of $C$ is greater than or equal to $4$.
  In this case, for any $\delta\in \Z_2^{48}$ of weight $2$,
  the $\tau_\delta$-twisted orbifold construction yields a VOA
  \[
    V(\tau_\delta)= \bigoplus_{\al\in D} \left( V^\al \oplus
    V_{C+\delta}\fusion_{V_C} V^\al \right),
    \quad D=\{\al\in  C^\perp \mid  \la \al, \delta\ra=0\},
  \]
  which is isomorphic to the Leech lattice $V_\Lambda$ and a pair
  $(C\sqcup (\delta+C), D)$ will be the structure codes of $V_\Lambda$.
  Note that the weight one subspace of $V(\tau_\delta)$ forms an abelian
  Lie algebra with respect to the bracket $[a,b]=a_{(0)} b$.
  We shall give more details on this point in our next work \cite{LY}.
\end{remark}

\section{Frame stabilizers and order four symmetries}

In Section 5, we have seen that structure codes $(C,D)$ of a framed
VOA $V=\oplus_{\alpha\in D}V^\alpha$ satisfy certain duality
conditions. The main property is that for any $\alpha \in D$, the
subcode $C_\alpha$ contains a doubly even self-dual subcode w.r.t.\
$\alpha$ and $V^\alpha$ is a simple current $V^0$-module. However,
it is shown in Corollary \ref{cor:4.16} that $V^\alpha$ is a simple
current module without the assumption on the doubly even property.
In this section, we shall discuss the role of the doubly even
property. It turns out that by relaxing the doubly even property, we
can obtain a refinement of the 1/16-word decomposition and define an
automorphism of order four in the pointwise frame stabilizer.

We begin by defining the frame stabilizer and the pointwise
frame stabilizer of a framed VOA.

\begin{definition}\label{df:6.1}
  Let $V$ be a framed VOA with a frame $F=\vir(e^1)\tensor \cds \tensor \vir(e^n)$.
  The {\it frame stabilizer} of $F$ is the subgroup of all automorphisms of $V$
  which stabilizes the frame $F$ setwise.
  The {\it pointwise frame stabilizer} is the subgroup of $\aut(V)$ which
  fixes $F$ pointwise.
  The frame stabilizer and the pointwise frame stabilizer of $F$ are denoted
  by $\stab_V (F)$ and $\pstab_V(F)$, respectively.
\end{definition}

Let $(C,D)$ be the structure code of $V$ with respect to $F$, i.e.,
$$
  V=\bigoplus_{\alpha\in D} V^\alpha,\quad  \tau(V^\alpha)=\alpha \quad
  \text{ and}\quad V^0= V_C.
$$
For any $\theta\in \pstab_V(F)$, it is easy to see that
$\tau_{e^i}=\tau_{\theta e^i}=\theta \tau_{e^i} \theta^{-1}$ and
thus $\theta$ centralizes $\la \tau_{e^1},\dots,\tau_{e^n}\ra$.
Therefore, the group $\tau(\Z_2^n)=\la \tau_{e^1},\dots,\tau_{e^n}\ra$
generated by the $\tau$-involutions is a central subgroup of $\pstab_V(F)$
isomorphic to $\Z_2^n/D^\perp$.
In addition, we have $\theta V^\alpha =V^\alpha$ for all $\alpha
\in D$ and hence $\theta|_{V^0}$ is  an automorphism of $V^0$.

The following results can be proved easily using the fusion rules.

\begin{lemma}\label{lem:6.2}
Let $V=\oplus_{\alpha \in D} V^\alpha$ be a framed VOA.
\\
(1) Let $\phi \in \aut (V^0)$ such that $\phi|_{F}=\mathrm{id}_{F}$.
Then $\phi \in \sigma(\Z_2^n)=\la \sigma_{e^1},\dots,\sigma_{e^n}\ra$.
\\
(2) Let $g\in \aut (V)$ such that $g|_{V^0}=\mathrm{id}_{V^0}$.
Then $g\in \tau(\Z_2^n)=\la \tau_{e^1},\dots,\tau_{e^n}\ra$.
\end{lemma}

\pf  (1) Consider the 1/2-word decomposition
$$
  V^0=\bigoplus_{\beta=(\beta_1,\dots,\beta_n)\in C} L(\shf,\beta_1/2)\tensor
  \cds \tensor L(\shf,\beta_n/2)
$$
of $V^0$ as an $F=\vir(e^1)\tensor \cds \tensor\vir(e^n)$-module.
Since $\phi|_{F}=\mathrm{id}_F$, it follows from Schur's lemma that
$\phi|_{V^0}$ acts on $L(\shf,\beta_1/2)\tensor \cds \tensor
L(\shf,\beta_n/2)$ by a non-zero scalar $a_\alpha$ for each $\alpha
\in C$. Moreover, it follows from the fusion rules of
$L(\shf,0)$-modules in \eqref{eq:3.1} that $a_\alpha a_\beta
=a_{\alpha+\beta}$ for all $\alpha,\beta\in C$. Thus the association
$C\ni \alpha\mapsto a_\alpha\in \C$ defines a character of $C$ and
hence there is a codeword $\xi\in \Z_2^n$ such that
$a_\alpha=(-1)^{\la \xi,\alpha\ra}$. Now it is easy to see that
$\phi|_{V^0}$ is realizable as a product of $\sigma_{e^i}$, $1\leq
e^i\leq n$, that is, $\phi|_{V^0}=\sigma_\xi$.

\medskip

(2) Since each $V^\alpha$, $\alpha\in D$, is an irreducible
$V^0$-module, it follows from Schur's lemma that $g$ acts on
$V^\alpha$ by a non-zero scalar $t_\alpha\in \C$.
Then again by the fusion rules of $L(\shf,0)$-modules
in \eqref{eq:3.1} we have $t_\alpha t_\beta=t_{\alpha+\beta}$
so that the map $\al \mapsto t_\alpha$ defines a character of $D$.
Therefore, there exists a codeword $\eta \in \Z_2^n$ such that
$g=\tau_\eta\in \la \tau_{e^1},\dots,\tau_{e^n}\ra$.
\qed
\vsb

As a corollary, we have the following theorem.

\begin{theorem}\label{thm:6.3}
  Let $V$ be a framed VOA with a frame
  $F=\vir(e^1)\tensor \cds \tensor \vir(e^n)$.
  For any $\theta\in \pstab_V(F)$, there exist
  $\xi$ and $\eta \in \Z_2^n$ such that
  $$
    \theta|_{V^0}= \sigma_\xi \qquad \text{ and } \qquad \theta^2 =\tau_\eta.
  $$
  In particular, we have $\theta^4=1$.
\end{theorem}
\vsb

Let $\theta\in \pstab_V(F)$.
Then $ \theta|_{V^0}=\sigma_\xi$ for some $\xi \in \Z_2^n$.
That means $\theta$ is an extension of a $\sigma$-involution on
$V^0$ to the whole framed VOA $V$.
In this section, we shall give a necessary and sufficient condition
on whether a $\sigma$-involution $\sigma_\xi$ can be extended to the
whole $V$.
Our argument is based on the representation theory of code VOAs
developed in Section 4 and 5.
\vsb

First let us consider $\theta\in \pstab_V(F)$ such that
$\theta|_{V^0}=\sigma_\xi \neq \mathrm{id}_{V^0}$, i.e.,
$\xi\notin C^\perp$.
Set
$C^0:=\{ \alpha \in C \mid \la \xi,\alpha\ra =0\}$ and
$C^1:=\{ \alpha\in C \mid \la \xi,\alpha \ra =1\}$.
Then $C^0$ is a subcode of $C$, $[C:C^0]=2$ and $C=C^0\sqcup C^1$.
Note also that $V_{C^0}$ is fixed by $\theta$ and $ \theta$ acts
by $-1$ on $V_{C^1}$.
In other words, $V^0=V_{C^0}\oplus V_{C^1}$ is the eigenspace
decomposition of $\theta$ on $V^0$.

Now assume that $\theta^2=\tau_\eta$ for some $\eta\in \Z_2^n$.
For each non-zero $\alpha \in D$, it is clear that $V^0\oplus V^\alpha$
is a subalgebra of $V$ and $\theta$ stabilizes $V^0\oplus V^\alpha$.
If $\alpha \in D\cap \la \eta\ra^\perp$, then $\theta^2$ acts as
an identity on $V^0\oplus V^\alpha$ and the eigenvalues of $\theta$
on $V^\al$ are $\pm 1$.
Let $V^{\alpha+}$ and $V^{\alpha-}$ be the eigenspaces of $\theta$
with eigenvalues $+1$ and $-1$, respectively.
Note that both $V^{\alpha\pm}$ are non-zero inequivalent
irreducible $V^{0+}$-submodules.
For if $V^{\alpha+}=0$, then $V^{0-}\cd V^{\alpha-}=V^{\alpha+}=0$, which
contradicts Proposition 11.9 of \cite{DL}.
Since the subalgebra
$V^0\oplus V^\alpha = V^{0+}\oplus V^{0-} \oplus V^{\alpha+} \oplus V^{\alpha-}$
affords a faithful action of a group $\Z_2\times \Z_2$ of order 4,
$V^{\alpha\pm}$ are inequivalent irreducible $V^{0+}$-submodules
by the quantum Galois theory \cite{DM1}.

If $\alpha \in D\setminus \la \eta\ra^\perp$,
then $\theta^2=-1$ on $V^\al$ and the eigenvalues of $\theta$
on $V^\alpha$ are $\pm \sqrt{-1}$.
Let $V^{\alpha\pm}$  be the eigenspace of $\theta$ of
eigenvalues of $\pm \sqrt{-1}$ on $V^\alpha$.
Then $V^\alpha=V^{\alpha+}\oplus V^{\alpha-}$.
$V^{\alpha\pm}$ are again non-zero inequivalent irreducible
$V^{0+}$-submodules.
The argument above actually shows that ${\theta}$ is
of order $2$ if and only if $\la \eta, D\ra =0$, i.e.,
$\eta\in D^\perp$.
\vsb

By the observation above, we have

\begin{lemma}\label{lem:6.4}
  For any $\alpha \in D$, let $V^{\alpha\pm}$ be defined as above.
  Then the dual $V^{0+}$-module $(V^{\alpha\pm})^*$ is
  isomorphic to $V^{\alpha\pm}$ if and only if $\alpha \in \la \eta\ra^\perp$.
  Otherwise,  $(V^{\alpha\pm})^*$ is isomorphic to $V^{\alpha\mp}$.
\end{lemma}

\pf
Since any framed VOA is self-dual, the sub VOA $V^0\oplus V^\alpha$ of $V$ is
also self-dual.
Since $V^{\alpha\pm}\cd V^{\alpha\pm}=V^{0+}$ if and only if
$\alpha\in \la \eta \ra^{\perp}$, the duality is as in the assertion.
\qed
\vsb

We have shown that for any $\theta\in \pstab_V(F)\setminus \tau(\Z_2^n)$,
$\abs{\theta}=2$ if and only if all irreducible $V^{0+}$-submodules
of $V$ are self-dual, and otherwise $\abs{\theta}=4$.
We rewrite this condition in terms of the structure codes as follows.

\begin{lemma}\label{lem:6.5}
  Let $\theta\in \pstab_V(F)$ such that
  $\theta|_{V^0}=\sigma_\xi$
  and $\theta^2=\tau_\eta$ for some $\xi\in \Z_2^n \setminus C^\perp$
  and $\eta\in \Z_2^n$.
  Set $C^0=\{ \alpha \in C \mid \la \xi,\alpha\ra =0\}$ and
  $C^1=\{ \alpha \in C \mid \la \xi,\alpha\ra =1\}$.
  \\
  (1) For $\alpha\in D\cap \la \eta\ra^\perp$, $(C^0)_\alpha$ contains
  a doubly even self-dual subcode w.r.t.\ $\alpha$.
  \\
  (2) For $\alpha\in D\setminus \la \eta\ra^\perp$,
  $(C^0)_\alpha$ contains a self-dual subcode w.r.t.\ $\alpha$,
  but $(C^0)_\alpha$ does not contain any doubly even self-dual
  subcode w.r.t.\ $\alpha$.
\end{lemma}

\pf
(1) For $\alpha \in D\cap \la \eta\ra^\perp$, let
$V^\alpha=V^{\alpha+}\oplus V^{\alpha-}$ be the eigenspace
decomposition such that $\theta$ acts on $V^{\alpha\pm}$ by $\pm 1$.
In this case the subspace $V^{0+}\oplus V^{\alpha+}$ forms a framed
sub VOA of $V$.
By Proposition \ref{prop:5.4}, $(C^0)_\alpha$ contains a doubly even
self-dual subcode w.r.t.\ $\alpha$.

(2) For $\alpha \in D\setminus \la \eta\ra^\perp$, let
$V^\alpha=V^{\alpha+}\oplus V^{\alpha-}$ be the eigenspace
decomposition such that $\theta$ acts on $V^{\alpha\pm}$
by $\pm \sqrt{-1}$.
In this case the restriction of $\theta$ on the sub VOA
$$
  V^0\oplus V^\alpha=V^{0+}\oplus V^{0-}\oplus V^{\alpha+} \oplus V^{\alpha-}
$$
is of order $4$. By the quantum Galois theory \cite{DM1},
$V^{\alpha+}$ and $V^{\alpha-}$ are inequivalent irreducible
$V^{0+}=V_{C^0}$-modules. By Lemma \ref{lem:6.4}, $V^{\alpha+}$ and
$V^{\alpha-}$ are dual to each other. Therefore, by Proposition
\ref{prop:4.9}, any maximal self-orthogonal subcode of
$(C^0)_\alpha$ is not doubly even. Let $H$ be a doubly even
self-dual subcode of $C_\alpha$ w.r.t.\ $\alpha$ and $H^0$ a maximal
self-orthogonal subcode of $(C^0)_\alpha$. Since $(C^0)_\alpha$ does
not contain a doubly even self-dual subcode w.r.t.\ $\alpha$,
$(C^0)_\alpha$ is a proper subgroup of $C_\alpha$ so that
$[C_\alpha:(C^0)_\alpha]=2$. It follows from Theorem \ref{thm:4.3}
that $V^\alpha$ is a direct sum of $[C:C_\alpha]$ inequivalent
$V(0)$-submodules with the  multiplicity $[C_\alpha:H]$. Similarly,
each of $V^{\alpha\pm}$ is a direct sum of $[C^0:(C^0)_\alpha]$
inequivalent irreducible $F$-submodules with the multiplicity
$[(C^0)_\alpha:H^0]$. Since $V^{\alpha+}$ and $V^{\alpha-}$ are dual
to each other, they are isomorphic as $F$-modules. Therefore, by
counting multiplicity of irreducible $F$-submodules of $V^\alpha$
and $V^{\alpha\pm}$, one has $[C_\alpha:H]=2[(C^0)_\alpha:H^0]$.
Combining with $\abs{C_\alpha}=2\abs{(C^0)_\alpha}$, we obtain
$\abs{H}=\abs{H^0}=2^{\wt(\alpha)/2}$. Therefore, $H^0$ is a
self-dual subcode of $(C^0)_\alpha$ w.r.t.\ $\alpha$. \qed \vsb

\begin{lemma}\label{lem:6.6}
  Let $C$ be an even code and $\beta \in C^\perp$.
  Assume that $V=V_C\oplus M_C(\beta,\gamma;\iota)$ forms a framed VOA
  and $C$ contains a subcode $E$ with index two such that
  $E$ contains a self-dual subcode w.r.t.\ $\beta$.
  Then $V$ decomposes into a direct sum of four inequivalent simple current
  $V_E$-submodules
  $$
    V=M_E(0,0) \oplus M_E(0,\delta)
    \oplus M_E(\beta,\gamma;\jmath)\oplus M_E(\beta,\gamma+\delta;\jmath),
  $$
  where $\delta$ is an element of $C$ such that
  $C=E\sqcup (E+\delta)$ and $\jmath:E\cap E^\perp\to \C^*$
  is a map such that $\jmath|_{E\cap C^\perp}=\iota$ and
  $(\al,\jmath(\al))\cdot (\be,\jmath(\be))
  =(\al+\be,\jmath(\al+\be))\in \pi_{\C^*}^{-1}(E)$
  for all $\al,\be\in E\cap E^\perp$.
  Moreover, all irreducible $V_E$-submodules of $V$ are self-dual
  if and only if   $E$ contains a doubly even self-dual subcode
  w.r.t.\ $\beta$.
\end{lemma}

\pf Let $\delta\in C$ such that $C=E\sqcup (E+\delta)$. Then the
decomposition $V_C=M_E(0,0)\oplus M_E(0,\delta)$ is obvious. Let $H$
be a self dual subcode of $E_\be$ w.r.t. $\beta$. Then $H$ is still
a maximal self-orthogonal subcode of $C_\be$. Let $\jmath: H \to
\C^*$ be an extension of $\iota: C\cap C^\perp \to \C^*$ such that
$(\al,\jmath(\al))\cdot (\be,\jmath(\be))=(\al+\be,\jmath(\al+\be))$
for all $\alpha, \beta\in H$. Then $\jmath|_{E\cap C^\perp}=\iota$
and the decomposition
$M_C(\beta,\gamma;\iota)=M_E(\beta,\gamma;\jmath)\oplus
M_E(\beta,\gamma+\delta;\jmath)$ follows from Corollary
\ref{cor:4.7}.  That all irreducible $V_E$-submodules are simple
currents follows from Corollary \ref{cor:4.16}. If $E$ contains a
doubly even self-dual subcode w.r.t.\ $\beta$, then all irreducible
$V_E$-submodules of $V$ are self-dual by Proposition \ref{prop:4.9}.
Conversely, if all irreducible $V_E$-submodules of $V$ are
self-dual, then $V_E\oplus M_E(\beta,\gamma;\jmath)$ forms a sub VOA
of $V$ so that $E$ contains a doubly even self-dual subcode w.r.t.\
$\beta$ by Proposition \ref{prop:5.4}. This completes the proof.
\qed

\begin{theorem}\label{thm:6.7}
  Let $V$ be a framed VOA with a frame
  $F=\vir(e^1)\tensor \cds \tensor \vir(e^n)$ and let $(C,D)$ be
  the corresponding structure codes.
  For a codeword $\xi\in \Z_2^n \setminus C^\perp$, there exists
  $\theta\in \pstab_V(F)$ such that $\theta|_{V^0}=\sigma_\xi$
  if and only if $\alpha\cd \xi\in C$ for all $\alpha \in D$.
  Moreover, $\abs{\theta}=2$ if and only if
  $\wt (\alpha\cd \xi)\equiv 0 \mod 4$ for all   $\alpha \in D$,
  and otherwise $\abs{\theta}=4$.
\end{theorem}

\pf
Let $\theta\in \pstab_V(F)$ such that $\theta|_{V^0}=\sigma_\xi$
with $\xi \in \Z_2^n\setminus C^\perp$.
Set $C^0:=C\cap \la \xi\ra^\perp$ and $C^1:=C\setminus C^0$.
Then by Lemma \ref{lem:6.5}, $(C^0)_\alpha$ contains a self-dual
subcode w.r.t.\ $\alpha$ for any $\alpha\in D$.
Since
$$
  (C^0)_\alpha =\{ \beta\in C \mid \la \beta,\xi\ra=0 \ \text{ and }\
  \supp(\beta)\subset \supp(\alpha) \}
$$
and $\la \beta, \xi\ra = \la \beta,\alpha\cd \xi\ra =0$ for all
$\beta\in (C^0)_\alpha$, $\alpha\cd \xi\in ((C^0)_\alpha)^\perp$ for
all $\alpha \in D$. Therefore, $\alpha\cd \xi$ is contained in all
self-dual subcodes of $(C^0)_\al $ w.r.t.\ $\alpha$ and hence
$\alpha\cd \xi\in (C^0)_\alpha \subset C$ as claimed.

Conversely, assume that a codeword $\xi\in \Z_2^n\setminus
C^\perp$ satisfies $\alpha\cd \xi\in C$ for all $\alpha \in D$.
Then $C^0=C\cap \la \xi\ra^\perp$ is a proper subcode of $C$ with
index 2.
By definition, $\alpha\cd \xi \in ((C^0)_\alpha)^\perp$ for all
$\alpha\in D$.
Therefore, any maximal self-orthogonal subcode of $(C^0)_\alpha$
contains $\alpha\cd \xi$.
Set
$D^0:=\{\alpha \in D \mid \wt(\alpha\cd\xi)\equiv 0 \mod 4\}$
and
$D^1:=\{\alpha \in D \mid \wt(\alpha\cd\xi) \equiv 2 \mod 4\}$.
It is clear that $D=D^0\sqcup D^1$. If $\alpha \in D^0$,
then there exists a doubly even self-dual subcode of $(C^0)_\alpha$
w.r.t.\ $\alpha$.
For, let $H$ be a doubly even self-dual subcode of $C_\alpha$,
which exists by Proposition \ref{prop:5.4}.
If $H$ is contained in $(C^0)_\alpha$, then we are done.
If not, then $\alpha\cd \xi\not \in H$ and
$H\cap (C^0)_\alpha=H\cap \la \xi\ra^\perp$ is a subcode
of $H$ with index 2 so that
$$
  (H\cap \la \xi\ra^\perp )\sqcup (H\cap \la\xi\ra^\perp +\alpha\cd \xi)
$$
gives a doubly even self-dual subcode of $(C^0)_\alpha$ w.r.t.\
$\alpha$. Similarly, we can show that $(C^0)_\alpha$ contains a
self-dual subcode w.r.t.\ $\alpha$ for any $\alpha\in D^1$. But in
this case any self-dual subcode of $(C^0)_\alpha$ w.r.t.\ $\alpha$
is not doubly even, as it always contains $\alpha\cd \xi$. We have
shown that for each $\alpha\in D$, $(C^0)_\alpha$ contains a
self-dual subcode w.r.t.\ $\alpha$ so that one has a
$V_{C^0}$-module decomposition $V^\alpha=V^{\alpha,1}\oplus
V^{\alpha,2}$ and $V^{\alpha,p}$, $p=1,2$, are simple current
$V_{C^0}$-submodules by Lemma \ref{lem:6.6}.

Let $\{ \alpha^1,\dots,\alpha^r\}$ be a linear basis of $D^0$.
For each $i$, $1\leq i\leq r$, choose an irreducible $V_{C^0}$-submodule
$U^{\alpha^i}$ of $V^{\alpha^i}$ arbitrary.
Then for $\alpha= \alpha^{i_1}+\cds +\alpha^{i_k} \in D^0$, set
\begin{equation}\label{eq:6.1}
  U^\alpha:= U^{\alpha^{i_1}}\fusion_{V_{C^0}} \cds \fusion_{V_{C^0}} U^{\alpha^{i_k}}.
\end{equation}
Since all $U^{\alpha^i}$, $1\leq i\leq r$, are simple current
self-dual $V_{C^0}$-modules, $U^\alpha$ is uniquely defined by
\eqref{eq:6.1} for all $\alpha \in D^0$.
Note that $U^0= V_{C^0}$.
Since $\oplus_{\alpha \in D^0} V^\alpha$ is a sub VOA of $V$,
$U^\alpha$ are irreducible $V_{C^0}$-submodules of $V^\alpha$ for
all $\alpha \in D^0$.
Therefore, we obtain a framed sub VOA
$U:=\oplus_{\alpha \in D^0} U^\alpha$ of $V$ with structure codes
$(C^0,D^0)$.
It is easy to see that
$V^\alpha = U^\alpha \oplus (V_{C^1}\fusion_{V_{C^0}} U^\alpha)$
for $\alpha \in D^0$ by Lemma \ref{lem:6.6}.

If $D=D^0$, then we have $V=U\oplus (V_{C^1}\fusion_{V_{C^0}} U)$
as a $V_{C^0}$-module.
In this case we define a linear automorphism $\theta_\xi$ on $V$ by
$$
  \theta_\xi:=
  \begin{cases}
    \ \ 1 & \text{ on }\ U,
    \vsb\\
    -1 & \text{ on }\ V_{C^1}\fusion_{V_{C^0}} U.
  \end{cases}
$$
Then it follows from Lemma \ref{lem:6.6} and Proposition \ref{prop:4.13}
that $\theta_\xi\in \pstab_V(F)$ and $\theta_\xi|_{V^0}=\sigma_\xi$.
Therefore, $\sigma_\xi$ can be extended to an involution on $V$.

If $D\ne D^0$, then $D=D^0\sqcup D^1$ with $D^1\ne \emptyset$.
In this case, take one $\beta \in D^1$ and an irreducible
$V_{C^0}$-submodule $W^\beta$ of $V^\beta$.
Since $W^\beta$ and all $U^\alpha$, $\alpha \in C^0$, are simple
current $V_{C^0}$-modules, we have a $V_{C^0}$-module decomposition
$$
  V^{\alpha+\beta}= (U^\alpha \fusion_{V_{C^0}} W^\beta) \oplus
  (V_{C^1}\fusion_{V_{C^0}} U^\alpha \fusion_{V_{C^0}} W^\beta)
$$
of $V^{\alpha+\beta}$ for all $\alpha \in C^0$ by Lemma \ref{lem:6.6}.
Since $(C^0)_{\alpha+\beta}$ contains no doubly even self-dual
subcode w.r.t.\ $\alpha+\beta$, the decomposition
$$
  V^0\oplus V^{\alpha+\beta}
  = V_{C^0}
    \oplus V_{C^1}
    \oplus (U^\alpha \fusion_{V_{C^0}} W^\beta)
    \oplus (V_{C^1}\fusion_{V_{C^0}} U^\alpha \fusion_{V_{C^0}} W^\beta)
$$
induces an order four automorphism on a sub VOA
$V^0\oplus V^{\alpha+\beta}$ of $V$ by Lemma \ref{lem:6.6} and
Proposition \ref{prop:4.13}.
Set
$$
  W:= \bigoplus_{\alpha\in C^0} U^\alpha \fusion_{V_{C^0}} W^\beta .
$$
Then we have obtained the following decomposition of $V$ as
a $V_{C^0}$-module:
$$
  V= U\oplus (V_{C^1}\fusion_{V_{C^0}} U) \oplus W
     \oplus (V_{C^1} \fusion_{V_{C^0}} W).
$$
We define a linear automorphism $\theta_\xi$ on $V$ by
$$
  \theta_\xi:=
  \begin{cases}
    \ \ 1 & \text{ on }\ U,
    \vsb\\
     -1  & \text{ on }\ V_{C^1}\fusion_{V_{C^0}}U,
    \vsb\\
    \sqrt{-1} & \text{ on }\ W,
    \vsb\\
    -\sqrt{-1} & \text{ on }\ V_{C^1}\fusion_{V_{C^0}} W.
  \end{cases}
$$
Then it follows from the argument above that $\theta_\xi\in
\pstab_V(F)$ and $\theta_\xi|_{V^0}= \sigma_\xi$.
Therefore, $\sigma_\xi$ gives rise to an automorphism of order 4.

Summarizing, we have shown that there exists
$\theta\in \pstab_V(F)$ such that $\theta|_{V^0}=\sigma_\xi$
if and only if $\alpha\cd \xi\in C$ for all $\alpha \in D$.
It remains to show that for such $\theta$, $\abs{\theta}=2$
if and only if $\wt(\alpha\cd \xi )\equiv 0 \mod 4$.
But this is almost obvious by the preceding argument.
\qed

\begin{remark}\label{rem:6.8}
  Let $V=\oplus_{\alpha\in D} V^\alpha$  be a framed VOA with
  structure codes $(C,D)$.
  It was conjectured in \cite[Conjecture 1]{M3} that for any codeword
  $\beta\in C$, the $\sigma$-type involution $\sigma_\be\in \aut(V^0)$
  can be extended to an automorphism of $V$, that means there
  exists $g\in \pstab_V(F)$ such that $g|_{V^0}=\sigma_\beta$.
  By the theorem above, we know that this is not correct;
  we have to take a codeword $\beta$ such that
  $\alpha\cd \beta\in C$ for all $\alpha \in D$.
\end{remark}

Motivated by Theorem \ref{thm:6.7}, we define
$P:=\{ \xi \in \Z_2^n\mid \alpha\cd \xi \in C\ \text{for all}\ \alpha\in D\}$.
It is clear that $P$ is a linear subcode of $C$.
Moreover, we have

\begin{lemma}\label{lem:6.9}
  $C^\perp \subset P$.
\end{lemma}

\pf
Let $\delta\in C^\perp$. For $\alpha \in D$, one has
$\la \delta, C_\alpha\ra =0$ by definition and hence
$\la \alpha\cd \delta, C_\alpha\ra =0$.
Since $C_\alpha$ contains a self-dual subcode w.r.t.\ $\alpha$
by Theorem \ref{thm:5.5},
$\alpha\cd \delta\in C_\alpha \subset C$.
\qed
\vsb

For each codeword $\xi\in P$, there exists
$\theta_\xi\in \pstab_V(F)$ such that $\theta_\xi|_{V^0}=\sigma_\xi$
by Theorem \ref{thm:6.7}.
However, the construction of $\theta_\xi$ in the proof
of Theorem \ref{thm:6.7} is not unique since we have to choose
a linear basis of $D^0$ and irreducible $V_{C^0}$-submodules.
Nevertheless, the following holds.

\begin{lemma}\label{lem:6.10}
  Let $\theta, \phi\in \pstab_V(F)$ such that
  $\theta|_{V^0}=\phi|_{V^0}=\sigma_{\xi}$.
  Then $\phi=\theta \tau_{\eta}$ for some $\eta\in \Z_2^n$.
\end{lemma}

\pf
Since $\theta|_{V^0}=\phi|_{V^0}$, we have
$\theta^{-1}\phi|_{V^0}=\mathrm{id}_{V^0}$.
By Lemma \ref{lem:6.2}, there exists $\eta\in \Z_2^n$ such that
$\theta^{-1}\phi=\tau_{\eta}$ and hence $\phi=\theta \tau_{\eta}$
as desired.
\qed

In other words,  $\theta_\xi$ is only determined modulo $\tau$-involutions.
We have also seen in  Lemma \ref{lem:6.2} that
$\theta_\xi\in \tau(\Z_2^n)$ if and only if $\xi\in C^\perp$.
Since $C^\perp \subset P$ by Lemma \ref{lem:6.9}, the association
$\xi+C^\perp \mapsto \theta_\xi\tau(\Z_2^n)$ defines a group
isomorphism between $P/C^\perp$ and $\pstab_V(F)/\tau(\Z_2^n)$.
Therefore, we have the following central extension:
$$
\begin{array}{ccccccccc}
  1 & \longrightarrow &\tau(\Z_2^n) & \longrightarrow & \pstab_V(F)
  & \longrightarrow & \pstab_V(F)/\tau(\Z_2^n) & \longrightarrow &  1
  \vsb\\
   && \downarrow\! \wr && ||  && \downarrow\! \wr &&
  \vsb\\
  1 &\longrightarrow& \Z_2^n/D^\perp &\longrightarrow& \pstab_V(F)
  &\longrightarrow&  P/C^\perp &\longrightarrow& 1
\end{array}
$$
The commutator relation in $\pstab_V(F)$ can also be described as follows.

\begin{theorem}\label{thm:6.11}
  For $\xi^1,\xi^2\in P$, let $\theta_{\xi^i}, i=1,2,$ be extensions
  of $\sigma_{\xi^i}$   to $\pstab_V(F)$.
  Then $[\theta_{\xi^1},\theta_{\xi^2}]=1$ if and only if
  $\la \alpha\cd \xi^1,\alpha\cd \xi^2\ra= 0$ for all $\alpha \in D$.
\end{theorem}

\pf
Since the case $\theta_{\xi^1}\in \theta_{\xi^2}\tau(\Z_2^n)$ is
trivial, we assume that $\sigma_{\xi^1}\ne \sigma_{\xi^2}$.
For $i=1,2$, set $C^{0,\xi^i}:=\{ \alpha \in C \mid \la \alpha,\xi^i\ra =0\}$
and $E:=C^{0,\xi^1}\cap C^{0,\xi^2}$.
Then $C^{0,\xi^i}$ are subcodes of $C$ with index 2 and $E$ is
a subcode of $C$ with index 4.
Let $\delta^1,\delta^2\in C$ such that
$C^{0,\xi^i}=E\sqcup (E+\delta^i)$.
By definition,
$\alpha\cd \xi^1$, $\alpha\cd \xi^2\in (E_\alpha)^\perp$
for all $\alpha \in D$ so that $E_\alpha$ contains a self-dual subcode
w.r.t.\ $\alpha$ if and only if $\la \alpha\cd \xi^1,\alpha\cd \xi^2\ra =0$.
We have seen that $\theta_{\xi^i}$ acts semisimply on each $V^\alpha$,
$\alpha\in D$, with two eigenvalues, and these eigenspaces are inequivalent
irreducible $V_{C^{0,\xi^i}}$-submodules.
For an irreducible $V_E$-submodule $W$ of $V^\alpha$, the subspace
$W+(V_{E+\delta^i}\cd W)$ forms a $V_{C^{0,\xi^i}}$-submodule so that
$\theta_{\xi^i}$ acts on $W$ by an eigenvalue.
Therefore, $\theta_{\xi^1}$ commutes with $\theta_{\xi^2}$ if and only if
$V^\alpha$ splits into a direct sum of 4 irreducible $V_E$-submodules for
all $\alpha \in D$.
Let $m_\alpha$ be the number of irreducible $V_E$-submodules of $V^\alpha$.
For $\alpha\in D$, let $H_\alpha$ and $H^0_\alpha$ be maximal self-orthogonal
subcodes of $C_\alpha$ and $E_\alpha$, respectively.
By the structure of irreducible modules over a code VOA shown in
Theorem \ref{thm:4.3}, $V^\alpha$ is a direct sum of $[C:H]$
irreducible $F$-submodules.
Moreover, again by Theorem \ref{thm:4.3} any irreducible $V_E$-submodule
of $V^\alpha$ is a direct sum of $[E:H^0_\alpha]$ irreducible $F$-submodules.
By counting the number of irreducible $F$-submodules of $V^\alpha$,
we have $m_\alpha \abs{H_\alpha}=4\abs{H^0_\alpha}$ as $[C:E]=4$.
Thus, $E_\alpha$ contains a self-dual subcodes w.r.t.\ $\alpha$
if and only if $m_\alpha=4$.
Hence, $\theta_{\xi^1}$ commutes with $\theta_{\xi^2}$ if and only if
$\la \alpha\cd \xi^1,\alpha\cd \xi^2\ra =0$.
\qed
\vsb

We have shown that the structure of $\pstab_V(F)$ is determined by
Theorems \ref{thm:6.7} and \ref{thm:6.11} only in terms of the
structure codes $(C,D)$.

\begin{remark}\label{rem:6.12}
  In \cite{Y2,Y3}, one of the authors has shown that for any
  Ising vector $e\in V^\natural$, we have no automorphism
  $g\in \aut(V^\natural)$ such that $g$ restricted on
  $(V^\natural)^{\la \tau_e\ra}$ is equal to $\sigma_e$.
  Thanks to Theorem \ref{thm:6.7}, we can give a simpler proof
  of this.
  As shown in \cite{DMZ}, the moonshine VOA $V^\natural$
  is framed.
  Take any Ising frame $F=\vir(e^1)\tensor \cds \tensor \vir(e^{48})$
  of $V^\natural$ and set $\xi=(10^{47})\in \Z_2^{48}$.
  Since $\xi$ is odd, there is no extension of
  $\sigma_\xi=\sigma_{e^1}$ to $\pstab_{V^\natural}(F)$
  by Theorem \ref{thm:6.7}.
  Since all the Ising vectors of $V^\natural$ are conjugate under
  $\aut(V^\natural)=\M$ (cf.\ \cite{C,M1}), $e$ and $e^1$ are
  conjugate.
  Therefore, there is no extension of $\sigma_e$ to $V^\natural$.
\end{remark}

At the end of this section, we give a brief description of the frame
stabilizer $\stab_V(F)$. Its structure is also discussed in
\cite{DGH}. It is clear that $\pstab_V(F)$ is a normal subgroup of
$\stab_V(F)$. Let $g\in \stab_V(F)$. Then $g$ induces a permutation
$\mu_g\in \mathrm{S}_n$ on the set of Ising vectors $\{
e^1,\dots,e^n\}$ of $F$, namely $g e^i = e^{\mu_g(i)}$. Since $g$
preserves the 1/16-word decomposition $V=\oplus_{\alpha\in D}
V^\alpha$, it follows that $g V^\alpha =V^{\mu_g (\alpha)}$ with
$\mu_g (\alpha)=(\alpha_{\mu_g(1)},\dots,\alpha_{\mu_g(n)})$. In
particular, $g$ restricted on $V^0$ defines an element of
$\aut(V^0)=\aut(V_C)$ which is a lift of $\aut(C)$. Therefore, every
element of $\stab_V(F)$ is a lift of $\aut(C)\cap \aut(D)$.
Conversely, we know that for any $\mu \in \aut(C)$, there exists
$\tilde{\mu}\in \aut(V_C)$ such that $\tilde{\mu}e^i= e^{\mu (i)}$
for $1\leq i\leq n$ by Theorem 3.3 of \cite{Sh}. It is shown in
Lemma 3.15 of \cite{SY} that if $\tilde{\mu}$ lifts to an element of
$\aut(V)$ then $\{ (V^\alpha)^{\tilde{\mu}} \mid \alpha \in D\}$
coincides  with $\{ V^\alpha \mid \alpha \in D\}$ as a set of
inequivalent irreducible $V_C$-modules. Therefore, there exists a
lift of $\tilde{\mu}\in \aut(V_C)=\aut(V^0)$ to an element of
$\aut(V)$ if and only if the subgroup $\{ V^\alpha \mid \alpha \in
D\}$ of the group formed by all the simple current $V_C$-module in
the fusion algebra is invariant under the conjugation action of
$\tilde{\mu}$. If such a lift of $\tilde{\mu}$ exists, it is unique
modulo $\pstab_V(F)$. For if $\tilde{\mu}$ and $\tilde{\mu}'$ are
two lifts of $\mu$, $\tilde{\mu}^{-1}\tilde{\mu}'$ fixes $F$
pointwise, showing $\tilde{\mu}^{-1}\tilde{\mu}'\in \pstab_V(F)$.
Thus, the factor group $\stab_V(F)/\pstab_V(F)$ is isomorphic to a
subgroup of $\aut(C)\cap \aut(D)$ which gives a slight refinement of
(3) of Theorem 2.8 in \cite{DGH}. The $V_C$-module structure of
$(V^\alpha)^{\tilde{\mu}}$ involves some extra information other
than $C$ and $D$, namely, if $V^\alpha \simeq
M_C(\alpha,\gamma;\iota_{\al})$ for $\al\in D$ then
$(V^\alpha)^{\tilde{\mu}}\simeq
M_C(\mu^{-1}\alpha,\gamma';\iota_{\mu^{-1}\al})$ for some codeword
$\gamma'\in \Z_2^n$, and this $\gamma'$ depends not only on $C$ and
$D$ but also on $\gamma$,  $\iota_{\al}$ and $\iota_{\mu^{-1}\al}$.
We do not have a general result for the lifting property of
$\aut(C)\cap \aut(D)$ at present.

\section{4A-twisted orbifold construction}

Let $V^\natural$ be the moonshine VOA constructed in \cite{FLM}. In
this section, we shall apply Theorem \ref{thm:6.7} to define a
4A-element of the Monster $\M=\aut(V^\natural)$ and exhibit that the
4A-twisted orbifold construction of the moonshine VOA $V^\natural$
will be $V^\natural$ itself.

By \cite{DGH,M3}, we can take an Ising frame
$F=\vir(e^1)\tensor \cds \tensor \vir(e^{48})$ of
$V^\natural$ such that the associated structure codes
$(\mathcal{C},\mathcal{D})$ are as follows:
$$
  \mathcal{C}=\mathcal{D}^\perp,\quad
  \mathcal{D}=\Span_{\Z_2}\{ (1^{16}0^{32}), (0^{32}1^{16}),
  (\alpha,\alpha,\alpha) \mid \alpha \in \RM(1,4)\},
$$
where $\RM(1,4)\subset \Z_2^{16}$ is the first order Reed-Muller code
defined by the generator matrix
$$
\begin{bmatrix}
  1111 & 1111 & 1111 & 1111
  \\
  1111 & 1111 & 0000 & 0000
  \\
  1111 & 0000 & 1111 & 0000
  \\
  1100 & 1100 & 1100 & 1100
  \\
  1010 & 1010 & 1010 & 1010
\end{bmatrix} .
$$
Note that
\begin{equation}\label{eq:7.1}
  \mathcal{C}=\{ (\alpha,\beta,\gamma) \in \Z_2^{48} \mid
  \alpha,\beta,\gamma\in \Z_2^{16}\ \text{are even and}\
  \alpha+\beta+\gamma \in \RM(2,4)\} .
\end{equation}

\begin{remark}\label{rem:7.1}
  The weight enumerator of $\RM(1,4)$ is $X^{16}+30X^8Y^8+Y^{16}$.
\end{remark}

Let $V^\natural=\oplus_{\alpha \in \mathcal{D}} (V^\natural)^\alpha$ be
the 1/16-word decomposition.
Set
$$
  \mathcal{P}:=\{ \xi \in \Z_2^{48} \mid \alpha\cd \xi \in \mathcal{C}\
  \text{for all}\ \alpha\in \mathcal{D}\}.
$$
Then for each $\xi\in\mathcal{P}$, one can define an automorphism
$\theta_{\xi}\in \pstab_{V^\natural}(F)$ such that
$\theta_{\xi}|_{V^0}=\sigma_\xi$ by  Theorem \ref{thm:6.7}.
Note also that
$\mathcal{D}=\mathcal{C}^\perp< \mathcal{P} < \mathcal{C}$
and $\sigma_{\xi^1}=\sigma_{\xi^2}$ if and only if
$\xi^1+\xi^2\in \mathcal{C}^\perp = \mathcal{D}$.

\begin{lemma}\label{lem:7.2}
  Let $\mathcal{C}$, $\mathcal{D}$ and $\mathcal{P}$ be defined as above.
  Then
  $$
    \mathcal{P}=\{ (\alpha,\beta,\gamma) \in \Z_2^{48} \mid
    \alpha,\beta,\gamma\in \RM(2,4) \ \text{and}\ \alpha+\beta+\gamma
    \in \RM(1,4)\},
  $$
  where $\RM(2,4)=\RM(1,4)^\perp\subset \Z_2^{16}$ is the second order
  Reed-Muller code of length 16.
\end{lemma}

\pf Set $E:=\{ (\alpha,\beta,\gamma) \mid \alpha,\beta,\gamma\in
\RM(2,4)\ \text{and}\ \alpha+\beta+\gamma\in \RM(1,4)\}$. We shall
first show that $E\subset \mathcal{P}$. It is clear that if $d_1\cd
\xi$, $d_2\cd\xi\in \mathcal{C}$ then $(d_1+d_2)\cd \xi = d_1\cd \xi
+ d_2\cd \xi \in \mathcal{C}$. Thus, we only need to show that $d\cd
\xi\in \mathcal{C}$ for $\xi\in E$ and  $d$ in a generating set of
$\mathcal{D}$.

Let $\xi=(\xi^1,\xi^2,\xi^3)\in E$ with each $\xi^i\in \Z_2^{16}$.
Then $\xi^1, \xi^2,\xi^3\in \RM(2,4)$.
Hence, $\alpha\cd \xi \in \mathcal{C}$ for
$\alpha=(1^{16}0^{32}), (0^{16}1^{16}0^{16})$ and
$(0^{32}1^{16})\in \mathcal{D}$.
Take any two codewords $\beta,\gamma\in \RM(1,4)$ with weight 8.
Then the weight of $\beta\cd \gamma$ is either 0, 4 or 8
by Remark \ref{rem:7.1}.
By the definition of $E$, we have $\xi^1+\xi^2+\xi^3\in \RM(1,4)$.
Therefore,
$$
  \la \beta,\gamma\cd (\xi^1+\xi^2+\xi^3)\ra
  \equiv \wt(\beta \cd \gamma \cd (\xi^1+\xi^2+\xi^3))\equiv 0\mod 2
$$
and hence $\gamma\cd (\xi^1+\xi^2+\xi^3)\in \RM(1,4)^\perp =\RM(2,4)$.
Since $\xi^i\in \RM(2,4)=\RM(1,4)^\perp$, all $\gamma\cd \xi^i$,
$i=1,2,3$, are even codewords.
Therefore,
$(\gamma,\gamma,\gamma)\cd \xi=(\gamma\cd \xi^1,\gamma\cd
\xi^2,\gamma\cd \xi^3) \in \mathcal{C}$ for all
$\gamma\in \RM(1,4)$.
As $\mathcal{D}$ is generated by the elements of the form:
$(1^{16}0^{32}), (0^{16}1^{16}0^{16}), (0^{32}1^{16})$ and
$(\gamma,\gamma,\gamma)$ with $\gamma\in \RM(1,4)$,
we have $E\subset \mathcal{P}$.

Conversely, assume $(\alpha^1,\alpha^2,\alpha^3) \in \mathcal{P}$
with $\alpha^i\in \Z_2^{16}$.
Then one has $(\alpha^1,\alpha^2,\alpha^3)\cd \beta \in \mathcal{C}$
for $\beta=(1^{16}0^{32})$, $(0^{16}1^{16}0^{16})$ and
$(0^{32}1^{16})\in \mathcal{D}$ so that $\alpha^i\in \RM(2,4)$
for $i=1,2,3$.
Moreover, for $(\gamma,\gamma,\gamma)\in \mathcal{D}$ with
$\gamma\in \RM(1,4)$, $(\alpha^1,\alpha^2,\alpha^3)\cd
(\gamma,\gamma,\gamma) \in \mathcal{C}$ is an even codeword.
Then it follows from \eqref{eq:7.1} that
$(\alpha^1+\alpha^2+\alpha^3)\cd \gamma\in \RM(2,4)$ and thus
$\alpha^1+\alpha^2+\alpha^3\in \RM(1,4)$. Hence $E=\mathcal{P}$.
\qed
\vsb

Take
\begin{equation}\label{eq:7.2}
  \xi=
  (1100 \, 0000 \, 1100 \, 0000 \, 0110 \, 0000\, 0110\, 0000 \, 1010 \, 0000 \, 1010 \, 0000)
  \in \Z_2^{48} .
\end{equation}
Then $\xi\in \mathcal{P}$ by Lemma \ref{lem:7.2}.
Set
$\mathcal{D}^0:=\{ \alpha \in \mathcal{D} \mid \wt(\alpha \cd \xi)\equiv 0 \mod 4\}$
and
$\mathcal{D}^1:=\{ \alpha \in \mathcal{D} \mid \wt(\alpha \cd \xi)\equiv 2 \mod 2\}$.
It is easy to see that
\[
  \mathcal{D}^0=\Span_{\Z_2}\left \{
  (1^{16}0^{32}),\ (0^{32}1^{16}),\
  (\alpha,\alpha,\alpha) \left |
\begin{matrix}
  \alpha = (1^{16}),\ (\{1^40^4\}^2),\\
   \ (\{1^20^2\}^4)\
  \text{or}\  (\{10\}^8)
\end{matrix}   \right. \right\}
\]
and $\mathcal{D}^1=(\{1^80^8\}^3)+\mathcal{D}^0$.
Therefore, the index $[\mathcal{D}:\mathcal{D}^0]$ is 2 and in this case
the involution $\sigma_\xi \in \aut((V^\natural)^0)$ can be extended to
an automorphism  $\theta_\xi \in \pstab_{V^\natural}(F)$ of order 4
by Theorem \ref{thm:6.7}.
We also set
$\mathcal{C}^0:=\{ \alpha \in \mathcal{C} \mid \la \alpha,\xi\ra =0\}$
and
$\mathcal{C}^1:=\{ \alpha \in \mathcal{C} \mid \la \alpha,\xi\ra =1\}$.
Let us consider a subgroup of $\mathcal{D}\times \C^*$ defined by
\begin{equation}\label{eq:7.3}
  \tilde{\mathcal{D}}
  := (\mathcal{D}^0\times \{ \pm 1\})\sqcup \l(\mathcal{D}^1\times
     \l\{\pm \sqrt{-1}\r\}\r).
\end{equation}
For $(\alpha,u) \in \mathcal{D}\times \C^*$, set
$(V^\natural)^{(\alpha, u)}:= \{ x\in (V^\natural)^\alpha \mid \theta_\xi x = u x\}$.
Then we have a $\tilde{\mathcal{D}}$-graded decomposition
\begin{equation}\label{eq:7.4}
  V^\natural=\bigoplus_{(\alpha,u)\in \tilde{\mathcal{D}}}
  (V^\natural)^{(\alpha,u)},\q
  (V^\natural)^{(0,1)}=V_{\mathcal{C}^0},\q
  (V^\natural)^{(0,-1)}=V_{\mathcal{C}^1},
\end{equation}
where $V_{\mathcal{C}^0}$ denotes the code VOA associated to
$\mathcal{C}^0$ and $V_{\mathcal{C}^1}$ is its module.
Since $(\mathcal{C}^0)_\alpha$ contains a self-dual subcode
w.r.t.\ $\alpha\in \mathcal{D}$ by Lemma \ref{lem:6.5},
all $(V^\natural)^{(\alpha,u)}$, $(\alpha,u)\in \tilde{\mathcal{D}}$,
are simple current $V_{\mathcal{C}^0}$-modules by
Corollary \ref{cor:4.16}.
Therefore, $V^\natural$ is a $\tilde{\mathcal{D}}$-graded simple current
extension of a code VOA $(V^\natural)^{(0,1)}=V_{\mathcal{C}^0}$.

By direct computation, it is not difficult to obtain
the following lemma.

\begin{lemma}\label{lem:7.3}
  For any non-zero $\alpha \in \mathcal{D}^0$, the subset
  $(\mathcal{C}^1)_\alpha$ is not empty.
  In other words, $[\mathcal{C}_\al, (\mathcal{C}^0)_\al]=2$.
\end{lemma}

\begin{remark}\label{rem:7.4}
  For a general framed VOA with structure codes $(C,D)$, it is
  possible that the set $(C^1)_\alpha$ is empty for some non-zero
  $\alpha \in D$.
  For example, we can take
  $D=\{ (0^{16}), (1^8\,0^8), (0^8\, 1^8), (1^{16})\}$,
  $C=D^\perp$, $\xi=(1^2 \, 0^{14})$ and $\al=(0^8\, 1^{8})$.
  Then, $C_\al=(C^0)_\al$ and $(C^1)_\alpha$ is empty.
  Note that $(C,D)$ can be realized as the structure codes of
  the lattice VOA $V_{E_8}$ (cf.\ \cite{DGH}).
\end{remark}

\begin{theorem}\label{thm:7.5}
  $\theta_\xi$ is a 4A-element of the Monster.
\end{theorem}

\pf
We shall compute the {\it McKay-Thompson series} of $\theta_\xi$
defined by
$$
  T_{\theta_\xi}(z):= \tr_{V^\natural}\, \theta_\xi\, q^{L(0)-1},\q
  q=e^{2\pii z}.
$$
Recall the notion of the conformal character of a module
$M=\oplus_{n\geq 0}M_{n+h}$ over a VOA $V$:
$$
  \ch_M(q):= \tr_M q^{L(0)-c/24}=\sum_{n=0}^\infty \dim_\C M_{n+h}q^{n+h-c/24}.
$$
It is clear that
$$
\begin{array}{ll}
  T_{\theta_\xi}(z)
  = &\dsum_{\alpha\in \mathcal{D}^0} \ch_{(V^\natural)^{(\alpha,1)}}(q)
  - \dsum_{\alpha\in \mathcal{D}^0} \ch_{(V^\natural)^{(\alpha,-1)}}(q)
  \vsb\\
  & +\sqrt{-1} \dsum_{\alpha\in \mathcal{D}^1} \ch_{(V^\natural)^{(\alpha,\sqrt{-1})}}(q)
  -\sqrt{-1} \dsum_{\alpha\in \mathcal{D}^1} \ch_{(V^\natural)^{(\alpha,-\sqrt{-1})}}(q).
\end{array}
$$
Let $\alpha \in \mathcal{D}^1$. Since $(V^\natural)^{(\alpha,\pm
\sqrt{-1})}$ are dual to each other, their conformal characters are
the same. Let $\alpha \in \mathcal{D}^0$ be a non-zero codeword. By
Lemma \ref{lem:7.3}, there exists a codeword in
$(\mathcal{C}^1)_\alpha$. Then by Corollary \ref{cor:4.7} one sees
that $(V^\natural)^{(\alpha,1)}$ and $(V^\natural)^{(\alpha,-1)}$
are isomorphic $F$-modules. Therefore, they have the same conformal
characters. Then
$$
\begin{array}{ll}
  T_{\theta_\xi}(z)
  &= \ch_{(V^\natural)^{(0,1)}}(q)-\ch_{(V^\natural)^{(0,-1)}}(q)
  \vsb\\
  &= \ch_{V_{\mathcal{C}^0}}(q)-\ch_{V_{\mathcal{C}^1}}(q)
  \vsb\\
  &= 2\ch_{V_{\mathcal{C}^0}}(q)-\ch_{V_{\mathcal{C}}}(q).
\end{array}
$$
The conformal character of a code VOA can be easily computed.
The following conformal characters are well-known (cf.\ \cite{FFR}):
$$
  \ch_{L(1/2,0)}(q)\pm  \ch_{L(1/2,1/2)}(q)
  = q^{-1/48}\PI_{n=0}^\infty (1\pm q^{n+1/2}).
$$
Since $\mathcal{C}=\mathcal{D}^\perp$ and
$\mathcal{C}^0=(\mathcal{D}+\la \xi\ra)^\perp$,
the weight enumerators of these codes are calculated
by the MacWilliams identity \cite{McS}:
$$
\begin{array}{lll}
  W_{\mathcal{C}}(x,y)
    &=&
    \dfr{1}{\abs{\mathcal{D}}} W_{\mathcal{D}}(x+y,x-y),
  \vsb\\
  W_{\mathcal{C}^0}(x,y)
    &=&
    \dfr{1}{\abs{\mathcal{D}+\la\xi\ra}} W_{\mathcal{D}+\la\xi\ra}(x+y,x-y),
\end{array}
$$
where
$$
\begin{array}{lll}
  W_{\mathcal{D}}(x,y)
  &=&
  x^{48} + 3 x^{32} y^{16} + 120 x^{24}y^{24} + 3x^{16}y^{32} + y^{48},
  \vsb\\
  W_{\mathcal{D}+\la\xi\ra}(x,y)
  &=&
  x^{48} + 2 x^{36} y^{12} + 3 x^{32} y^{16} + 30 x^{28} y^{20} +
    184 x^{24} y^{24} + 30 x^{20} y^{28}
  \vsb\\
  &&  + 3 x^{16} y^{32} + 2 x^{12}y^{36} + y^{48}.
\end{array}
$$
Now set
$$
\begin{array}{lll}
  f(x,y) &:=&  W_{\mathcal{D}+\la\xi\ra}(x,y) -  W_{\mathcal{D}}(x,y)
  \vsb\\
  &\ =& 2 x^{36} y^{12} + 30 x^{28} y^{20} + 64 x^{24} y^{24} +
      30 x^{20} y^{28} + 2 x^{12} y^{36}.
\end{array}
$$
Then one has
$$
\begin{array}{lll}
  T_{\theta_\xi}(z)
  &=&
  2\ch_{V_{\mathcal{C}^0}}(q)-\ch_{V_{\mathcal{C}}}(q)
  \vsb\\
  &=&
  \l[ 2W_{\mathcal{C}^0}(x,y) -
    W_{\mathcal{C}}(x,y)\r]_{x=\ch_{L(1/2,0)}(q),\ y=\ch_{L(1/2,1/2)}(q)}
  \vsb\\
  &=&  \dfr{1}{2^7}\l[ W_{\mathcal{D}+\la\xi\ra}(x+y,x-y)
     - W_{\mathcal{D}}(x+y,x-y)\r]_{ {x=\ch_{L(1/2,0)}(q)},\atop {\ y=\ch_{L(1/2,1/2)}(q)}}
  \vsb\\
  &=& \ds \dfr{1}{2^7}q^{-1} f\l(\PI_{n=0}^\infty (1+q^{n+1/2}),
    \PI_{n=0}^\infty (1-q^{n+1/2})\r)
  \vsb\\
  &=&
  q^{-1}+276q+2048q^2+\cds.
\end{array}
$$
Therefore, $\theta_\xi$ is a 4A-element of the Monster by \cite{ATLAS}.
\qed
\vsb

Next, we shall construct the irreducible 4A-twisted
$V^\natural$-module. Let us consider irreducible
$V_{\mathcal{C}^0}$-modules whose 1/16-word is $\xi$ defined in
\eqref{eq:7.2}.

\begin{lemma}\label{lem:7.6}
  (1) $(\mathcal{C}^0)_\xi$ is a self-dual subcode w.r.t.\ $\xi$.
  \\
  (2) For any $\gamma\in \Z_2^n$, the dual of
  $M_{\mathcal{C}^0}(\xi,\gamma;\iota)$
  is isomorphic to
  $M_{\mathcal{C}^0}(\xi,\gamma+\kappa;\iota)$
  with $\kappa=(\{10^7\}^2 0^{32})\in \Z_2^{48}$.
\end{lemma}

\pf By a direct computation, one can show that $\mathcal{C}_\xi
=(\mathcal{C}^0)_\xi$ is generated by the following generator
matrix:
$$
\begin{bmatrix}
  1100 0000 & 1100  0000 & 0000  0000 & 0000  0000 & 0000  0000 & 0000  0000
  \\
  0000  0000 & 0000  0000 & 0110  0000 & 0110  0000 & 0000  0000 & 0000  0000
  \\
  0000  0000 & 0000  0000 & 0000  0000 & 0000  0000 & 1010  0000 & 1010  0000
  \\
  1000  0000 & 1000  0000 & 0100  0000 & 0100  0000 & 0000  0000 & 0000  0000
  \\
  1000  0000 & 1000  0000 & 0000  0000 & 0000  0000 & 1000  0000 & 1000  0000
  \\
  1100  0000 & 0000  0000 & 0110  0000 & 0000  0000 & 1010  0000 & 0000  0000
\end{bmatrix}
$$
From this,  it is easy to see that $(\mathcal{C}^0)_\xi$ is a
self-dual code w.r.t.\ $\xi$. If we set
\begin{equation}\label{eq:7.5}
  \kappa=(1000\, 0000\, 1000\, 0000\, 0000\, 0000\, 0000\, 0000\, 0000\, 0000\, 0000\, 0000)
  \in \Z_2^{48},
\end{equation}
then $\supp(\kappa)\subset \supp(\xi)$ and we have $(-1)^{\la
\alpha,\kappa\ra}= (-1)^{\wt(\alpha)/2}$ for all $\alpha\in
(\mathcal{C}^0)_\xi$. Therefore, the dual of
$M_{\mathcal{C}^0}(\xi,\gamma;\iota)^*$ is isomorphic to
$M_{\mathcal{C}^0}(\xi,\gamma+\kappa;\iota)$ by Proposition
\ref{prop:4.9}. \qed \vsb

By the lemma above, we know that $(\mathcal{C}^0)_\xi$ is self-dual
w.r.t.\ $\xi$ and thus $(\mathcal{C}^0)_\xi$ equals to its own
radical. From now on, we shall fix a map $\iota :
(\mathcal{C}^0)_\xi\to \C^*$ such that the section map
$(\mathcal{C}^0)_\xi \ni \alpha \mapsto (\alpha,\iota(\alpha))\in
\pi_{\C^*}^{-1}(\mathcal{C}^0)$ is a group homomorphism. We  shall
simply denote $M_{\mathcal{C}^0}(\xi,\gamma;\iota)$ by
$M_{\mathcal{C}^0}(\xi,\gamma)$ and set
$W:=M_{\mathcal{C}^0}(\xi,0)$.

\begin{lemma}\label{lem:7.7}
  All $(V^\natural)^{(\alpha,u)} \fusion_{V_{\mathcal{C}^0}} W$,
  $(\alpha,u)\in \tilde{\mathcal{D}}$, are inequivalent irreducible
  $V_{\mathcal{C}^0}$-modules.
\end{lemma}

\pf
Suppose $(V^\natural)^{(\alpha,u)}\fusion_{V_{\mathcal{C}^0}} W \simeq
(V^\natural)^{(\beta,v)}\fusion_{V_{\mathcal{C}^0}} W$ with
$(\alpha,u), (\beta,v)\in \tilde{\mathcal{D}}$.
Since $((V^\natural)^{(\alpha,u)})^*\simeq (V^\natural)^{(\alpha,u)^{-1}}
= (V^\natural)^{(\alpha,u^{-1})}$, we have
$$
\begin{array}{ll}
  W
  &= (V^\natural)^{(\alpha,u^{-1})}\fusion_{V_{\mathcal{C}^0}} 
     (V^\natural)^{(\alpha,u)} \fusion_{V_{\mathcal{C}^0}} W
  \vsb\\
  &= (V^\natural)^{(\alpha,u^{-1})}\fusion_{V_{\mathcal{C}^0}} 
     (V^\natural)^{(\beta,v)} \fusion_{V_{\mathcal{C}^0}} W
  \vsb\\
  &= (V^\natural)^{(\alpha+\beta,u^{-1}v)}\fusion_{V_{\mathcal{C}^0}} W
\end{array}
$$
in the fusion algebra.
By considering 1/16-word decompositions, one has $\alpha=\beta$
and $u^{-1}v\in \{ \pm 1\}$.
Let $\delta\in \mathcal{C}$ be such that
$\mathcal{C}^1=\mathcal{C}^0+\delta$.
If $u^{-1}v=-1$, then
$(V^\natural)^{(\alpha+\beta,u^{-1}v)} =(V^\natural)^{(0,-1)}
=V_{\mathcal{C}^1} =M_{\mathcal{C}^0}(0,\delta)$
so that by Lemma \ref{lem:4.12} one has
$$
  W=
  (V^\natural)^{(\alpha+\beta,u^{-1}v)} \fusion_{V_{\mathcal{C}^0}} W
  = M_{\mathcal{C}^0}(0,\delta)\fusion_{V_{\mathcal{C}^0}} M_{\mathcal{C}^0}(\xi,0)
  = M_{\mathcal{C}^0}(\xi,\delta).
$$
It is shown in Lemma \ref{lem:7.6} that $(\mathcal{C}^0)_\xi$ contains
a self-dual subcode w.r.t.\ $\xi$.
Then $W$ is not isomorphic to $M_{\mathcal{C}^0}(\xi,\delta)$
by Lemma \ref{lem:4.6} and we obtain a contradiction.
Therefore, $u^{-1}v=1$ and hence $(\alpha,u)=(\beta,v)$.
This completes the proof.
\qed
\vsb

By Theorem \ref{thm:2.4} and the lemma above, the space
\begin{equation}\label{eq:7.6}
  V^\natural(\theta_\xi) := V^\natural \fusion_{V_{\mathcal{C}^0}} W
  = \bigoplus_{(\alpha,u)\in \tilde{\mathcal{D}}} (V^\natural)^{(\alpha,u)}
    \fusion_{(V^\natural)^{(0,1)}} M_{\mathcal{C}^0}(\xi,0)
\end{equation}
carries a unique structure of an irreducible $\chi_W$-twisted
$V^\natural$-module for some
$\chi_W\in (\tilde{\mathcal{D}})^*\subset \M$.
It is clear that $\chi_W\in \pstab_{V^\natural}(F)$.
Since the top weight of $W$ is 3/4,
$V^\natural(\theta_\xi)$ is neither 2A-twisted nor 2B-twisted
$V^\natural$-module so that $\chi_W$ is an element of $\M$
of order 4.
Hence, there exists $\xi'\in \mathcal{P}$ such that
$\chi_W^2|_{(V^\natural)^0}=\sigma_{\xi'}$.
By the construction of $V^\natural(\theta_\xi)$, we know that
$\mathcal{C}^0 \subset \{ \alpha \in \mathcal{C}\mid \la
\alpha,\xi'\ra=0\} \subsetneq \mathcal{C}$
so that $\sigma_\xi=\sigma_{\xi'}$ and
$\chi_W= \theta_\xi \cdot \tau_\eta$ for some $\eta\in \Z_2^{48}$.
Since the definition of $\theta_{\xi}$ is only unique up
to a product of $\tau$-involutions, we can take $\theta_\xi=\chi_W$.
Note that the definition of $V^\natural(\theta_\xi)$ in
\eqref{eq:7.6} depends only on the decomposition of $V^\natural$
as a $V_{\mathcal{C}^0}$-module.
Replacing $\theta_\xi$ by $\theta_\xi \tau_\beta$ with
$\beta\in \Z_2^{48}$ will only change the labeling of the
$V_{\mathcal{C}^0}$-modules in \eqref{eq:7.4} and does not
affect the isotypical decomposition of $V^\natural$ as a
$V_{\mathcal{C}^0}$-module.
Thus, we have constructed the irreducible $\theta_\xi$-twisted
$V^\natural$-module.

\begin{theorem}\label{thm:7.8}
  $V^\natural(\theta_\xi)$ defined in \eqref{eq:7.6} is an
  irreducible 4A-twisted module over $V^\natural$.
\end{theorem}

\begin{remark}\label{rem:7.9}
  By Lemma \ref{lem:7.6} and Corollary \ref{cor:4.7},
  we see that the top weight of   the 4A-twisted module
  is 3/4 and the dimension of the top level is 1.
\end{remark}

Let us consider the dual module $W^*$ of $W$.
It is clear that $W$ and $W^*$ are inequivalent since
$\kappa\not\in \mathcal{C}^0$.
All $(V^\natural)^{(\alpha,u)}\fusion_{V_{\mathcal{C}^0}} W^*$,
$(\alpha,u)\in \tilde{\mathcal{D}}$, are inequivalent
$V_{\mathcal{C}^0}$-modules, and so the space
\begin{equation}\label{eq:7.7}
  V^\natural(\theta_\xi^3)
  := V^\natural \fusion_{V_{\mathcal{C}^0}} W^*
  = \bigoplus_{(\alpha,u)\in \tilde{\mathcal{D}}} (V^\natural)^{(\alpha,u)}
    \fusion_{(V^\natural)^{(0,1)}} M_{\mathcal{C}^0}(\xi,\kappa)
\end{equation}
has a unique irreducible $\chi_{W^*}$-twisted
$V^\natural$-module structure with a linear character
$\chi_{W^*}\in \tilde{\mathcal{D}}^*$ (cf.\ Theorem \ref{thm:2.4}).
It is also shown \cite{DLM2} that the dual of $\chi_W$-twisted module
forms a $\chi_W^{-1}$-twisted module. The dual $V^\natural$-module
of $V^\natural(\theta_\xi)$ contains $W^*$ as a
$V_{\mathcal{C}^0}$-submodule, and $V^\natural(\theta_\xi^3)$ is
uniquely determined by $W^*$, so $\chi_{W^*}$ is actually equal to
$\theta_\xi^3=\theta_\xi^{-1}$ by our choice of $\theta_\xi$.
Therefore, $V^\natural(\theta_\xi^3)$ is the irreducible
$\theta_\xi^3$-twisted $V^\natural$-module.

In order to perform the 4A-twisted orbifold construction
of $V^\natural$, we classify the irreducible representations
of $(V^\natural)^{\la \theta_\xi\ra}$.
By \eqref{eq:7.4}, the fixed point subalgebra
$(V^\natural)^{\la \theta_\xi\ra}$ is a framed VOA with
structure codes $(\mathcal{C}^0,\mathcal{D}^0)$.

\begin{proposition}\label{prop:7.10}
  There are 16 inequivalent irreducible
  $(V^\natural)^{\la \theta_\xi\ra}$-modules.
  Every irreducible $(V^\natural)^{\la \theta_\xi\ra}$-module
  is a submodule of an irreducible $\theta_\xi^i$-twisted
  $V^\natural$-module for $0\leq i\leq 3$.
  Among them, 8 irreducible modules have integral top weights.
\end{proposition}

\pf
Since $V^\natural$ and $(V^\natural)^{\la \theta_\xi\ra}$ are
simple current extensions of the code VOA $V_{\mathcal{C}^0}$,
$V^\natural$ is a $\Z_4$-graded simple current extension of
$(V^\natural)^{\la \theta_\xi\ra}$ by Proposition \ref{prop:2.3}.
Then by Theorem \ref{thm:2.4},  every irreducible
$(V^\natural)^{\la \theta_\xi\ra}$-module is a submodule of a
$\theta_\xi^i$-twisted module for $0\leq i\leq 3$.
It is shown in \cite{DLM3} that $V^\natural$ has a unique
irreducible $\theta_\xi^i$-twisted module for $0\leq i\leq 3$ as
$V^\natural$ is holomorphic.
Therefore, we only have to show that each irreducible
$\theta_\xi^i$-twisted $V^\natural$-module decomposes into
a direct sum of 4 inequivalent irreducible
$(V^\natural)^{\la \theta_\xi\ra}$-submodules.
Since $V^\natural(\theta_\xi)$ and $V^\natural(\theta_\xi^3)$ are
dual to each other (cf.\ \cite{DLM2}), we know each of them
has four inequivalent irreducible
$(V^\natural)^{\la \theta_\xi\ra}$-submodules.
So we shall consider the $\theta_\xi^2$-twisted module.
Let $\kappa\in \Z_2^{48}$ be defined as in \eqref{eq:7.5}.
Then one can easily verify that $\theta_\xi^2=\tau_\kappa$.
We consider an irreducible $V_{\mathcal{C}^0}$-module
$M_{\mathcal{C}^0}(0,\kappa)$.
We claim that $(V^\natural)^{(\alpha,u)}\fusion_{V_{\mathcal{C}^0}}
M_{\mathcal{C}^0}(0,\kappa)$, $(\alpha,u)\in \tilde{\mathcal{D}}$,
are inequivalent irreducible $V_{\mathcal{C}^0}$-modules.
The irreducibility of
$(V^\natural)^{(\alpha,u)}\fusion_{V_{\mathcal{C}^0}} M_{\mathcal{C}^0}(0,\kappa)$
is clear since all $(V^\natural)^{(\alpha,u)}$ are simple
current $V_{\mathcal{C}^0}$-modules.
If
$(V^\natural)^{(\alpha,u)}\fusion_{V_{\mathcal{C}^0}}
M_{\mathcal{C}^0}(0,\kappa) \simeq (V^\natural)^{(\beta,v)}
\fusion_{V_{\mathcal{C}^0}} M_{\mathcal{C}^0}(0,\kappa)$,
then by the 1/16-word decomposition,
we have $\alpha=\beta$ and $u^{-1}v\in \{ \pm 1\}$.
Let $\delta\in \Z_2^{48}$ such that $\mathcal{C}^1=\mathcal{C}^0+\delta$.
If $u^{-1}v=-1$, then
$(V^\natural)^{(\alpha+\beta,u^{-1}v)}=M_{\mathcal{C}^0}(0,\delta)$
and one has
$$
\begin{array}{ll}
  M_{\mathcal{C}^0}(0,\kappa)
  & = (V^\natural)^{(\alpha,u^{-1})}\fusion_{V_{\mathcal{C}^0}} (V^\natural)^{(\beta,v)}
      \fusion_{V_{\mathcal{C}^0}} M_{\mathcal{C}^0}(0,\kappa)
  \vsb\\
  &= (V^\natural)^{(\alpha+\beta,u^{-1}v)} \fusion_{V_{\mathcal{C}^0}} 
     M_{\mathcal{C}^0}(0,\kappa)
  \vsb\\
  &= M_{\mathcal{C}^0}(0,\delta)\fusion_{V_{\mathcal{C}^0}} M_{\mathcal{C}^0}(0,\kappa)
  \vsb\\
  &= M_{\mathcal{C}^0}(0,\kappa+\delta).
\end{array}
$$
But this is a contradiction by Lemma \ref{lem:4.6}.
Therefore, all $(V^\natural)^{(\alpha,u)}\fusion_{V_{\mathcal{C}^0}}
M_{\mathcal{C}^0}(0,\kappa)$, $(\alpha,u)\in \tilde{\mathcal{D}}$,
are inequivalent irreducible $V_{\mathcal{C}^0}$-modules.
Now by Theorem \ref{thm:2.4} and Lemma \ref{lem:4.12} the space
\begin{equation}\label{eq:7.8}
  V^\natural(\theta_\xi^2)
  := \bigoplus_{(\alpha,u)\in\tilde{\mathcal{D}}} (V^\natural)^{(\alpha,u)}
     \fusion_{V_{\mathcal{C}^0}} M_{\mathcal{C}^0}(0,\kappa)
\end{equation}
forms an irreducible $\tau_\kappa=\theta_\xi^2$-twisted
$V^\natural$-module.
Thus $V^\natural(\theta_\xi^2)$ splits into four irreducible
$(V^\natural)^{\la \theta_\xi\ra}$-submodules as follows.
$$
\begin{array}{lll}
  V^\natural(\theta_\xi^2)
  =
  \ds \l\{\bigoplus_{\alpha\in \mathcal{D}^0} (V^\natural)^{(\alpha,1)}
    \fusion_{V_{\mathcal{C}^0}} M_{\mathcal{C}^0}(0,\kappa)\r\}
  \bigoplus \l\{\bigoplus_{\alpha\in \mathcal{D}^0} (V^\natural)^{(\alpha,-1)}
    \fusion_{V_{\mathcal{C}^0}} M_{\mathcal{C}^0}(0,\kappa)\r\}
  \vsb\\
  \ds
  \ \bigoplus \l\{\bigoplus_{\alpha\in \mathcal{D}^1} (V^\natural)^{(\alpha,\sqrt{-1})}
    \fusion_{V_{\mathcal{C}^0}}  M_{\mathcal{C}^0}(0,\kappa)\r\}
  \bigoplus \l\{\bigoplus_{\alpha\in \mathcal{D}^1} (V^\natural)^{(\alpha,-\sqrt{-1})}
    \fusion_{V_{\mathcal{C}^0}}  M_{\mathcal{C}^0}(0,\kappa)\r\} .
\end{array}
$$
Therefore, all irreducible $\theta_\xi^i$-twisted
$V^\natural$-modules are direct sums of 4 inequivalent
irreducible $(V^\natural)^{\la \theta_\xi\ra}$-submodules
and we have obtained 16 irreducible
$(V^\natural)^{\la \theta_\xi\ra}$-modules.
It remains to show that these 16 modules are inequivalent.
Since every irreducible
$(V^\natural)^{\la \theta_\xi\ra}$-module can be uniquely extended
to a $\theta_\xi^i$-twisted $V^\natural$-module
by Theorem \ref{thm:2.4}, these 16 irreducible modules are
actually inequivalent.

Among these 16 irreducible $(V^\natural)^{\la \theta_\xi\ra}$-modules,
we have 4 modules having integral top weights from $V^\natural$,
2 from $\theta_\xi^2$-twisted $V^\natural$-module, 1 from
$\theta_\xi$-twisted and 1 from $\theta_\xi^3$-twisted modules,
respectively.
This completes the proof.
\qed
\vsb

We have also shown that every irreducible
$\theta_\xi^i$-twisted $V^\natural$-module has a $\Z_4$-grading
which agrees with the action of $\theta_\xi$ on $V^\natural$.
By this fact, we adopt the following notation.

For $u\in \C^*$ satisfying $u^4=1$, we set
$V^\natural(1,u):=\{ a\in V^\natural \mid \theta_\xi a = u a\}$.
For $i=1$ or $3$, we define $V^\natural(\theta_\xi^i,1)$ to be
the unique irreducible $(V^\natural)^{\la\theta_\xi\ra}$-submodule
of $V^\natural(\theta_\xi^i)$ which has an integral top weight.
They can be defined explicitly as follows.
$$
\begin{array}{l}
  V^\natural(\theta_\xi,1)
  := \bigoplus_{\alpha\in \mathcal{D}^1} (V^\natural)^{(\alpha,-\sqrt{-1})}
     \fusion_{V_{\mathcal{C}^0}} M_{\mathcal{C}^0}(\xi,0),
  \vsv\\
  V^\natural(\theta_\xi^3,1)
  := \bigoplus_{\alpha\in \mathcal{D}^1} (V^\natural)^{(\alpha,\sqrt{-1})}
     \fusion_{V_{\mathcal{C}^0}} M_{\mathcal{C}^0}(\xi,\kappa).
\end{array}
$$
For $i=2$, there are two irreducible
$(V^\natural)^{\la \theta_\xi\ra}$-submodules
in $V^\natural(\theta_\xi^2)$ having integral top weights.
We shall define
\begin{equation}\label{eq:7.9}
\begin{array}{l}
  V^\natural(\theta_\xi^2,1)
  := \bigoplus_{\alpha\in \mathcal{D}^0} (V^\natural)^{(\alpha,-1)}
     \fusion_{V_{\mathcal{C}^0}} M_{\mathcal{C}^0}(0,\kappa),
  \vsv\\
  V^\natural(\theta_\xi^2,-1)
  := \bigoplus_{\alpha\in \mathcal{D}^0} (V^\natural)^{(\alpha,1)}
     \fusion_{V_{\mathcal{C}^0}} M_{\mathcal{C}^0}(0,\kappa).
\end{array}
\end{equation}
In addition, we define
\begin{equation}\label{eq:7.10}
  V^\natural(\theta_\xi^i,u)
  := V^\natural(1,u)\fusion_{(V^\natural)^{\la \theta_\xi\ra}} V^\natural(\theta_\xi^i,1)
  \q \text{for }\  1\leq i\leq 3.
\end{equation}
Set $G:= \la \theta_\xi  \ra \times \{ u\in \C^* \mid u^4=1\}$.
Then $\{V^\natural(g,u)|\,(g,u)\in G\}$ is the set of all inequivalent
irreducible $(V^\natural)^{\la \theta_\xi\ra}$-modules.

\begin{proposition}\label{prop:7.11}
  The fusion algebra associated to $(V^\natural)^{\la \theta_\xi\ra}$
  is isomorphic to the group algebra of $G$.
  The isomorphism is given by $V^\natural(g,u)\mapsto (g,u)$.
\end{proposition}

\pf
Since the structure codes of $(V^\natural)^{\la \theta_\xi\ra}$ is
$(\mathcal{C}^0,\mathcal{D}^0)$, $(V^\natural)^{\la \theta_\xi\ra}$
is a $\mathcal{D}^0$-graded simple current extension of $V_{\mathcal{C}^0}$.
So we have the following fusion rules:
$$
  V^\natural(1,u)\fusion_{(V^\natural)^{\la\theta_\xi\ra}} V^\natural(1,v)
  = V^\natural(1,uv) \q \text{for}\ u,v\in \C^*,\ u^4=v^4=1.
$$
Since all $V^\natural(g,u)$, $(g,u)\in G$, are
$\mathcal{D}^0$-stable, we can use Proposition \ref{prop:2.5}.
The following fusion rules of $V_{\mathcal{C}^0}$-modules
are already known:
$$
\begin{array}{ll}
  M_{\mathcal{C}^0}(0,\kappa)\fusion_{V_{\mathcal{C}^0}} M_{\mathcal{C}^0}(0,\kappa)
    = M_{\mathcal{C}^0}(0,0),
  & M_{\mathcal{C}^0}(0,\kappa)\fusion_{V_{\mathcal{C}^0}} M_{\mathcal{C}^0}(\xi,0)
    = M_{\mathcal{C}^0}(\xi,\kappa),
  \vsb\\
  M_{\mathcal{C}^0}(\xi,0)\fusion_{V_{\mathcal{C}^0}} M_{\mathcal{C}^0}(\xi,0)
    = M_{\mathcal{C}^0}(0,\kappa),
  & M_{\mathcal{C}^0}(\xi,0)\fusion_{V_{\mathcal{C}^0}} M_{\mathcal{C}^0}(\xi,\kappa)
    = M_{\mathcal{C}^0}(0,0).
\end{array}
$$
Therefore, we have the following fusion rules of
$(V^\natural)^{\la\theta_\xi\ra}$-modules:
$$
\begin{array}{ll}
  V^\natural(\theta^2,1)\fusion_{(V^\natural)^{\la\theta_\xi\ra}} V^\natural(\theta^2,1)
  = V^\natural(1,1),
  & V^\natural(\theta^2,1)\fusion_{(V^\natural)^{\la\theta_\xi\ra}} V^\natural(\theta,1)
  = V^\natural(\theta^3,1),
  \vsb\\
  V^\natural(\theta,1)\fusion_{(V^\natural)^{\la\theta_\xi\ra}} V^\natural(\theta,1)
  = V^\natural(\theta^2,1),
  & V^\natural(\theta,1)\fusion_{(V^\natural)^{\la\theta_\xi\ra}} V^\natural(\theta^3,1)
  = V^\natural(1,1).
\end{array}
$$
Since the fusion algebra is commutative and associative,
the remaining fusion rules are deduced from the above
and we can establish the isomorphism.
\qed
\vsb

A {\it $\theta_\xi$-twisted orbifold construction} of $V^\natural$
refers to a construction of a $\Z_4$-graded (simple current)
extension of the $\theta_\xi$-fixed point subalgebra
$(V^\natural)^{\la \theta_\xi\ra}$ by using the irreducible
submodules of $V^\natural(\theta_\xi^i)$ with integral weights.
By Proposition \ref{prop:7.10}, such modules are denoted by
$$
  V^\natural(1,\pm 1),\q
  V^\natural(1,\pm \sqrt{-1}),\q
  V^\natural(\theta_\xi^2,\pm 1),\q
  V^\natural(\theta_\xi,1)\q \text{and} \q
  V^\natural(\theta_\xi^3,1).
$$
By the fusion rules in Proposition \ref{prop:7.11}, there are
three possible extensions of $(V^\natural)^{\la \theta_\xi\ra}$,
namely,
\begin{equation}\label{eq:7.11}
\begin{array}{lll}
  V^\natural
  &=& V^\natural(1,1)\oplus V^\natural(1,-1) \oplus V^\natural(1,\sqrt{-1})
  \oplus V^\natural(1,-\sqrt{-1}),
  \vsb\\
  V_{2B}
  &=& V^\natural(1,1)\oplus V^\natural(1,-1) \oplus V^\natural(\theta_\xi^2,1)
  \oplus V^\natural(\theta_\xi^2,-1),
  \vsb\\
  V_{4A}
  &=& V^\natural(1,1) \oplus V^\natural(\theta_\xi,1)\oplus V^\natural(\theta_\xi^2,1)
  \oplus V(\theta_\xi^3,1).
\end{array}
\end{equation}
Consider the fixed point subalgebra
$(V^\natural)^{\la \theta_\xi^2\ra}= V^\natural(1,1)\oplus V^\natural(1,-1)$.
Since $\theta_\xi$ is a 4A-element of $\M$, its square $\theta_\xi^2$ belongs
to the 2B conjugacy class of $\M$ \cite{ATLAS}.
Therefore, by the original construction of the moonshine VOA in \cite{FLM},
the subalgebra $(V^\natural)^{\la \theta_\xi^2\ra}$
is isomorphic to the fixed point subalgebra $V_\Lambda^+$ of
the Leech lattice VOA $V_\Lambda$.
It is shown in \cite{D} that $V_\Lambda^+$ has four
inequivalent irreducible modules which are denoted by
$V_\Lambda^{\pm}$ and $V_\Lambda^{T\pm}$ in \cite{FLM}.
Since the top weights of irreducible $V_\Lambda^+$-modules
belong to $\Z/2$, the inequivalent irreducible
$(V^{\natural})^{\la \theta_\xi^2\ra}$-modules are given
by the list below:
\begin{equation}\label{eq:7.12}
\begin{array}{ll}
  V^\natural(1,1)\oplus V^\natural(1,-1),
  & V^\natural(1,\sqrt{-1})\oplus V^\natural(1,-\sqrt{-1}),
  \vsb\\
  V^\natural(\theta_\xi^2,1)\oplus V^\natural(\theta_\xi^2,-1),
  & V^\natural(\theta_\xi^2,\sqrt{-1}) \oplus V^\natural(\theta_\xi^2,-\sqrt{-1}).
\end{array}
\end{equation}
It is shown in \cite{H1} that there are two inequivalent simple
extensions of $V_\Lambda^+$; one is the moonshine VOA
$V^\natural=V_\Lambda^+\oplus V_\Lambda^{T+}$ and the other is
$V_\Lambda=V_\Lambda^+\oplus V_\Lambda^-$.
So we have the following isomorphisms:
\begin{equation}\label{eq:7.13}
\begin{array}{ll}
  V_\Lambda^+ \simeq V^\natural(1,1)\oplus V^\natural(1,-1),
  & V_\Lambda^{T+} \simeq V^\natural(1,\sqrt{-1})\oplus V^\natural(1,-\sqrt{-1}).
 \end{array}
\end{equation}
We shall prove that $\theta_\xi^2$-twisted orbifold construction $V_{2B}$
in \eqref{eq:7.11} is isomorphic to the Leech lattice VOA $V_\Lambda$.
For this, it is enough to show the following.

\begin{lemma}\label{lem:7.12}
  The top weight of $V^\natural(\theta_\xi^2,-1)$ is 1 and
  the dimension of the top level is 24.
\end{lemma}

\pf Recall the 1/16-word decomposition \eqref{eq:7.9} of
$V^\natural(\theta_\xi^2,-1)$ as a module over $V_{\mathcal{C}^0}$.
It contains a $V_{\mathcal{C}^0}$-submodule
$$
  (V^\natural)^{(0,1)}\fusion_{V_{\mathcal{C}}} M_{\mathcal{C}^0}(0,\kappa)
  = M_{\mathcal{C}^0}(0,\kappa) =V_{\mathcal{C}^0+\kappa}.
$$
By a straightforward computation we see that there are
exactly 24 weight two codewords in the coset $\mathcal{C}^0+\kappa$
and the support of each is one of the following:
\begin{equation}\label{11}
\begin{split}
  &\{1, 9\},\ \{2, 10\},\ \{3, 11\},\ \{4, 12\},\ \{5, 13\},\ 
  \{6, 14\},\ \{7, 15\},\ \{8, 16\},
  \vsb\\
  &\{17, 25\}, \{18, 26\}, \{19, 27\}, \{20, 28\}, \{21, 29\}, 
  \{22, 30\},\{23, 31\}, \{24, 32\},
  \vsb\\
 &\{33, 41\}, \{34, 42\}, \{35, 43\}, \{36, 44\}, \{37, 45\}, 
 \{38, 46\},\{39, 47\}, \{40, 48\} .
\end{split}
\end{equation}
Therefore, the top weight of $V^\natural(\theta^2,-1)$ is 1.
Then by the list of irreducible $V_\Lambda^+$-modules
in \eqref{eq:7.12},
the dimension of the top level of $V^\natural(\theta^2,-1)$ must be 24.
Or, one can directly check that all
$(V^\natural)^{(\alpha,1)}\fusion_{V_{\mathcal{C}^0}} M_{\mathcal{C}^0}(0,\kappa)$
has the top weight greater than 1 for any non-zero
$\alpha \in \mathcal{D}^0$ by considering their $F$-module structures.
\qed
\vsb

Since the top weights of $V_\Lambda^-$ and $V_\Lambda^{T-}$ are 1 and 3/2,
we have the remaining isomorphisms as follows.
\begin{equation}\label{eq:7.15}
\begin{array}{ll}
   V_\Lambda^- \simeq V^\natural(\theta_\xi^2,1)\oplus V^\natural(\theta_\xi^2,-1),
  & V_\Lambda^{T-} \simeq V^\natural(\theta_\xi^2,\sqrt{-1})
    \oplus V^\natural(\theta_\xi^2,-\sqrt{-1}).
\end{array}
\end{equation}
Therefore, the $\theta_\xi^2$-twisted orbifold $V_{2B}$ is
isomorphic to $V_\Lambda$.

\begin{remark}\label{rem:7.13}
  One can identify $V_{2B}$ with the Leech lattice VOA without
  the isomorphisms   in \eqref{eq:7.13} and \eqref{eq:7.15}.
  For, $V_{2B}$ can be defined as a framed VOA with structure codes
  $(\mathcal{C}^0\sqcup \mathcal{C}^1 \sqcup (\mathcal{C}^0+\kappa)\sqcup
  (\mathcal{C}^1+\kappa), \mathcal{D}^0)$ by Theorem \ref{thm:5.11},
  which is holomorphic by Corollary \ref{cor:5.10}.
  It follows from \eqref{11} that $V_{2B}$ contains a free bosonic
  VOA associated to a vector space of rank 24.
  Then $V_{2B}$ is isomorphic to a lattice VOA associated to an even
  unimodular lattice by \cite{LiX}.
  Since the weight one subspace of $V_{2B}$ is 24-dimensional,
  $V_{2B}$ is actually   isomorphic to the lattice VOA associated
  to the Leech lattice.
  Similarly, one can also show that
  $(V^\natural)^{\la \theta_\xi\ra}=V^\natural(1,1)$
  is isomorphic to a $\Z_2$-orbifold $V_L^+$ of a lattice VOA $V_L$
  for certain sublattice $L$ of $\Lambda$.
  For, we know that $V^\natural(1,1)\oplus V^\natural(\theta_\xi^2,-1)$
  forms a sub VOA of $V_{2B}$ isomorphic to a lattice VOA again
  by \cite{LiX}.
  Since $V^\natural(1,1)$ is a $\Z_2$-fixed point subalgebra of
  $V^\natural(1,1)\oplus V^\natural(\theta_\xi^2,-1)$ under an
  involution acting on $V^\natural(\theta_\xi^2,-1)$ by $-1$,
  we have the isomorphism as claimed.
  This isomorphism is first pointed out by Shimakura
  from a different view point.
\end{remark}

Next we consider the proper $\theta_\xi$-twisted orbifold
construction $V_{4A}$ in \eqref{eq:7.11}.
Let $\alpha \in \tilde{\mathcal{D}}^1$ be arbitrary.
We can find the following $V_{\mathcal{C}^0}$-submodule in $V_{4A}$.
$$
  U
  := V_{\mathcal{C}^0}\oplus V_{\mathcal{C}^1+\kappa}
  \oplus \{ (V^\natural)^{(\alpha,-\sqrt{-1})} \fusion_{V_{\mathcal{C}^0}}
  M_{\mathcal{C}^0}(\xi,0) \}
    \oplus \{ (V^\natural)^{(\alpha,\sqrt{-1})} \fusion_{V_{\mathcal{C}^0}}
  M_{\mathcal{C}^0}(\xi,\kappa) \} .
$$

\begin{lemma}\label{lem:7.14}
  There exists a unique structure of a framed VOA on $U$.
\end{lemma}

\pf
Set $U^0:= V_{\mathcal{C}^0}\oplus V_{\mathcal{C}^1+\kappa}$ and
$$
  U^1
  :=\{ (V^\natural)^{(\alpha,-\sqrt{-1})} \fusion_{V_{\mathcal{C}^0}}
  M_{\mathcal{C}^0}(\xi,0)\}
  \oplus \{ (V^\natural)^{(\alpha,\sqrt{-1})} \fusion_{V_{\mathcal{C}^0}}
  M_{\mathcal{C}^0}(\xi,\kappa)\} .
$$
Then $U^0$ is a code VOA associated to
$\mathcal{C}^0\sqcup (\mathcal{C}^1+\kappa)$
and $U^1$ is an irreducible $U^0$-module with the 1/16-word $\alpha+\xi$.
Clearly the top weight of $U^1$ is integral.
The dual code of $\mathcal{C}^0\sqcup (\mathcal{C}^1+\kappa)$
is given by $\mathcal{D}^0\sqcup (\mathcal{D}^1+\xi)$ and
it is straightforward to verify that $\mathcal{D}^0\sqcup
(\mathcal{D}^1+\xi)$ is triply even, i.e., $\wt(\alpha)$ is
divisible by 8 for any $\alpha \in \mathcal{D}^0\sqcup
(\mathcal{D}^1+\xi)$.
Therefore, $(\mathcal{C}^0\sqcup (\mathcal{C}^1+\kappa))_{\alpha+\xi}$
contains a doubly even self-dual subcode w.r.t.\ $\alpha+\xi$
by Remark \ref{rem:5.13}.
Then by Lemma \ref{lem:5.1}, $U=U^0\oplus U^1$
forms a framed VOA.
\qed
\vsb

By the lemma above, we can apply the extension property
of simple current extensions in Theorem \ref{thm:2.7} to
define a framed VOA structure on $V_{4A}$ with structure codes
$(\mathcal{C}^0\sqcup (\mathcal{C}^1+\kappa), \mathcal{D}^0\sqcup
(\mathcal{D}^1+\xi))$.
We know that $V_{4A}$ is holomorphic by Corollary \ref{cor:5.10}.
We shall prove that $V_{4A}$ is isomorphic to the moonshine VOA
$V^\natural$.
On $V_{2B}= V^\natural(1,1)\oplus V^\natural(1,-1) \oplus
V^\natural(\theta_\xi^2,1) \oplus V^\natural(\theta_\xi^2,-1)$,
define  $\psi_1, \psi_2\in \aut(V_{2B})$ by
$$
 \psi_1 :=
  \begin{cases}
    \ \ 1 \text{ on } V^\natural(1,1)\oplus V^\natural(1,-1),
    \vsb\\
    -1 \text{ on } V^\natural(\theta_\xi^2,1) \oplus V^\natural(\theta_\xi^2,-1),
  \end{cases}
$$
and
$$
\psi_2:=
  \begin{cases}
     \ \ 1 \text{ on } V^\natural(1,1)\oplus V^\natural(\theta_\xi^2,1),
     \vsb\\
     -1 \text{ on } V^\natural(1,-1) \oplus V^\natural(\theta_\xi^2,-1) .
  \end{cases}
$$
Then both of $\psi_1$, $\psi_2$ are involutions on $V_{2B}$
by the fusion rules in Proposition \ref{prop:7.11}.

\begin{lemma}\label{lem:7.15}
  The fixed point subalgebras $V_{2B}^{\la \psi_1\ra}$ and
  $V_{2B}^{\la \psi_2\ra}$ are isomorphic to the $\Z_2$-orbifold
  subalgebra $V_\Lambda^+$ of the Leech lattice VOA.
\end{lemma}

\pf
We have shown that $V_{2B}$ is isomorphic to the Leech lattice
VOA $V_\Lambda$.
By \eqref{eq:7.13}, the fixed point subalgebra
$V_{2B}^{\la \psi_1\ra}$ is isomorphic to the $\Z_2$-orbifold
$V_\Lambda^+$. So it remains to prove that $V_{2B}^{\la \psi_2\ra}$
is isomorphic to $V_{2B}^{\la \psi_1\ra}$.
Since the weight one subspace $(V_{2B})_1$ of $V_{2B}$ is a subspace
of $V^\natural(\theta^2,-1)$ by Lemma \ref{lem:7.12},
both $\psi_1$ and $\psi_2$ acts as $-1$ on $(V_{2B})_1$.
The weight one subspace of $V_{2B}$ generates a sub VOA
isomorphic to the free bosonic VOA $M_{\C \Lambda}(0)$ associated
to the linear space $\C \Lambda=\C\tensor_\Z \Lambda$.
Since $\psi_1\psi_2^{-1}$ trivially acts on the weight
one subspace of $V_{2B}$, $\psi_1\psi_2^{-1}$ commute with
the action of $M_{\C \Lambda}(0)$ on $V_{2B}$.
Therefore, $\psi_1\psi_2^{-1}$ is a linear character
$\rho_h=\exp(2\pii h_{(0)})\in \aut(V_{2B})$ induced
by a weight one vector $h \in (V_{2B})_1$.
Since $\psi_1=\psi_2=-1$ on the weight one subspace,
we have $\psi_i \rho_h=\rho_{-h}\psi_i=\rho_h^{-1} \psi_i$
for $i=0,1$.
Then
$\psi_1=\rho_h \psi_2= \rho_{h/2}\rho_{h/2}\psi_2
=\rho_{h/2}\psi_2 \rho_{h/2}^{-1}$
so that $\psi_1$ and $\psi_2$ are conjugate in $\aut(V_{2B})$.
From this we have the desired isomorphism
$\rho_{h/2}:(V_{2B})^{\la \psi_2\ra}\simto (V_{2B})^{\la \psi_1\ra}$.
\qed

\begin{corollary}\label{cor:7.16}
  There exists $\rho \in \aut((V^\natural)^{\la \theta_\xi\ra})$
  such that $V^\natural(1,-1)^\rho \simeq V^\natural(\theta_\xi^2,1)$.
\end{corollary}

\pf
Since $V^\natural(1,1)\subset \rho_{h/2} (V_{2B})^{\la \psi_2\ra}
\cap (V_{2B})^{\la \psi_2\ra} = (V_{2B})^{\la \psi_1\ra}\cap
(V_{2B})^{\la \psi_2\ra}$,
$\rho_{h/2}$ keeps $(V^\natural)^{\la\theta_\xi\ra}= V^\natural(1,1)$
invariant.
Thus the restriction of $\rho_{h/2}$ on $(V^\natural)^{\la\theta_\xi\ra}$
is the desired automorphism.
\qed
\vsb

By the classification of irreducible modules over $V_\Lambda^+$
and the fusion rules in Proposition \ref{prop:7.11},
the irreducible untwisted $(V_{2B})^{\la \psi_2\ra}$-modules
are as follows.
$$
\begin{array}{ll}
  V_\Lambda^+ \simeq V^\natural(1,1)\oplus V^\natural(\theta_\xi^2,1),
  & V_\Lambda^{T+} \simeq V^\natural(\theta_\xi,1)\oplus V^\natural(\theta_\xi^3,1),
  \vsb\\
  V_\Lambda^- \simeq V^\natural(1,-1)\oplus V^\natural(\theta_\xi^2,-1),
  & V_\Lambda^{T-} \simeq V^\natural(\theta_\xi,-1) \oplus V^\natural(\theta_\xi^3,-1).
\end{array}
$$
Actually, the isomorphisms above are induced by
$\rho\in \aut((V^\natural)^{\la\theta_\xi\ra})$ defined
in Corollary \ref{cor:7.16}.
By the isomorphisms above, the space
$$
  V_{4A}=V^\natural(1,1)\oplus V^\natural(\theta_\xi,1)
  \oplus V^\natural(\theta_\xi^2,1)
  \oplus V^\natural(\theta_\xi^3,1)
$$
is a $\Z_4$-graded simple current extension of
$(V^\natural)^{\la \theta_\xi\ra}$ and isomorphic to
$V^\natural=V_\Lambda^+\oplus V_\Lambda^{T+}$ as
a $(V^\natural)^{\la \theta_\xi\ra}$-module.
Since both $V_{4A}$ and $V^\natural$ are simple current
extensions of $(V^\natural)^{\la \theta_\xi\ra}$,
these two VOA structures are isomorphic.
Therefore, we have obtained our main result in this section.

\begin{theorem}\label{thm:7.17}
  The VOA $V_{4A}$ obtained by the 4A-twisted orbifold
  construction of $V^\natural$ is isomorphic to $V^\natural$.
\end{theorem}


\end{document}